%% file: I-rom-arxiv.tex
\documentclass[review]{elsarticle}

\usepackage{hyperref}

\usepackage{amsmath,color,url}

\definecolor{mygreen}{RGB}{28,172,0} %
\usepackage{graphicx,float}
\usepackage{tikz,amsfonts,bbm}
\usepackage{tabls}
\usepackage{multirow}

\usepackage{subcaption}
\journal{arXiv.org}

\begin{document}

\input{notation}

\begin{frontmatter}

\title{An Energy-Based Lengthscale for Reduced Order Models of Turbulent Flows}

\author[addresschanghong]{Changhong Mou}

\author[addresselia]{Elia Merzari}

\author[addressomer]{Omer San}

\author[addresstraian]{Traian Iliescu\corref{correspondingauthor}}
\cortext[correspondingauthor]{Corresponding author}
\ead{iliescu@vt.edu}
\ead[url]{sites.google.com/view/iliescu/}

\address[addresschanghong]{Department of Mathematics, University of Wisconsin-Madison, Madison, WI 53706, USA}
\address[addresselia]{Department of Nuclear Engineering, The Pennsylvania State University, University Park, PA 16802, USA}
\address[addressomer]{School of Mechanical and Aerospace Engineering, Oklahoma State University, Stillwater, OK 74078, USA}
\address[addresstraian]{Department of Mathematics, Virginia Tech, Blacksburg, VA 24061, USA}

\begin{abstract}
In this paper, we propose a novel reduced order model (ROM) lengthscale that is constructed by using energy distribution arguments.
The new energy-based ROM lengthscale is fundamentally different from the current ROM lengthscales, which are built by using dimensional arguments.
To assess the novel, energy-based ROM lengthscale, we compare it with a standard, dimensionality-based ROM lengthscale in two fundamentally different types of models:
(i) the mixing-length ROM (ML-ROM), which is a ROM closure model; and  
(ii) the evolve-filter-relax ROM (EFR-ROM), which is a regularized ROM.
We test the four combinations (i.e., ML-ROM and EFR-ROM equipped with the energy-based and dimensionality-based lengthscales) in the numerical simulation of the turbulent channel flow at $Re_{\tau} = 395$.
The numerical investigation 
yields the following conclusions:
(i) The new energy-based ROM lengthscale is significantly (almost two orders of magnitude) larger than the standard dimensionality-based ROM lengthscale.
As a result, the energy-based lengthscale yields more stable ML-ROMs and EFR-ROMs than the dimensionality-based lengthscale.
(ii) The 
energy-based 
lengthscale displays the correct asymptotic behavior with respect to the ROM dimension, whereas the 
dimensionality-based 
lengthscale does not.
(iii) 
The energy-based lengthscale yields ML-ROMs and (when significant filtering is effected) EFR-ROMs whose parameters are less sensitive (i.e., more robust) than the parameters of the ML-ROMs and EFR-ROMs based on the dimensionality-based lengthscale.
The novel energy-based lengthscale could enable the development of better scale-aware ROM strategies for flow-specific applications and 
is expected to have long term applications in nuclear reactor thermal-hydraulics.
\end{abstract}


\end{frontmatter}


\section{Introduction}

Reduced order models (ROMs) are models whose dimensions are orders of magnitude lower than the dimensions of full order models (FOMs), i.e., computational models constructed by using classical numerical discretizations (e.g., finite element or finite volume methods).
Because of their relatively low-dimensionality, ROMs can be used as efficient alternatives to FOMs in computationally intensive applications, e.g., flow control, shape optimization, and uncertainty quantification.

In recent years, ROMs have received a great deal of interest in nuclear engineering applications to generate models that can account for the fundamental physics of key phenomena while maintaining a low computational cost \cite{merzari2017, fick2018}. Examples are efforts to develop modal and machine-learning based 
ROMs for thermal stratification \cite{liusam}, which is recognized as critical for the licensing of liquid metal reactors. ROMs have also found applications in the modeling of parameterized coupled thermal-hydraulics and reactor physics problems \cite{vergari2020, vergari2021}. Finally,  another notable effort 
is the use of ROMs to develop flow acceleration for the advection-diffusion equation \cite{merzari2011}.

The Galerkin ROM (G-ROM) framework has been often used in the numerical simulation of fluid flows~\cite{HLB96,ahmed2021closures}.
The G-ROM is constructed as follows:
First, in an offline phase, the FOM is used to produce snapshots, which are then utilized to construct a low-dimensional (i.e., $r \ll N$) ROM basis $\{ \varphi_{1}, \ldots, \varphi_{r}\}$, where $r$ is the ROM dimension and $N$ is the FOM dimension.
There are several strategies for constructing the ROM basis, e.g., the proper orthogonal decomposition (POD)~\cite{HLB96,volkwein2013proper} and the reduced basis method~\cite{hesthaven2015certified,quarteroni2015reduced}.
In this paper, we exclusively use the POD to construct the ROM basis.
Next, the ROM basis is used together with a Galerkin projection to build the G-ROM, which 
can be written as follows:
\begin{eqnarray}
\overset{\bullet}{\ba}
	= {\bf F}(\ba), 	
	\label{eqn:g-rom}
\end{eqnarray}
where
$\ba$ is the vector of coefficients in the ROM approximation $\sum_{i=1}^{r} a_{i}(t) \bphi_{i}(\bx)$ of the variable of interest, 
$\overset{\bullet}{\ba}$ denotes the vector of time derivatives of $\ba$, 
and the vector $\bf F$ comprises the ROM operators that are preassembled in the offline phase. 
In the online phase, the G-ROM~\eqref{eqn:g-rom} is employed  for parameter 
values and/or time intervals that are different from those used in the training stage. 

We emphasize that, when FOMs (i.e., classical numerical methods, e.g., finite element, finite volume, or finite difference methods) are used for the spatial discretization of the fluid flow equations (e.g., the Navier-Stokes equations), the FOM lengthscale is generally defined as the meshsize, $h$, of the spatial mesh.
In contrast, the following {\it natural question} is, to our knowledge, {\it still open}:
\begin{center}
		\fbox{
		\addtolength{\linewidth}{3\fboxsep}
		\addtolength{\linewidth}{2\fboxrule}
		\begin{minipage}{0.42\linewidth}
			\vspace*{0.3cm}
			What is the ROM lengthscale, $\delta$?
			\vspace*{0.3cm}
		\end{minipage}
		}
\end{center}
To formulate the above question mathematically, we first assume that the following FOM and ROM variables are given, as is generally the case when ROMs are applied in practical settings:
\begin{itemize} \itemsep 5pt
    \item ROM variables, e.g., the ROM dimension, $r$, the total number of available ROM basis functions, $R$, the ROM basis functions, $\bphi_i$, and the corresponding eigenvalues, $\lambda_i$, in the eigenvalue problem solved to construct the POD basis.
    \item FOM variables, e.g., the FOM mesh size, $h$, the FOM solution, $\bu^{FOM}$, and the computational domain characteristic lengthscale, $L$.
\end{itemize}

Given these input FOM and ROM variables, we then try to answer the above question, i.e., to express the ROM lengthscale, $\delta$, as a function of the given FOM and ROM variables.

To motivate the need for a ROM lengthscale, we point out that the lengthscale is a fundamental notion in engineering, geophysical, and biomedical applications, where it is used to characterize the resolved spatial scales (i.e., to determine the size of the spatial scales approximated by the computational model) in flows around cars or airplanes, ocean or atmospheric flows, or blood flow in an artery, respectively.
To further motivate the need for a ROM lengthscale, we emphasize that there are ROMs that use a lengthscale in their very definition.
For example, in under-resolved simulations, i.e., when the number of ROM basis functions is not enough to accurately represent the turbulent flow dynamics, ROM closures and stabilizations are often used (see~\cite{ahmed2021closures} for a review).
There are several examples of ROM closure models that are defined by using a ROM lengthscale: the mixing-length ROM~\cite{AHLS88,HLB96,wang2012proper}, the Smagorinsky ROM~\cite{wang2012proper}, the dynamic SGS ROM~\cite{wang2012proper}, and the eddy viscosity variational multiscale ROM~\cite{wang2012proper}.
There are also several examples of ROM stabilizations that are defined by using a ROM lengthscale: the evolve-filter-relax ROM~\cite{wells2017evolve,gunzburger2019evolve,girfoglio2021pod, girfoglio2021pressure,strazzullo2022consistency} 
and the Leray ROM~\cite{wells2017evolve,sabetghadam2012alpha,kaneko2020towards,gunzburger2020leray,iliescu2018regularized}. 
A ROM lengthscale is needed to construct any of these two fundamentally different types of ROMs (i.e., closures and stabilizations).
A ROM lengthscale could also be useful in the preprocessing strategy advocated in~\cite{aradag2011filtered,farcas2022filtering} as a means to filter out the noise in the input data.

In this paper, we propose a {\it novel ROM lengthscale}, which is constructed by using energy distribution arguments.
To assess the new energy-based ROM lengthscale, we compare it with the classical ROM lengthscale used in~\cite{AHLS88,HLB96,wang2012proper}, which is based on fundamentally different, dimensional arguments.
To compare the two ROM lengthscales, we utilize them to build two different types of ROMs for under-resolved simulations:
(i) the mixing length ROM~\cite{AHLS88,HLB96,wang2012proper}, which is a ROM closure model that augments the standard G-ROM with a correction term; and
(ii) the evolve-filter-relax ROM~\cite{wells2017evolve}, which is a regularized ROM (Reg-ROM) that leverages ROM spatial filtering to increase the ROM stability and accuracy.
We test both the ROM closure and the Reg-ROM equipped with both ROM lengthscales in the numerical simulation of the turbulent channel flow at $Re_{\tau}=395$.
We note that a preliminary numerical investigation of the new ROM lengthscale was performed in~\cite{mou2022numerical}.

The rest of the paper is organized as follows:
In Section~\ref{sec:g-rom}, we outline the standard G-ROM, ML-ROM, and EFR-ROM.
In Section~\ref{sec:rom-lengthscale}, we define the new, energy-based lengthscale and 
the standard dimensionality-based lengthscale.
In Section~\ref{sec:numerical-results}, we present results for our investigation of the ML-ROM and EFR-ROM 
equipped with the two lengthscales in numerical simulation of the turbulent channel flow at $Re_{\tau}=395$.
Finally, in Section~\ref{sec:conclusions}, we draw conclusions and outline directions of future research.

\section{
{Reduced Order Models}}
    \label{sec:g-rom}

In this section, we outline the construction of the standard Galerkin ROM (G-ROM), mixing-length ROM (ML-ROM), and evolve-filter-relax ROM (EFR-ROM).    
As a mathematical model, we consider the incompressible Navier-Stokes equation (NSE): 
\begin{eqnarray}
\frac{\partial \bu}{\partial t} - Re^{-1} \Delta\bu +\bigl(\bu\cdot\nabla\bigr)\bu+\nabla p \,&=\, \bff,\label{eqn:nse-1} \\
\nabla\cdot\bu\,&=\,0, \label{eqn:nse-2}
\end{eqnarray}
where  $\bu = [u_1,u_2,u_3]^\top$ is the velocity vector field, $p$ the pressure field, $Re$ the Reynolds number, and $\bff$ the forcing term.
The NSE are equipped with appropriate boundary and initial conditions.

\subsection{Galerkin ROM (G-ROM)}

To build the G-ROM, we consider the centering trajectory of the flow,
\begin{eqnarray}
\bU(\bx) = \frac{1}{T}\int_t^{t+T} \bu(\bx,t)dt,
\end{eqnarray}
and 
assume that the ROM velocity approximation can be written as follows:
\begin{eqnarray}
\bu_r(\bx,t) \, 
=\,\bU(\bx)+\sum_{j=1}^ra_j(t)\bphi_j(\bx)\, ,\label{eqn:rom-soln}
\end{eqnarray}
where $\{\bphi_j\}_{j=1}^r$ are the ROM basis functions and $\ba = [a_1,\cdots,a_r]^\top$ are the sought ROM coefficients.
In our numerical experiments, we use the proper orthogonal decomposition (POD)~\cite{HLB96} to construct the ROM basis, but other ROM bases could be used~\cite{brunton2019data,hesthaven2015certified,quarteroni2015reduced}.
The next step in the G-ROM construction is to replace $\bu$ with $\bu_r$ in
~\eqref{eqn:nse-1} and project the resulting equations onto the space spanned by the ROM basis, $\{\bphi_j\}_{j=1}^r$.
This yields the G-ROM:
\begin{eqnarray}
\left(\frac{\partial \bu_r}{\partial t},\bphi_i \right)
+\bigl((\bu_r\cdot\nabla)\bu_r,\bphi_i\bigr)
+ Re^{-1} \bigl(\nabla\bu_r,\nabla\bphi_i\bigr) = \bigl(\bff,\bphi_i\bigr)\,,
\ \ i =1,\cdots,r\, , \label{eqn:g-rom-1}
\end{eqnarray}
where $( \cdot , \cdot )$ denotes the $L^2$ inner product.
The G-ROM can be 
written as the following dynamical system for the vector of time coefficients, $\ba(t)$:
\begin{eqnarray}
	\overset{\bullet}{\ba} 
	= \bb+A  \ba 
	+ \ba^\top B  \ba \, , 
	\label{eqn:g-rom-U}
\end{eqnarray}
where
\begin{align}
&\bb_i \,=\, \bigl(\bphi_i,\bff\bigr)-\bigl(\bphi_i,\bU\cdot\nabla\bU\bigr)
- Re^{-1} \bigl(\nabla\bphi_i,\nabla\bU\bigr)\, ,\\[0.3cm]
&\bA_{im} \,=\, -\bigl(\bphi_i,\bU\cdot\nabla\bphi\bigr) - \bigl(\bphi_i,\bphi_m\cdot\nabla\bU\bigr)
- Re^{-1} \bigl(\nabla\bphi_i,\nabla\bphi_m\bigr)\, ,\\[0.3cm]
&\bB_{imn} \, =\, -\bigl(\bphi_i,\bphi_m\cdot\nabla\bphi_n\bigr)\, .
\end{align}


\subsection{Mixing-Length ROM (ML-ROM)}
    \label{sec:ml-rom}

The G-ROM~\eqref{eqn:g-rom-U} is computationally efficient and relatively accurate in the numerical simulation of laminar flows.
However, the G-ROM generally yields inaccurate results in the numerical simulation of turbulent flows.
The main reason for the G-ROM's inaccuracy is that it is used in the under-resolved regime, i.e., when the number of ROM basis functions, $r$, is not large enough to accurately represent the complex dynamics of the turbulent flow.
Thus, in those cases, the G-ROM 
can be equipped with a ROM closure model, which models the effect of the discarded ROM modes $\{ \bphi_{r+1}, \ldots \}$ on the G-ROM dynamics.
In general, the G-ROM with a closure model can be written as 
\begin{eqnarray}
	\overset{\bullet}{\ba} 
	= \bb
	+ A  \ba 
	+ \ba^\top B  \ba 
	+ \btau \, , 
	\label{eqn:g-rom-closed}
\end{eqnarray}
where $\btau$ is the ROM closure model.
The current ROM closure models are carefully surveyed in~\cite{ahmed2021closures}.
Some of these ROM closure models are inspired from classical large eddy simulation (LES) closure modeling~\cite{BIL05,sagaut2006large}.
To construct these LES-ROM closure models, one needs to define a ROM lengthscale, which represents the size of the spatial scales modeled in the LES-ROM.
There are only a few ROM lengthscales in current use.
In Section~\ref{sec:rom-lengthscale}, we define a novel ROM lengthscale.
To assess this new ROM lengthscale, we consider one of the simplest ROM closure models, the ML-ROM~\cite{HLB96,wang2012proper}, in which the ROM closure term $\btau$ in~\eqref{eqn:g-rom-closed} is written as
\begin{eqnarray}
	\btau
	= - \bigl( \alpha \, U_{ML} \, \delta \bigr) \, S_r \, \ba \, ,
	\label{eqn:ml-v2}
\end{eqnarray}
where $\delta$ is one of the two ROM lengthscales defined in Section~\ref{sec:rom-lengthscale}, $U_{ML}$ is a characteristic velocity scale, $\alpha$ is a constant, and $S_r$ is the ROM stiffness matrix with entries $(S_r)_{ij} = \bigl(\nabla\bphi_i,\nabla\bphi_j\bigr), \ i, j = 1, \ldots, r$.
The ML-ROM model~\eqref{eqn:ml-v2} is a functional closure model, which aims at increasing the ROM viscosity in order to dissipate energy and mimic the effect of the discarded modes~\cite{CSB03}.  
The ML-ROM~\eqref{eqn:ml-v2} was first used in~\cite{AHLS88,HLB96} and was further investigated in~\cite{wang2012proper}.

\subsection{Evolve-Filter-Relax ROM (EFR-ROM)}
    \label{sec:efr-rom}

Regularized ROMs (Reg-ROMs)~\cite{kaneko2020towards,wells2017evolve} represent an alternative to ROM closures (e.g., the ML-ROM outlined in Section~\ref{sec:ml-rom}) in under-resolved simulations of turbulent flows.
Instead of adding a closure term, $\btau$, as in ROM closure modeling (see~\eqref{eqn:g-rom-closed}), Reg-ROMs are constructed by using ROM spatial filtering of various terms in the NSE to increase the ROM numerical stability.
Although regularized models have been used for decades in classical CFD~\cite{fischer2001filter,layton2012approximate,mullen1999filtering}, Reg-ROMs have only been recently developed~\cite{wells2017evolve}. 

The evolve-filter-relax ROM (EFR-ROM) is one of the most popular Reg-ROMs.
EFR-ROM is a modular ROM stabilization strategy that consists of three steps:
In the first step, which is called the evolve step, the standard G-ROM is used to advance the current EFR-ROM time iteration, $\ba^{n}$, to an intermediate approximation, $\bw^{n+1}$.
In the second step, which is called the filter step, the intermediate approximation, $\bw^{n+1}$, is filtered with the ROM differential filter~\cite{wells2017evolve}, which yields the filtered intermediate approximation, $\overline{\bw}^{n+1}$. 
In the third step, which is called the relax step, the EFR-ROM approximation at the next time step, $\ba^{n+1}$, is calculated as the convex combination of the intermediate approximation, $\bw^{n+1}$, and the filtered intermediate approximation, $\overline{\bw}^{n+1}$.
The EFR-ROM is summarized in the following algorithm:
\begin{eqnarray}
         &	\text{\bf (I)}& \text{\emph{ Evolve}:} \qquad 
        \ba^{n}
        \quad \stackrel{G-ROM~\eqref{eqn:g-rom-U}}{\longmapsto} \quad 
        \bw^{n+1}
            \nonumber \\[0.3cm]
        &	\text{\bf (II)} &\text{\emph{ Filter:}} \qquad
        \left( M_r + \gamma \, \delta^2 \, S_r \right) \overline{\bw}^{n+1}
        = \bw^{n+1}
        \nonumber \\[0.3cm]
        &	\text{\bf (III)} &\text{\emph{Relax:}} \qquad 
        \ba^{n+1}
        = (1 - \chi) \, \bw^{n+1}
        + \chi \, \overline{\bw}^{n+1} \, ,
        \nonumber
\end{eqnarray}
where $\chi \in [0,1]$ is a relaxation parameter. 
In Step (II), we use the ROM differential filter (DF) with an explicit ROM lengthscale, $\delta$, which represents the filtering radius. 
The DF acts as a spatial filter by eliminating the small scales (i.e., high frequencies) from the input data~\cite{BIL05}. 
We note that, in Step (II),  we modify the classical DF~\cite{BIL05} by introducing a new parameter $\gamma$.
This new parameter $\gamma$ has a role similar to that of the parameter $\alpha$ used in the ML-ROM~\eqref{eqn:ml-v2}:
It controls the amount of filtering used in the DF.
Step (III) is a relaxation step in which the EFR-ROM velocity approximation at the new time step is defined as a convex combination of the approximations obtained in Step (I) and Step (II).
The relaxation parameter $\chi$ diminishes the magnitude of the numerical diffusion \cite{ervin2012numerical, fischer2001filter,mullen1999filtering} and increases the accuracy~\cite{bertagna2016deconvolution,ervin2012numerical}.
The scaling $\chi \sim \Delta t$, where $\Delta t$ is the time step size, is a popular choice \cite{ervin2012numerical}, but higher values have also been used (see, e.g., \cite{bertagna2016deconvolution}).
To our knowledge, the EFR-ROM was first used  in~\cite{wells2017evolve} without the relaxation step, and in~\cite{gunzburger2019evolve} with the relaxation step (see also~\cite{strazzullo2022consistency} and references therein).
EFR-ROM was also investigated in~\cite{wells2017evolve,gunzburger2019evolve,girfoglio2021pod, girfoglio2021pressure,strazzullo2022consistency}.

\section{ROM Lengthscales}
    \label{sec:rom-lengthscale}

In this section, we present two different ROM lengthscales: 
In Section~\ref{sec:rom-lengthscale-dimensional}, we present a standard ROM lengthscale, denoted $\delta_1$, which is constructed by using  dimensional analysis arguments.
In Section~\ref{sec:rom-lengthscale-energy}, we propose a novel ROM lengthscale, denoted $\delta_2$, which is constructed by using energy balance arguments.
As explained in the introduction, both definitions aim at expressing the ROM lengthscale as a function of the following two types of input variables: 
(i) ROM variables (e.g., the ROM dimension, $r$, the total number of ROM basis functions, $R$, the eigenvalues, $\lambda_i$, and the ROM basis functions, $\bphi_i$).
(ii) FOM variables (e.g., the fine FOM mesh size, $h$, the FOM solution, $\bu^{FOM}$, and the computational domain characteristic lengthscale, $L$).
Given these input FOM and ROM variables, we then try to answer the following natural question:
{\it What is the ROM lengthscale, $\delta$?}

\subsection{Standard dimensionality-based ROM Lengthscale $\delta_1$} 
    \label{sec:rom-lengthscale-dimensional}

In this section, we use dimensional analysis to construct the first ROM lengthscale, $\delta_1$.
To this end, we follow the approach used in Section 3.2 in \cite{wang2012proper}, which, in turn, is based on the pioneering 
ML-ROM proposed in~\cite{AHLS88} for a turbulent pipe flow.

To construct the ROM lengthscale $\delta_1$, we first define the componentwise FOM velocity fluctuations,
i.e., the unresolved component of the velocity, which is computed by using FOM data:
\begin{eqnarray}
    {u^{'}_{i}}^{FOM} 
    = \sum_{j=r+1}^{R} a_j^{FOM} \, \varphi_j^i ,
    \qquad
    i = 1, 2, 3,
    \label{eqn:u-prime}
\end{eqnarray}
where $R$ is the total number of ROM modes, 
$\varphi_j^i$ are the componentwise ROM basis functions, and $a_j^{FOM} = (\bu^{FOM}, \bphi_j)$ are the ROM coefficients computed by using FOM data.
Using the componentwise FOM velocity fluctuations ${u^{'}_{1}}^{FOM}, {u^{'}_{2}}^{FOM}$, and ${u^{'}_{3}}^{FOM}$ in the $x, y$, and $z$ directions, respectively, we 
build the FOM velocity fluctuation vector field ${\bu'}^{FOM} =  [ {u^{'}_{1}}^{FOM}, {u^{'}_{2}}^{FOM}, {u^{'}_{3}}^{FOM} ]$.
Since ${\bu'}^{FOM}$ varies with time, we calculate the time averaged value of ${\bu'}^{FOM}$, i.e.,
\begin{eqnarray}
    \langle {\bu'}^{FOM} \rangle_t (\bx) 
    = \frac{1}{M} \sum_{k=1}^{M} {\bu'}^{FOM}(\bx, t_k) 
    = \frac{1}{M} \sum_{k=1}^{M} \sum_{l= r+1}^{R} \biggl(\bu^{FOM}(\cdot,t_k),\bphi_l(\cdot)\biggr)\bphi_l(\bx),
\end{eqnarray}
where $M$ is the number of snapshots. 

To construct the ROM lengthscale $\delta_1$, we adapt equation (22) in~\cite{wang2012proper} to our computational setting (i.e., the turbulent channel flow in Section~\ref{sec:numerical-results}):
\begin{eqnarray}
\delta_{1}
:= \left( \frac{\int_{0}^{L_1} \, \int_{0}^{L_2} \, \int_{0}^{L_3} 
\sum_{i=1}^{3} {u^{'}_{i}}^{FOM} {u^{'}_{i}}^{FOM}  
\, dx_1 \, dx_2 \, dx_3}
{\int_{0}^{L_1} \, \int_{0}^{L_2} \, \int_{0}^{L_3} 
\sum_{i=1}^{3} \sum_{j=1}^{3}  \frac{\partial {u^{'}_{i}}^{FOM}}{\partial x_j\hfill} \, \frac{\partial {u^{'}_{i}}^{FOM}}{\partial x_j\hfill}  
\, dx_1 \, dx_2 \, dx_3} \right)^{1/2} ,
\label{eqn:delta-1}
\end{eqnarray}
where 
$L_1, L_2$, and $L_3$ are the streamwise, wall-normal, and spanwise dimensions of the 
computational domain of the turbulent channel flow test problem, 
respectively.

Note that a quick dimensional analysis shows that the quantity defined in \eqref{eqn:delta-1}
has the units of a lengthscale:
\begin{eqnarray}
[ \delta_{1} ] 
= \left( 
\frac{\frac{m}{s} \, \frac{m}{s} \, m^3}
{\frac{1}{s} \, \frac{1}{s} \, m^3 }
\right)^{1/2}
= m \, .
\label{eqn:delta-1-dimensional-analysis}
\end{eqnarray}
We note that the ROM lengthscale, $\delta_1$, defined in~\eqref{eqn:delta-1}, depends on the FOM velocity fluctuation vector field, ${\bu'}^{FOM}$.

We also note that an alternative lengthscale was defined in equation (23) in~\cite{wang2012proper}.
Since this alternative lengthscale was not used in the numerical investigation in~\cite{wang2012proper} (because it was harder to implement), we do not consider it in this study.

\subsection{Novel Energy-Based ROM Lengthscale $\delta_2$} 
    \label{sec:rom-lengthscale-energy}

In this section, we use energy balancing arguments and propose a new ROM lengthscale, $\delta_2$.
Noticing that the ROM truncation level, $r$, has the role of dividing the kinetic energy of the system, we can require that the new ROM lengthscale, $\delta_2$, do the same.
Specifically, we require that the ratio of kinetic energy contained in the first $r$ ROM modes,  $\sum_{i=1}^{r} \lambda_i$, to the kinetic energy contained in the total number of ROM modes, $\sum_{i=1}^{R} \lambda_i$, 
be equal to the ratio of the kinetic energy that can be represented on an imaginary mesh of size $\delta_2$, $KE(\delta_2)$, to the kinetic energy that can be represented on the FOM mesh, $KE(h)$:
\begin{eqnarray}
\frac{\sum_{i=1}^{r} \lambda_i}{\sum_{i=1}^{R} \lambda_i} 
= \frac{KE(\delta_2)}{KE(h)} \, .
\label{energy_balancing}
\end{eqnarray}
\begin{remark}
We emphasize that the mesh of size $\delta_2$ is not used in the actual ROM construction.
Instead, this imaginary mesh is used to highlight the fundamental difference between the physical lengthscales of the space generated by the first $r$ ROM modes, $\bX^{r} = \text{span}\{ \bphi_{1}, \ldots, \bphi_{r} \}$, and the physical lengthscales of the space generated by all $R$ ROM modes, $\bX^{R} = \text{span}\{ \bphi_{1}, \ldots, \bphi_{r}, \bphi_{r+1}, \ldots, \bphi_{R} \}$.
Indeed, since the FOM data was used to build all the ROM modes, all the lengthscales of the space $\bX^{R}$ can be represented on the FOM mesh.
However, the imaginary mesh of size $\delta_2$ can represent only the lengthscales of the space $\bX^{r}$; it cannot be expected to represent the physical lengthscales of the space spanned by the higher index ROM modes $\{ \bphi_{r+1}, \ldots, \bphi_{R} \}$, which are generally associated with the small scales.
\end{remark}
To compute the ratio $\frac{KE(\delta_2)}{KE(h)}$ in \eqref{energy_balancing}, we transfer the problem to the usual Fourier space.
To this end, we first notice that $\delta_2$ defines a cutoff wavenumber:
\begin{eqnarray}
k_{\delta_2} := \frac{2 \, \pi}{\delta_2} \, .
\label{def_k_c}
\end{eqnarray}
We then notice that the kinetic energy in the system can be written in terms of the energy 
spectrum, $E(\cdot)$:
\begin{eqnarray}
KE(k) 
= \int_{k_0}^{k} E(k') \, dk' \, ,
\label{def_E}
\end{eqnarray}
where $\displaystyle k_0 = \frac{2 \, \pi}{L}$ is the Fourier wavenumber that corresponds to the computational domain characteristic lengthscale, $L$.
In the case of isotropic, homogeneous turbulence, we have the usual energy spectrum given by
Kolmogorov's theory~\cite{sagaut2006large, Pop00}
\begin{eqnarray}
E(k) 
\sim C \, \varepsilon^{2/3} \, k^{-5/3} \, .
\label{energy spectrum}
\end{eqnarray}
Thus, the condition imposed in \eqref{energy_balancing} can be written as
\begin{eqnarray}
\frac{\int_{k_0}^{k_{\delta_2}} E(k) \, dk}{\int_{k_0}^{k_h} E(k) \, dk}
= \frac{\sum_{i=1}^{r} \lambda_i}{\sum_{i=1}^{R} \lambda_i} 
\stackrel{\text{notation}}{=} \Lambda \, ,
\label{delta_energy}
\end{eqnarray}
where $\displaystyle k_h = \frac{2 \, \pi}{h}$ is the highest Fourier wavenumber that can be resolved on the given FOM meshsize, $h$.
The LHS of \eqref{delta_energy} can be evaluated by using \eqref{energy spectrum}:
\begin{eqnarray}
\int_{k_0}^{k_{\delta_2}} E(k) \, dk
= C \, \varepsilon^{2/3} \, \int_{k_0}^{k_{\delta_2}} k^{-5/3}\, dk
= C \, \varepsilon^{2/3} \, \frac{k_{\delta_2}^{-2/3} - k_{0}^{-2/3}}{-2/3} \, ,
\label{delta_energy_1}
\end{eqnarray}
and, similarly, 
\begin{eqnarray}
\int_{k_0}^{k_h} E(k) \, dk
= C \, \varepsilon^{2/3} \, \int_{k_0}^{k_h} k^{-5/3}\, dk
= C \, \varepsilon^{2/3} \, \frac{k_h^{-2/3} - k_{0}^{-2/3}}{-2/3} .
\label{delta_energy_2}
\end{eqnarray}
Plugging \eqref{delta_energy_1} and \eqref{delta_energy_2} back into \eqref{delta_energy}, simplifying, and rearranging,
we obtain
\begin{eqnarray}
k_{\delta_2}^{-2/3}
= \Lambda \, k_{h}^{-2/3}
+ \left( 1 - \Lambda \right) \, k_{0}^{-2/3} \, .
\label{delta_energy_3}
\end{eqnarray}
Since $1 \leq r \leq R$, $\Lambda$ satisfies the inequality $0 < \Lambda \leq 1$.
Thus,~\eqref{delta_energy_3} 
implies that $k_{\delta_2}^{-2/3}$ is a convex combination of $k_{h}^{-2/3}$ and $k_{0}^{-2/3}$.
Furthermore, as expected, $k_{\delta_2}$ satisfies the following asymptotic relations:
\begin{eqnarray}
\biggl[ (r \rightarrow R) \, \Longrightarrow \, (k_{\delta_2} \longrightarrow k_h) \biggr] 
\qquad \text{and} \qquad 
\biggl[ (r \rightarrow 0) \, \Longrightarrow \, (k_{\delta_2} \longrightarrow k_0) \biggr] \, .
\label{delta_energy_3b}
\end{eqnarray}
Using~\eqref{delta_energy_3} together with \eqref{def_k_c}, gives us a formula for $\delta_2$:
\begin{eqnarray}
\delta_2
= \frac{2 \, \pi}{k_{\delta_2}}
= 2 \, \pi \,
\left[ 
\Lambda \, \left( \frac{2 \, \pi}{h} \right)^{-2/3}
+ \left( 1 - \Lambda \right) \, \left( \frac{2 \, \pi}{L} \right)^{-2/3}
\right]^{3/2} 
= \left[ 
\Lambda \, h^{2/3}
+ \left( 1 - \Lambda \right) \, L^{2/3}
\right]^{3/2} \, .
\label{eqn:delta-2}
\end{eqnarray}
The new ROM lengthscale, $\delta_2$, defined in~\eqref{eqn:delta-2}, depends on the FOM mesh size, $h$, the ROM dimension, $r$, the total number of ROM basis functions, $R$, the eigenvalues, $\lambda_i$, and the computational domain characterisitic lengthscale, $L$.

\begin{remark}[Asymptotic Behavior]
We note that the novel ROM lengthscale, $\delta_2$, satisfies the following two natural asymptotic relations:
\begin{itemize} \itemsep5pt
    \item As $r$ approaches $R$, $\delta_2$ approaches $h$.
    
    This simply says that, as expected, when the ROM dimension, $r$, approaches the maximal ROM dimension, $R$, the ROM lengthscale, $\delta_2$, approaches the minimal FOM lengthscale, $h$ (i.e., the spatial mesh size).
    
    \item As $r$ approaches $1$, $\delta_2$ approaches $L$.
    
    This simply says that, as expected, when the ROM dimension, $r$, approaches the minimal ROM dimension, $1$, the ROM lengthscale, $\delta_2$, approaches the maximal FOM lengthscale, $L$ (i.e., the dimension of the computational domain).
    
\end{itemize}
\label{remark:delta-2-asymptotic}
\end{remark}

\section{Numerical Results}
    \label{sec:numerical-results}

In this section, we perform a numerical investigation of the two lengthscales discussed in Section~\ref{sec:rom-lengthscale}: the standard dimensionality-based ROM lengthscale, $\delta_1$, defined in~\eqref{eqn:delta-1}, and the new energy-based ROM lengthscale, $\delta_2$, defined in~\eqref{eqn:delta-2}.
To this end, we use two fundamentally different ROMs:
the ML-ROM presented in Section~\ref{sec:ml-rom} and the EFR-ROM presented in Section~\ref{sec:efr-rom}.
In each type of ROM, we use the two ROM lengthscales and compare the results.
To this end, in our numerical investigation, we test four types of models:
(i) ML-ROM1, which is the ML-ROM in which the velocity scale, $U_{ML}$, is set equal to $\delta_1$;
(ii) ML-ROM2, which is the ML-ROM in which the velocity scale, $U_{ML}$, is set equal to $\delta_2$;
(iii) EFR-ROM1, which is the EFR-ROM in which the filter radius, $\delta$, is set equal to $\delta_1$; and
(iv) EFR-ROM2, which is the EFR-ROM in which the filter radius, $\delta$, is set equal to $\delta_2$. 
Thus, to compare the two lengthscales, we compare the four ROMs (i.e., ML-ROM1, ML-ROM2, EFR-ROM1, and EFR-ROM2) in the numerical simulation of the 3D turbulent channel flow at $Re_{\tau} = 395$.
{\it 
We emphasize that the goal of this section is not to find the best ML-ROMs and EFR-ROMs.
Instead, we 
investigate whether the two lengthscales are different and, if so, quantify their differences and how these differences impact the ML-ROM and EFR-ROM results.
}

\subsection{Numerical Setting}

\paragraph{FOM}
The computational domain is a rectangular box, $\Omega = (-2\pi,2\pi)\times (0,2)\times (-2\pi/3,2\pi/3)$. 
We enforce no slip boundary conditions 
on the walls at $y=0$ and $y=2$, and periodic boundary conditions in the $x-$ and $z-$directions.
We also use 
the forcing term 
$\bff = [1,0,0]^\top$
and the Reynolds number 
$Re_\tau =395$ ($Re=13,750$). 
To generate the snapshots, we run 
an LES model using the rNS-$\alpha$ scheme~\cite{rebholz2017global,rebholz2017accurate} with the time step size $\Delta t
= 0.002$.

\begin{figure}[H]
\centering
    \includegraphics[width=.32\textwidth]{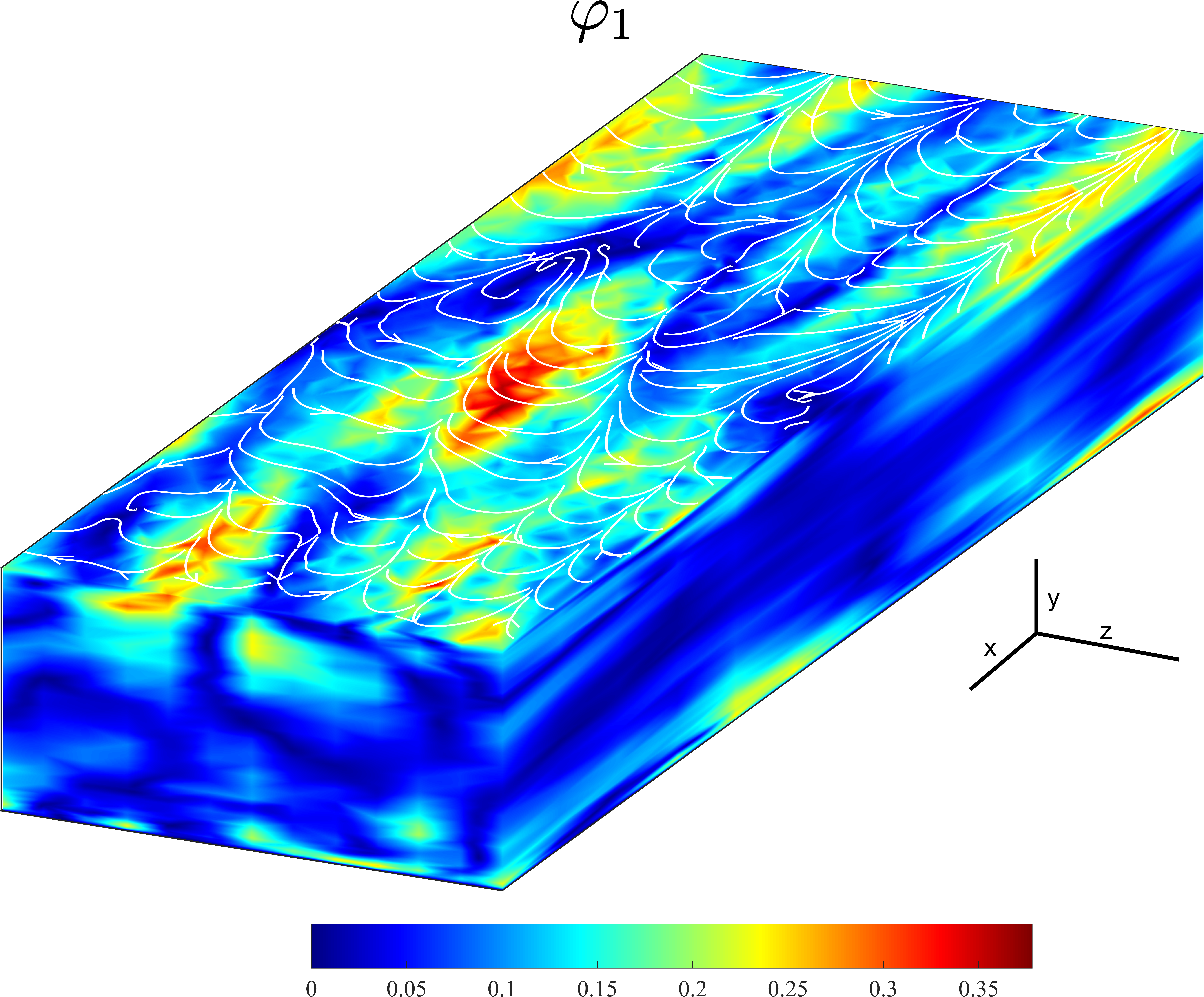}
    \includegraphics[width=.32\textwidth]{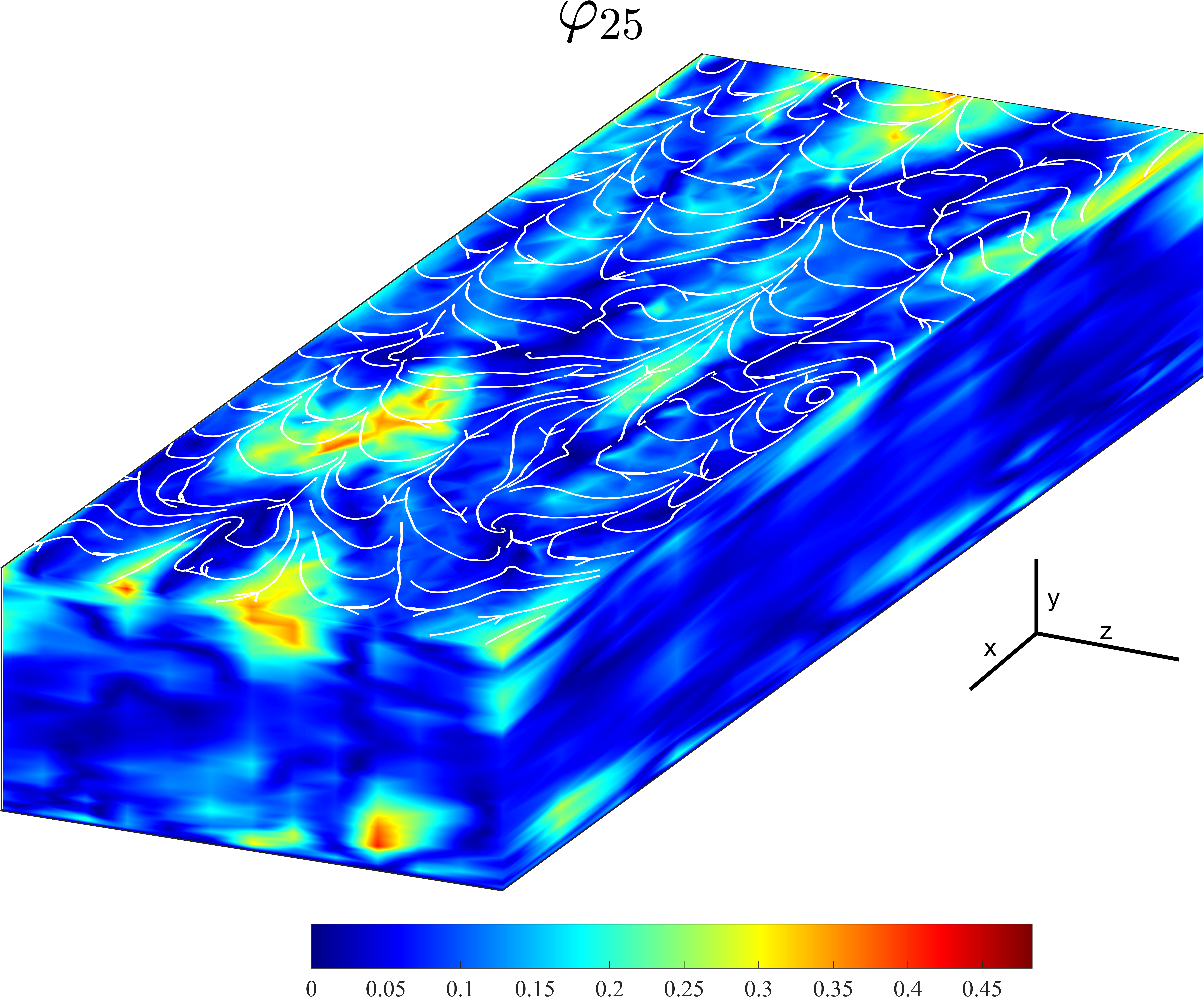}
    \includegraphics[width=.32\textwidth]{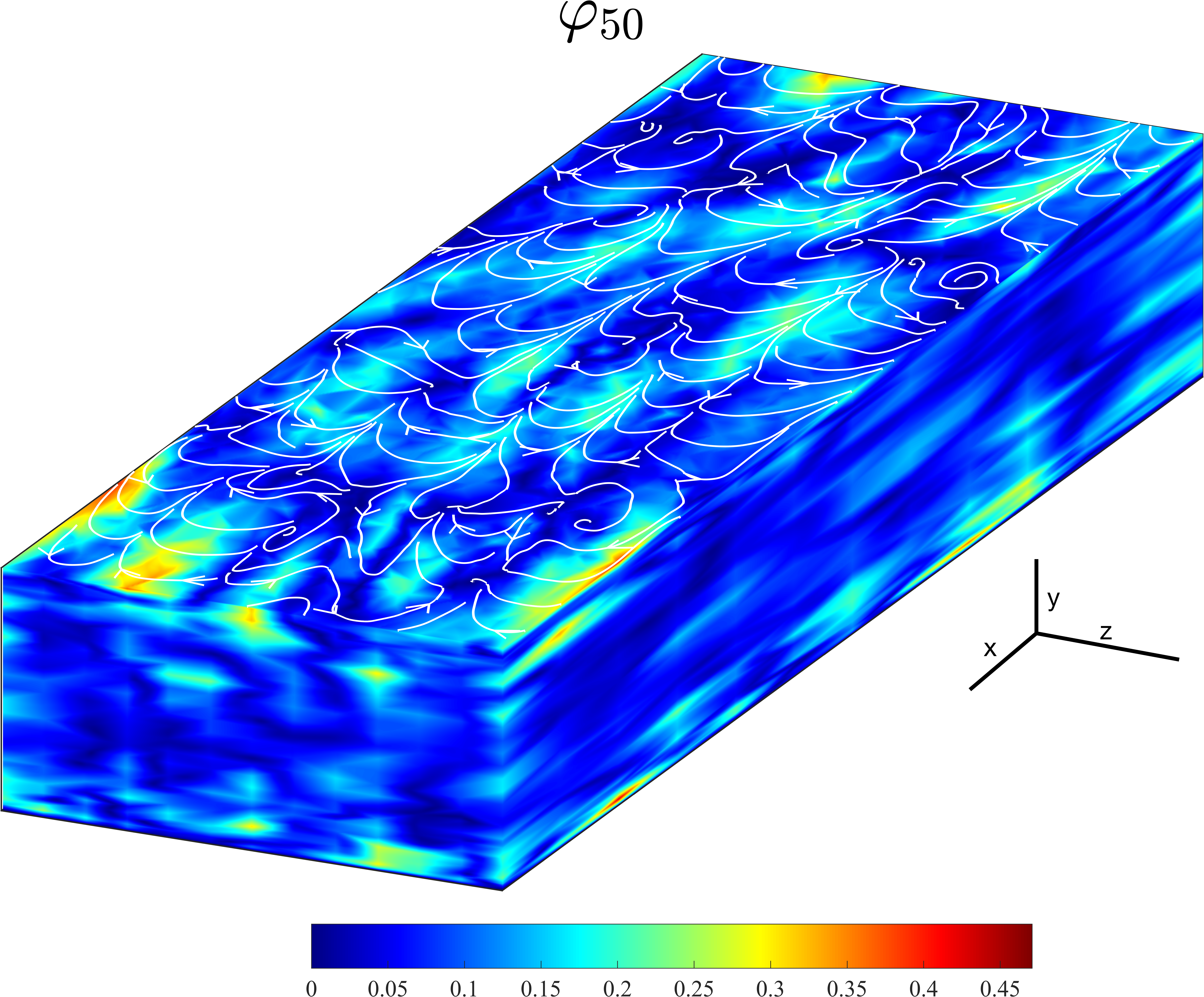}
    \caption{
    The 1st, 25th, and 50th POD modes
    }    
    \label{fig:pod}
\end{figure}

\paragraph{ROM}
We collect a total of $5000$ snapshots from $t = 60$ to $t = 70$ 
(i.e., from a time interval in which the flow is in the statistically steady regime) 
and use the POD to generate the ROM basis.
For illustrative purposes, we plot the magnitude fields of the ROM basis functions $\bphi_{1}, \bphi_{25}$, and $\bphi_{50}$ 
in Fig.~\ref{fig:pod}.

For the ROM time discretization, we utilize the commonly used linearized BDF2 temporal discretization with a time step size $\Delta t = 0.002$.
As the ROM initial conditions, we use the ROM projections of the LES approximations at $t = 60$ and $t = 60.002$.
For convenience, in our ROM simulations, $t=0$ corresponds to $t=60$ in the LES model.
 
To assess the ROMs' performance, we use two different criteria: 
(i) the time evolution of the kinetic energy, $E(t)$, and (ii) second-order statistics.

We define the ROM kinetic energy as follows:
\begin{eqnarray}
E(t) 
= \frac{1}{2}\int_\Omega\bigl(
u_1(\bx,t)^2+u_2(\bx,t)^2+u_3(\bx,t)^2 
\bigr) \, d\bx \, ,
\end{eqnarray}
where $u_1, u_2$, and $u_3$ are the components of the ROM velocity field approximation.

Following~\cite{rebholz2017accurate}, we consider the 
following two second-order statistics:
(i) the normalized root mean square (RMS) of the streamwise velocity component, $U_{RMS}$:
\begin{eqnarray}
U_{RMS} := \frac{\biggl|\widetilde{\mathbb{R}}_{11}-\frac{1}{3} \sum_{j=1}^3\widetilde{\mathbb{R}}_{jj}\biggr|^{1/2}}{u_{1,\tau}},
\end{eqnarray}
and (ii) the normalized streamwise-spanwise Reynolds stress tensor component, $\mathbb{R}_{12}$:
\begin{eqnarray}
\mathbb{R}_{12} := \frac{\widetilde{\mathbb{R}}_{12}}{u_{1,\tau}^2}.
\label{eqn:R_12}
\end{eqnarray}
In these second-order statistics, 
the Reynolds stress tensor components are calculated as follows: 
\begin{eqnarray}
\widetilde{\mathbb{R}}_{ij} = \left\langle \left\langle u_{i}u_{j}\right\rangle_s\right\rangle_t
-\left\langle \left\langle u_{i}\right\rangle_s\right\rangle_t \left\langle\left\langle u_{j}\right\rangle_s\right\rangle_t \, ,
\end{eqnarray}
where $\langle \cdot \rangle_s$ denotes spatial averaging, $\langle \cdot \rangle_t$ denotes time averaging, and $u_i$ are the components of the given ROM or FOM velocity field approximations.
The friction velocity, $u_{1,\tau}$, which is used in~\eqref{eqn:R_12}, is calculated by using the following formula:
\begin{align}
    u_{1,\tau} = 
\sqrt{\nu \, \frac{U_{\text{mean}}(y_{min})}{y_{min}} },
\end{align}
where 
$\nu$ is the kinematic viscosity,  
$y_{min}$ is the minimum positive $y$-value of the FOM mesh, and
$U_{\text{mean}} = \left\langle \left\langle 
u_1 \right\rangle_s\right\rangle_t$ the average ROM or FOM velocity flow profile.

\subsection{Numerical Results: Lengthscale Comparison}

In this section, we investigate the relative size of the two ROM lengthscales: (i) the standard dimensionality-based lengthscale, $\delta_1$, defined in~\eqref{eqn:delta-1}, and (ii) the novel energy-based ROM lengthscale, $\delta_2$, defined in~\eqref{eqn:delta-2}.
To calculate $\delta_2$ in equation~\eqref{eqn:delta-2}, we define the FOM global mesh size, $h$, as  $h=\max_{K\in\mathcal{K}}r_K$, where the mesh $\mathcal{K}$ is the set $\{K\}$ of tetrahedrons K and $r_K$ is the inradius of the local tetrahedron, $K$. 
For the test problem 
used in our numerical investigation, 
$h = 1.1 \times 10^{-1}$.
Furthermore, we define the maximal FOM lengthscale, $L$, as the vertical size of the computational domain, i.e., $L=2$.

In Table~\ref{table:delta}, we list the $\delta_1$ and $\delta_2$ values for $r$ values from $4$ to $50$.
These results show that the two ROM lengthscales have very different behaviors with respect to changes in $r$. 

\paragraph{Magnitude Behavior}
First, we notice that the $\delta_2$ magnitude is between 
one and 
two orders of magnitude larger than the $\delta_1$ magnitude.
Indeed, for the low $r$ values (i.e., $r = 4, 8, 16$), $\delta_2$ is more than an order of magnitude larger than $\delta_1$.
For the large $r$ values (i.e., $r = 32, 40, 50$), $\delta_2$ is still an order of magnitude larger than $\delta_1$, but the difference between the two lengthscales is smaller.

\paragraph{Asymptotic Behavior}
Second, the asymptotic behavior of the ROM lengthscales with respect to $r$ is fundamentally different:
As $r$ increases, the $\delta_1$ magnitude remains relatively unchanged.
In contrast, as $r$ increases, the $\delta_2$ magnitude decreases by almost one order of magnitude.
We emphasize that $\delta_2$ has the natural asymptotic behavior expected from a ROM lengthscale, as explained in Remark~\ref{remark:delta-2-asymptotic}:
(i) As $r$ 
increases toward its maximal value, i.e., $R$ (which is 
$R=804$ in our numerical investigation), $\delta_2$ approaches $h$.
Indeed, the results in Table~\ref{table:delta} show that for the maximal $r$ value (i.e., for $r = 50$), $\delta_2$ achieves the smallest value, $\delta_2 = 4.32 \times 10^{-1}$, which is the same order of magnitude as the FOM mesh size, $h = 1.1 \times 10^{-1}$. 
(ii) As $r$ approaches its minimal value 
(i.e., $1$), $\delta_2$ approaches $L$.
Indeed, the results in Table~\ref{table:delta} show that for the minimal $r$ value (i.e., for $r = 1$), $\delta_2$ achieves the largest value, $\delta_2 = 1.63 \times 10^{0}$, which is the same order of magnitude as the size of the computational domain, $L = 2.0 \times 10^{0}$.

\begin{table}[H]
    \centering
    \begin{tabular}{c|c c c c c c c c c c c c }
    \hline\hline
        $r$ & 4 & 8 & 16 & 32 &40  & 50 
    \\ \hline
        $\delta_1$ &4.64e-02 &4.65e-02 
        &4.68e-02  &4.68e-02 &4.66e-02   &4.62e-02  
    \\ 
        $\delta_2$
        &1.63e00 &1.41e+00 &1.08e+00  &6.84e-01 
        &5.56e-01  &4.32e-01        
    \\     
    \hline 
    \end{tabular}
    \caption{ROM lengthscale values for different $r$ values.
    }
    \label{table:delta}
\end{table}

\subsection{Numerical Results: ML-ROM Investigation}
    \label{sec:numerical-results-ml}

\begin{figure}[H]
\centering
    \includegraphics[width=.45\textwidth]{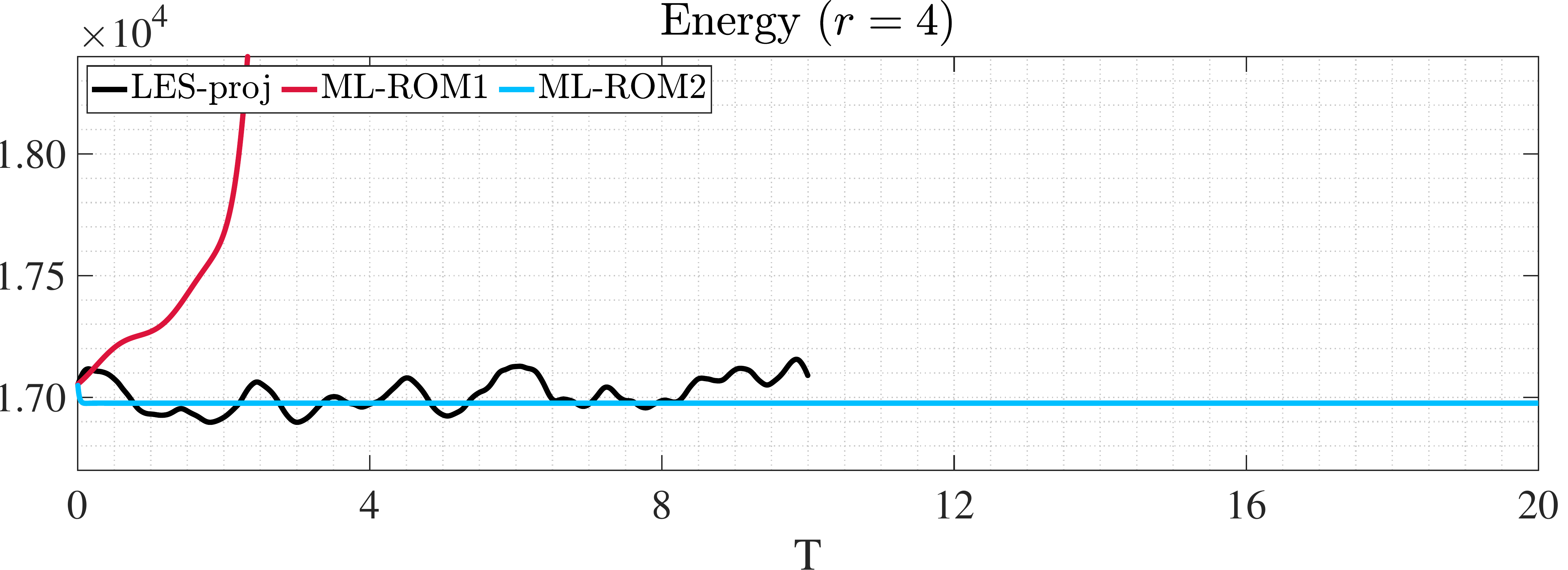}
    \includegraphics[width=.45\textwidth]{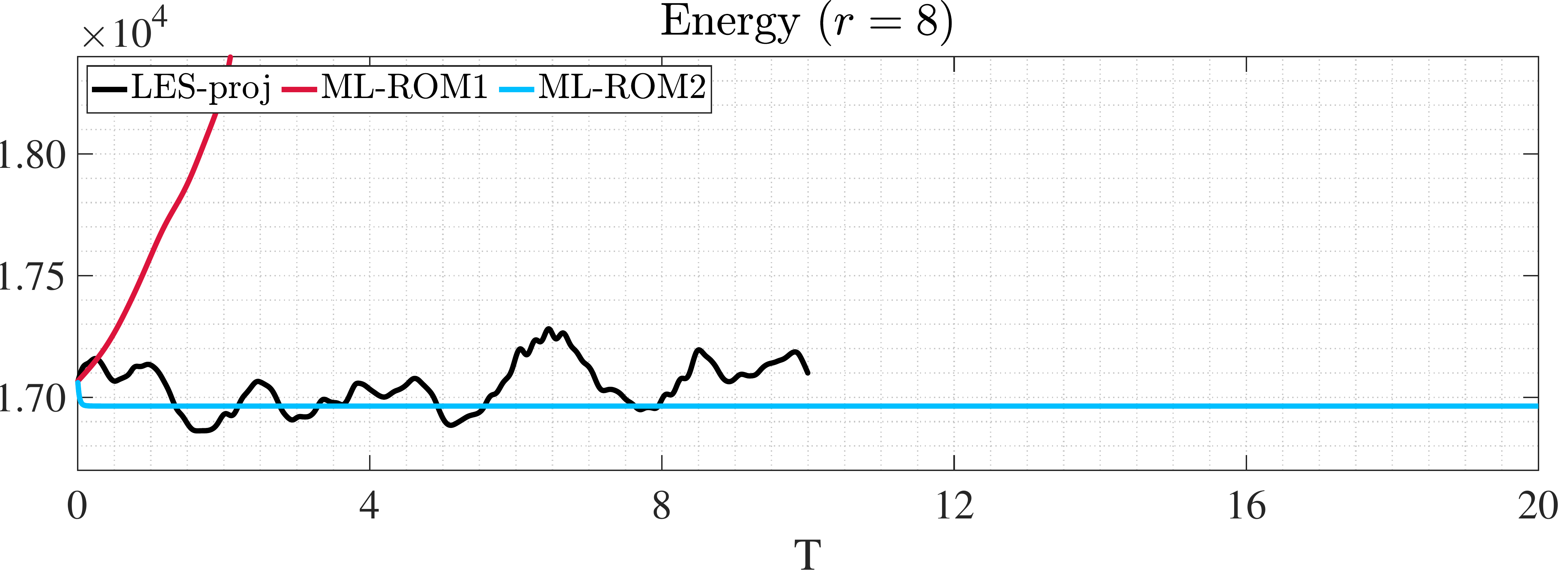}
    \includegraphics[width=.45\textwidth]{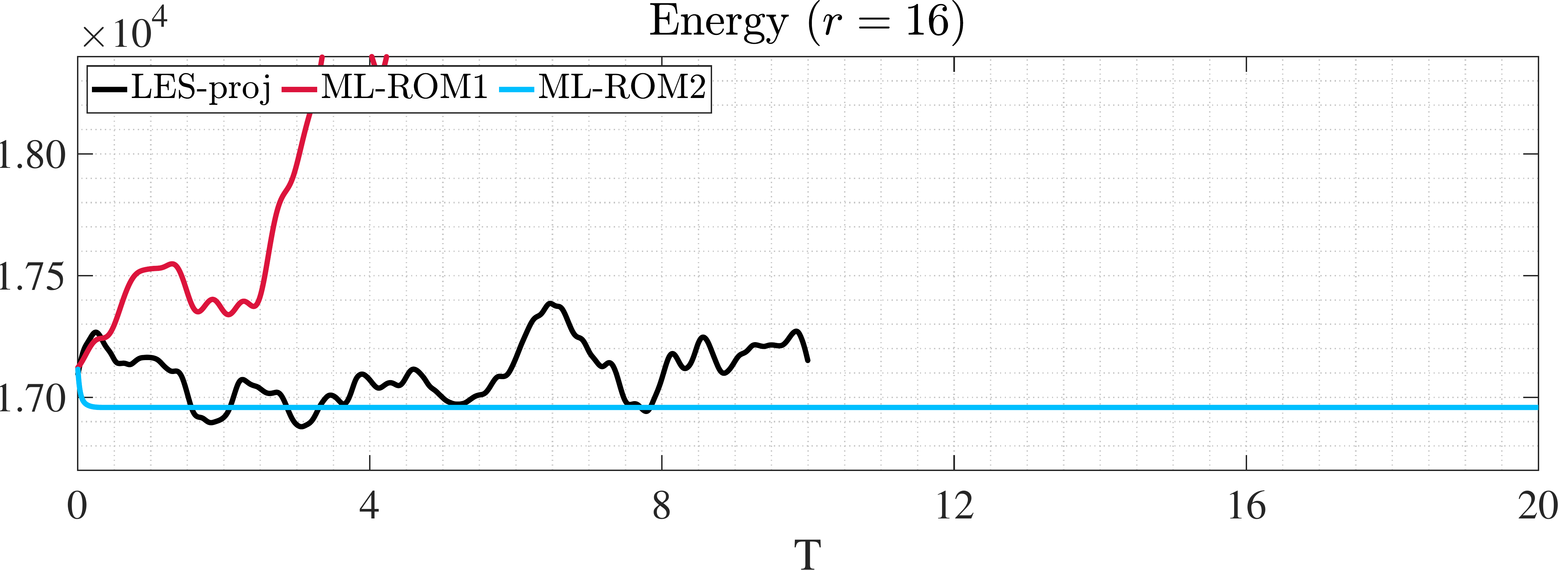}
    \includegraphics[width=.45\textwidth]{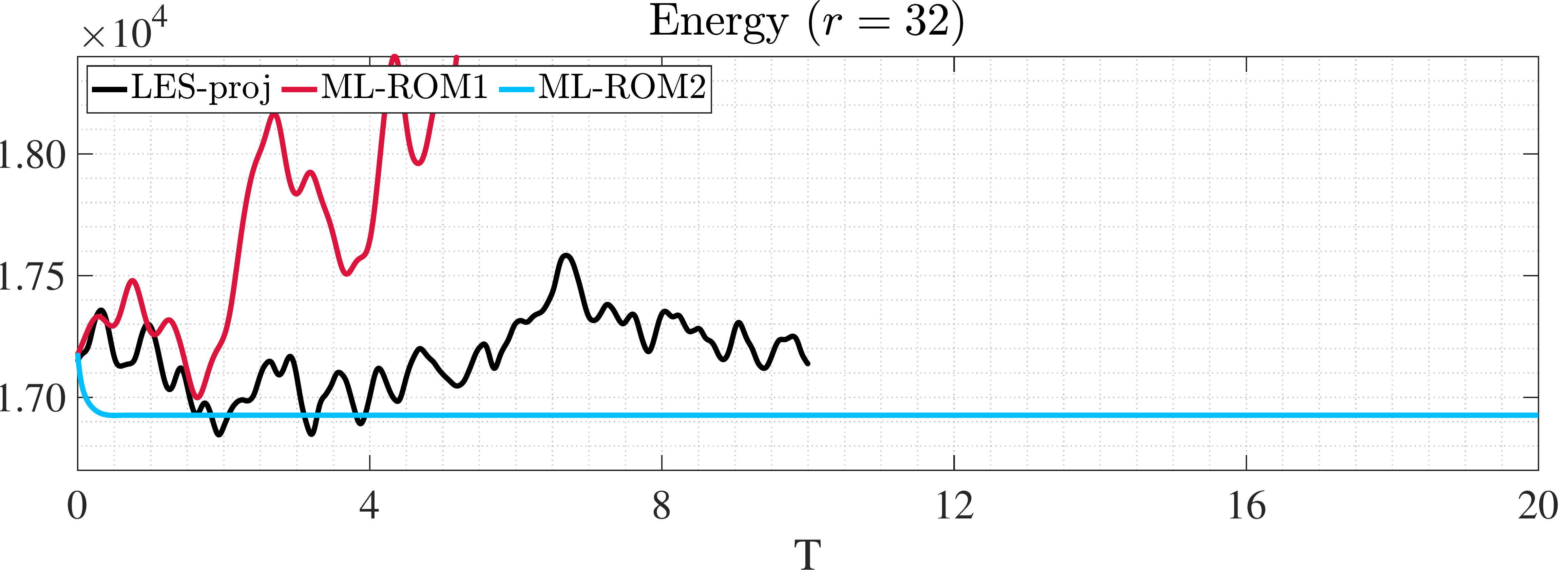}
    \includegraphics[width=.45\textwidth]{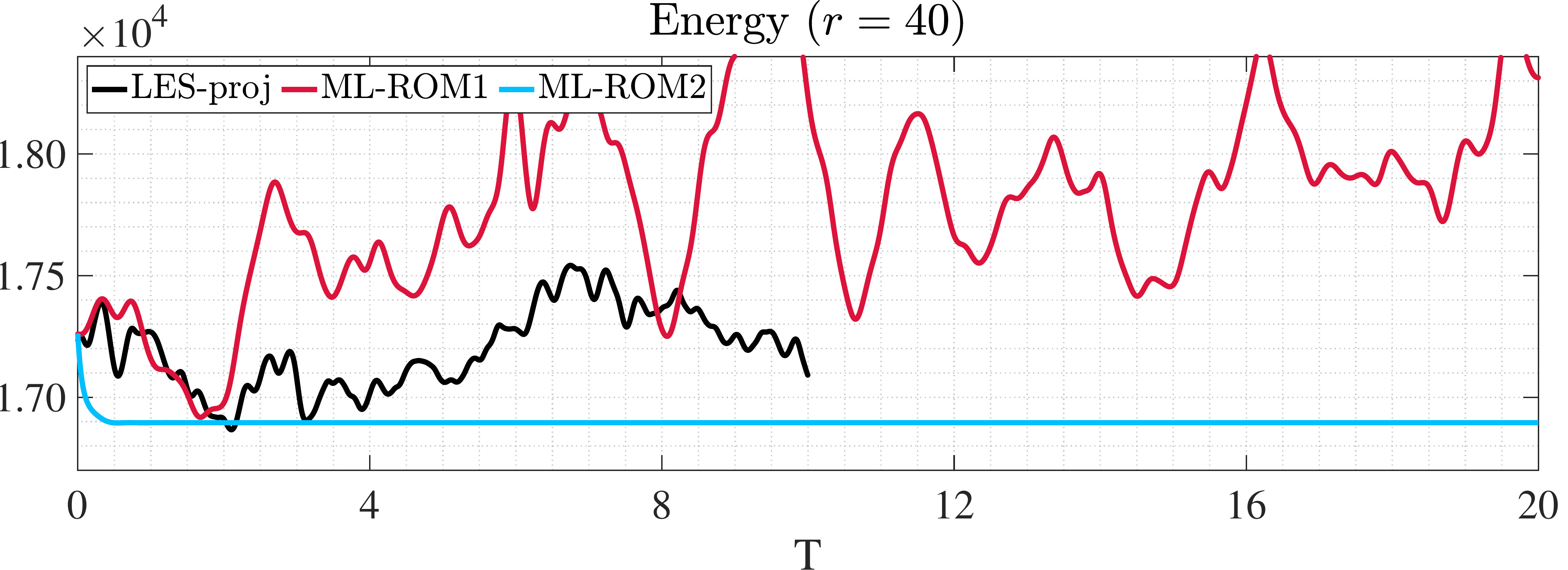}
    \includegraphics[width=.45\textwidth]{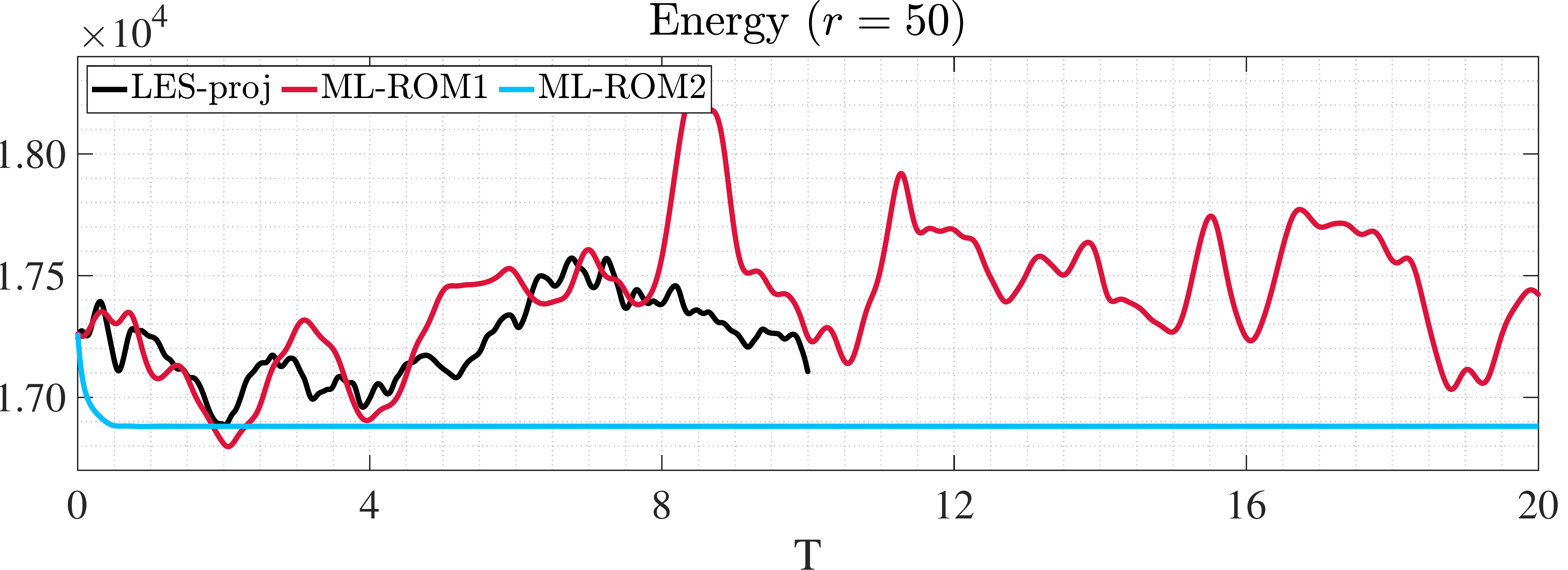}
    \caption{Time evolution of the ML-ROM kinetic energy for $\alpha=6\times 10^{-3}$
    }    
    \label{fig:ke-alpha-1}
\end{figure}

\begin{figure}[H]
\centering
    \includegraphics[width=.45\textwidth]{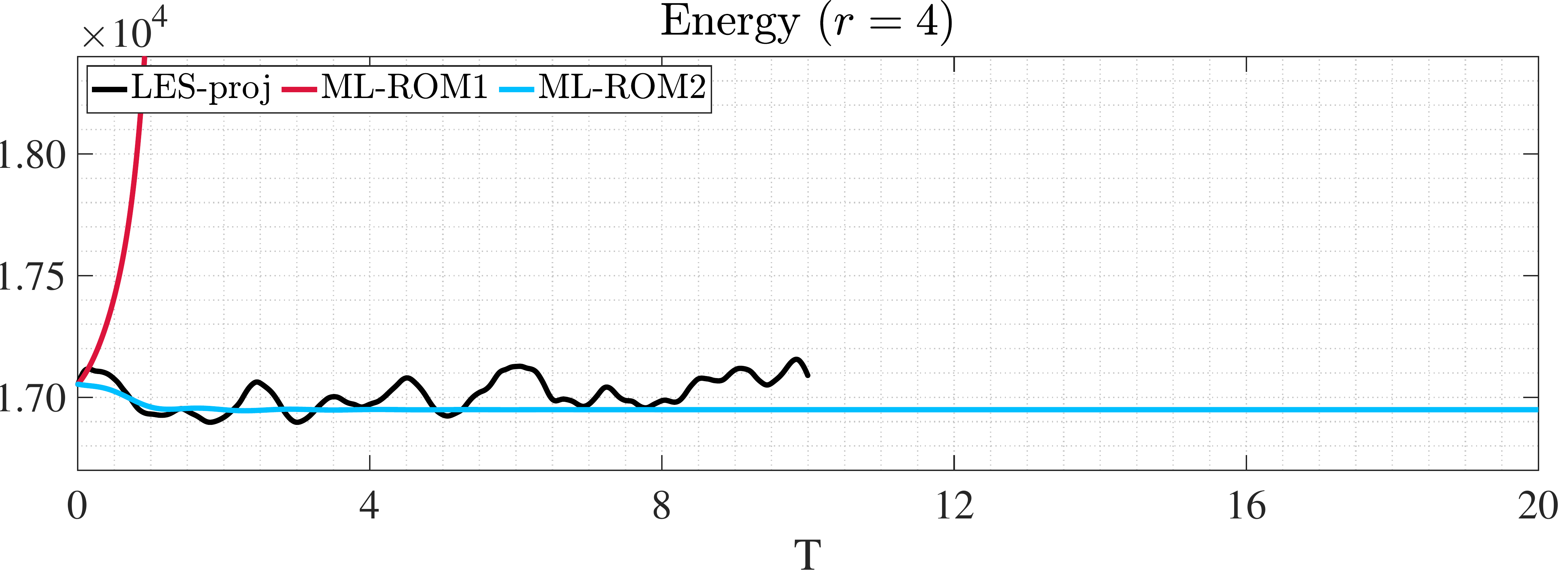}
    \includegraphics[width=.45\textwidth]{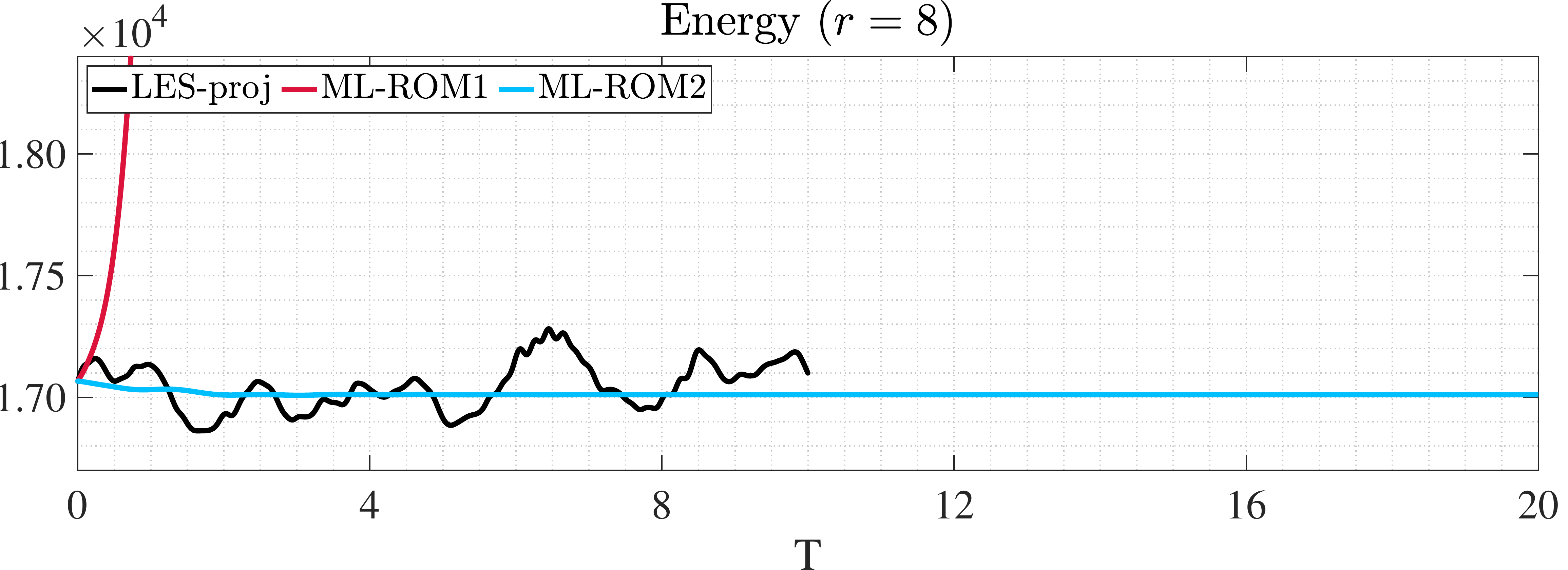}
    \includegraphics[width=.45\textwidth]{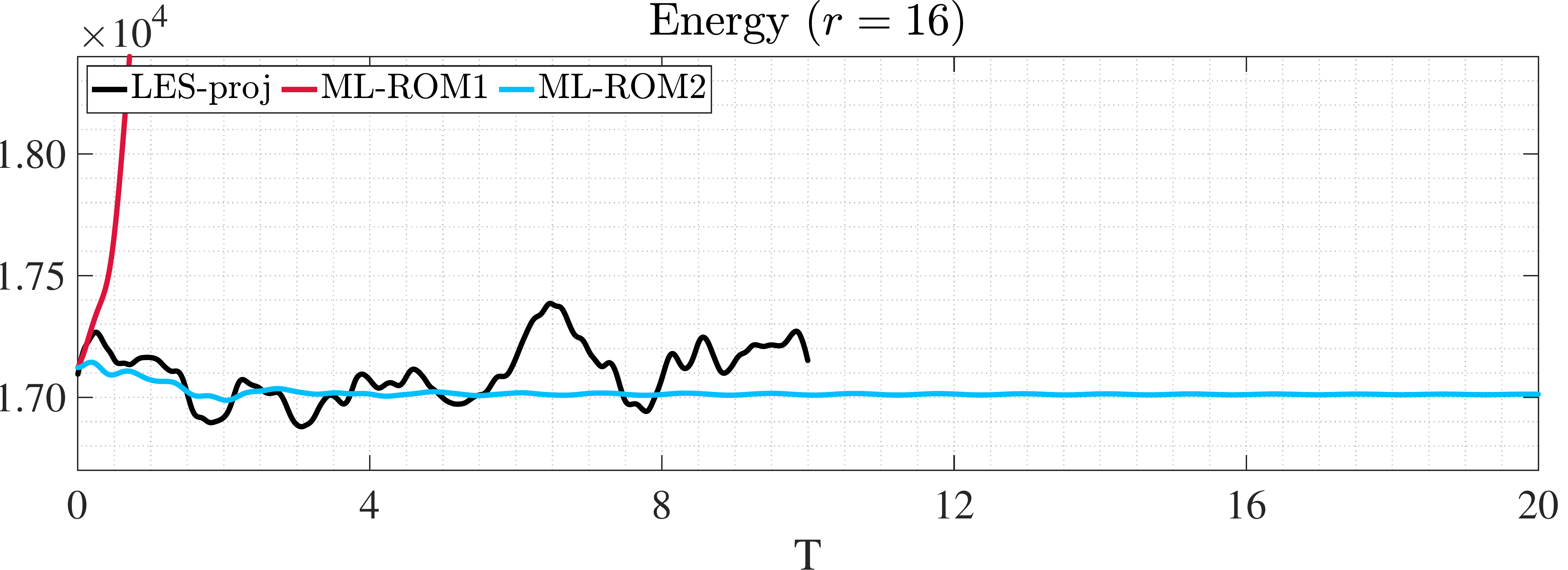}
    \includegraphics[width=.45\textwidth]{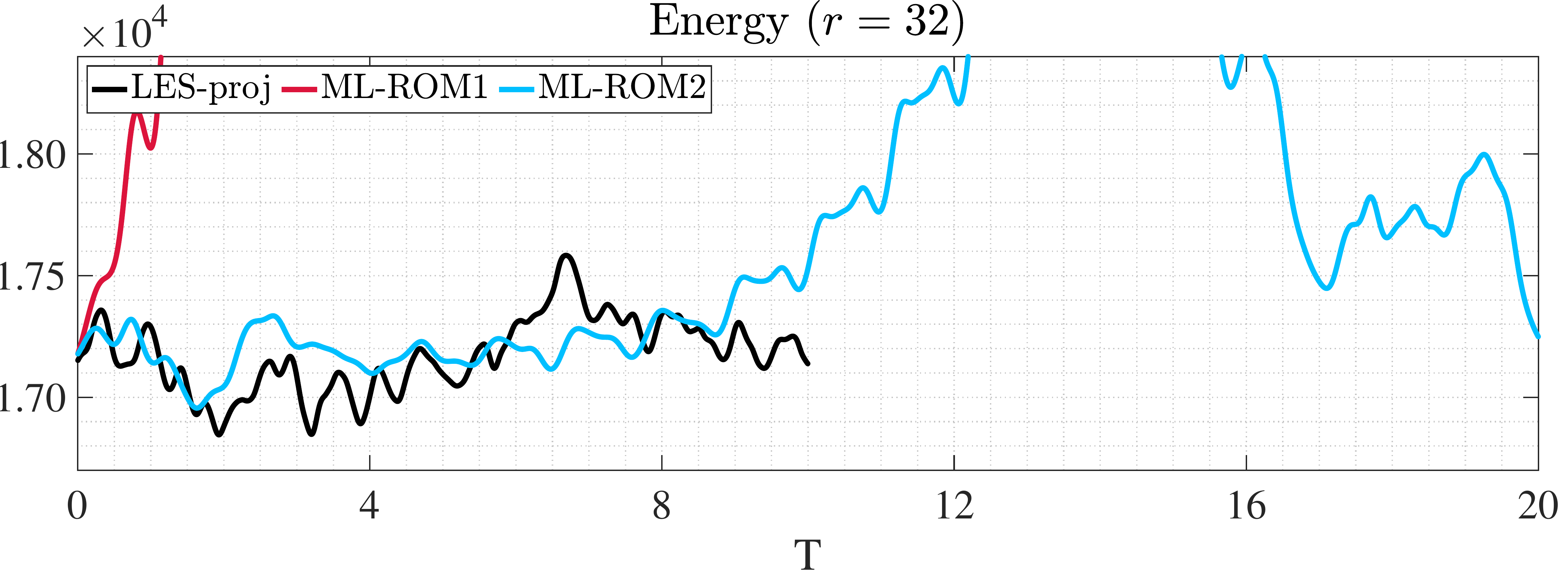}
    \includegraphics[width=.45\textwidth]{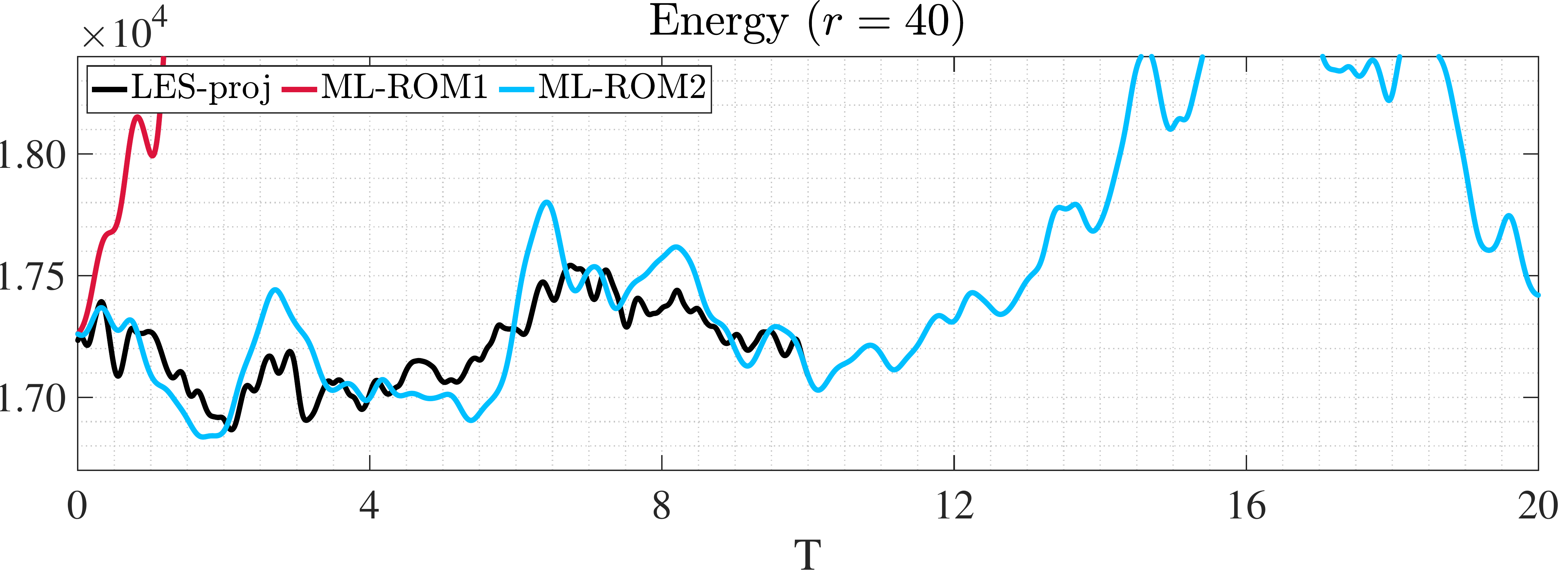}
    \includegraphics[width=.45\textwidth]{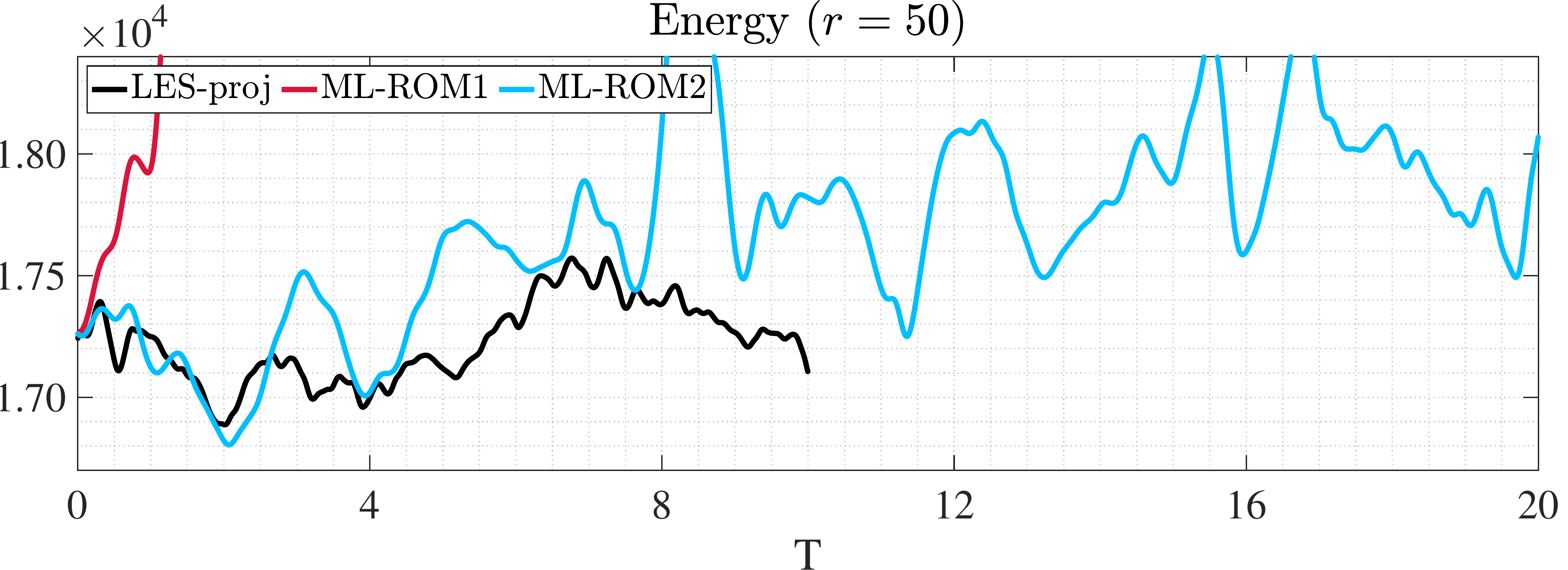}
    \caption{Time evolution of the ML-ROM kinetic energy for $\alpha=6\times 10^{-4}$
    }    
    \label{fig:ke-alpha-3}
\end{figure}

\begin{figure}[H]
\centering
     \begin{subfigure}[b]{0.48\textwidth}
         \centering
    \includegraphics[width=.45\textwidth]{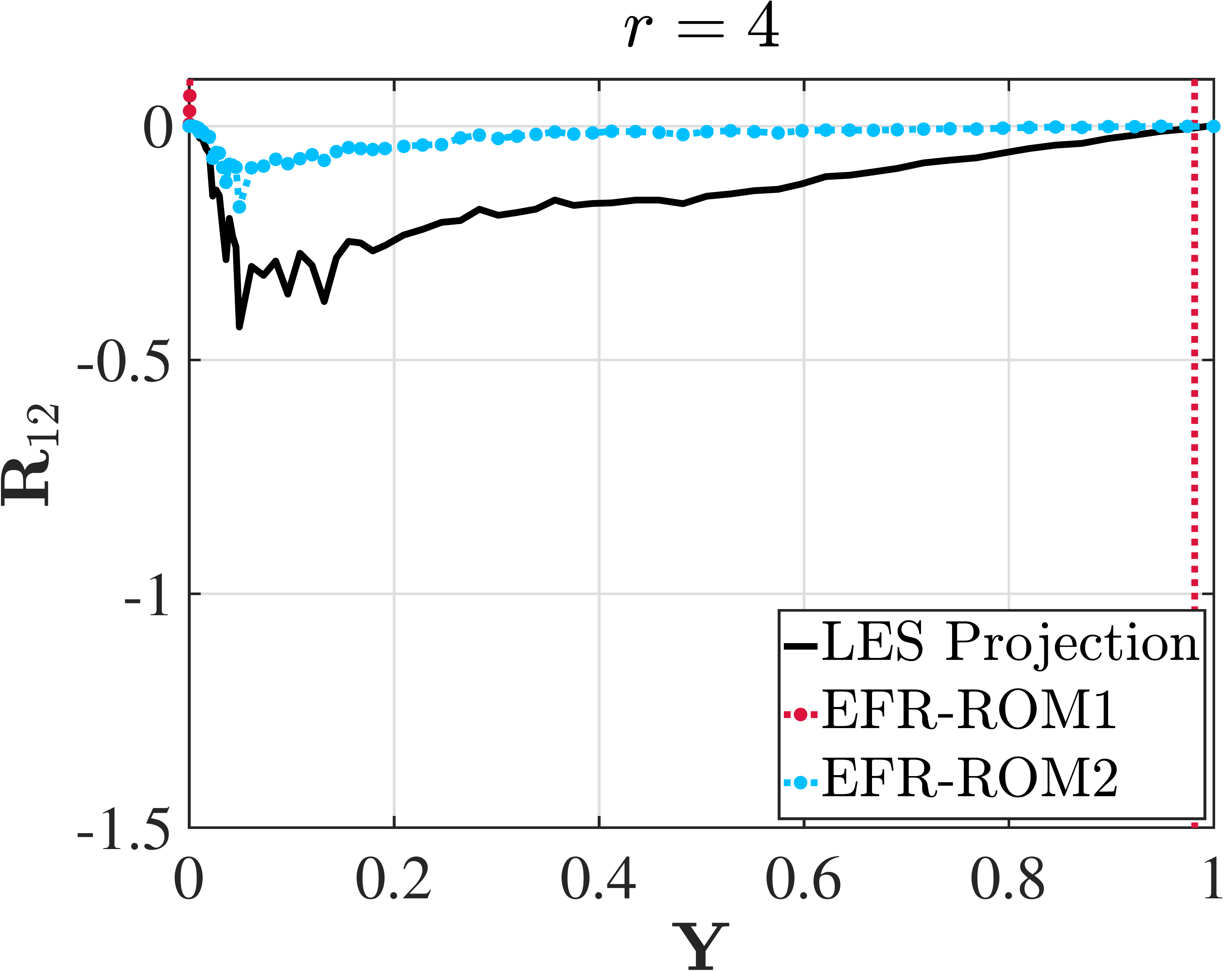}
    \includegraphics[width=.45\textwidth]{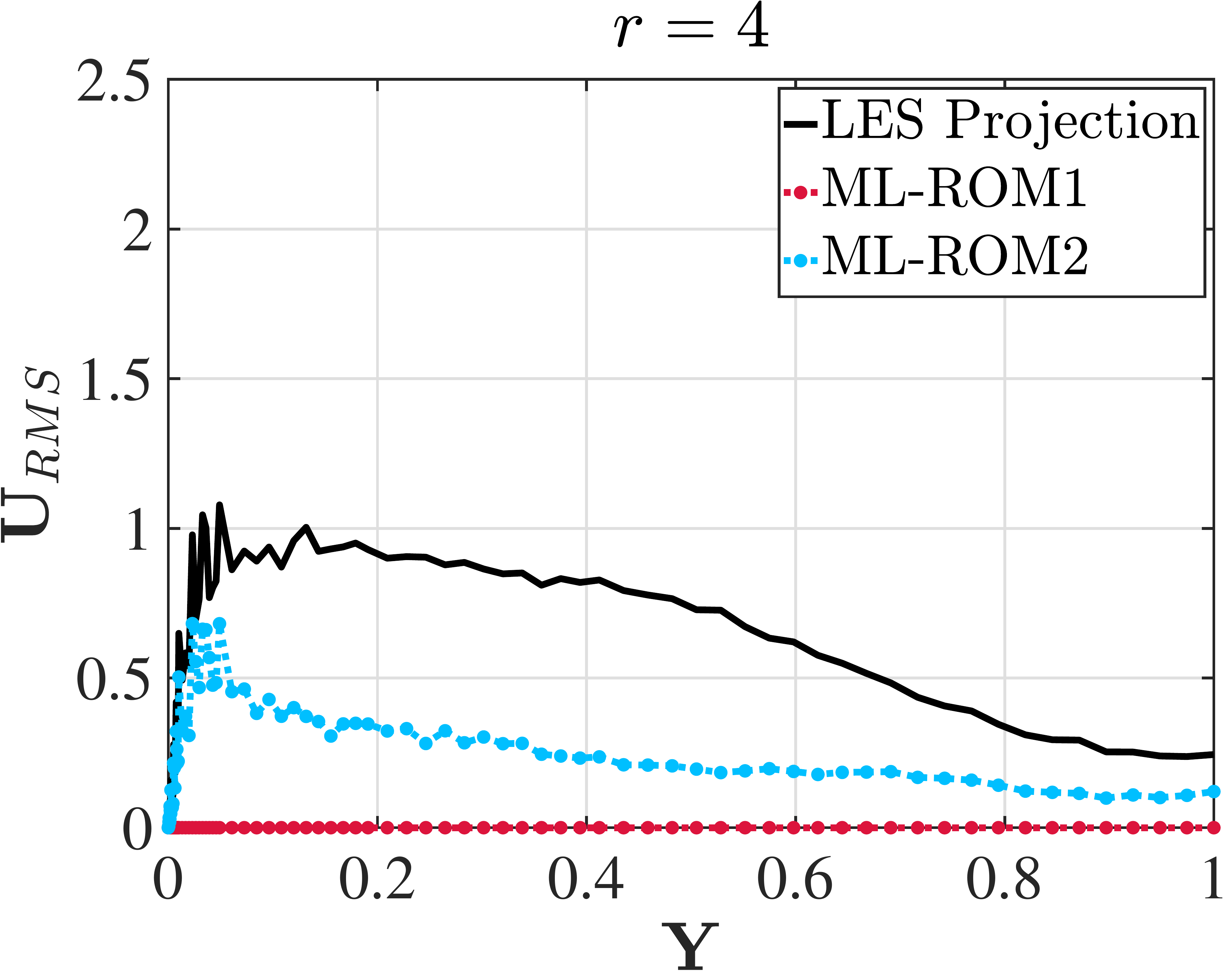}         \caption{$r=4$}
         \label{fig:stat-r-4}
     \end{subfigure}
     \begin{subfigure}[b]{0.48\textwidth}
         \centering
    \includegraphics[width=.45\textwidth]{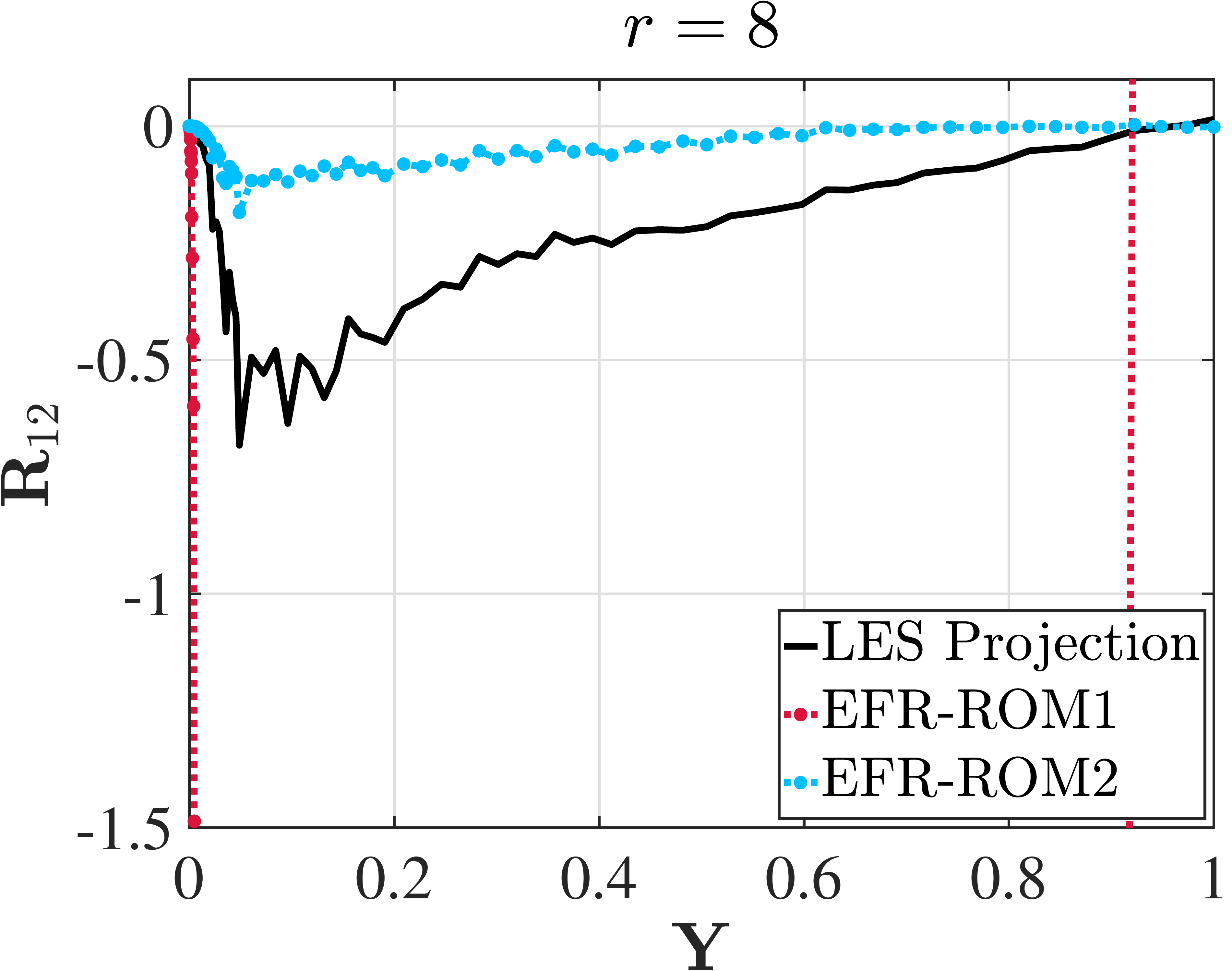}
    \includegraphics[width=.45\textwidth]{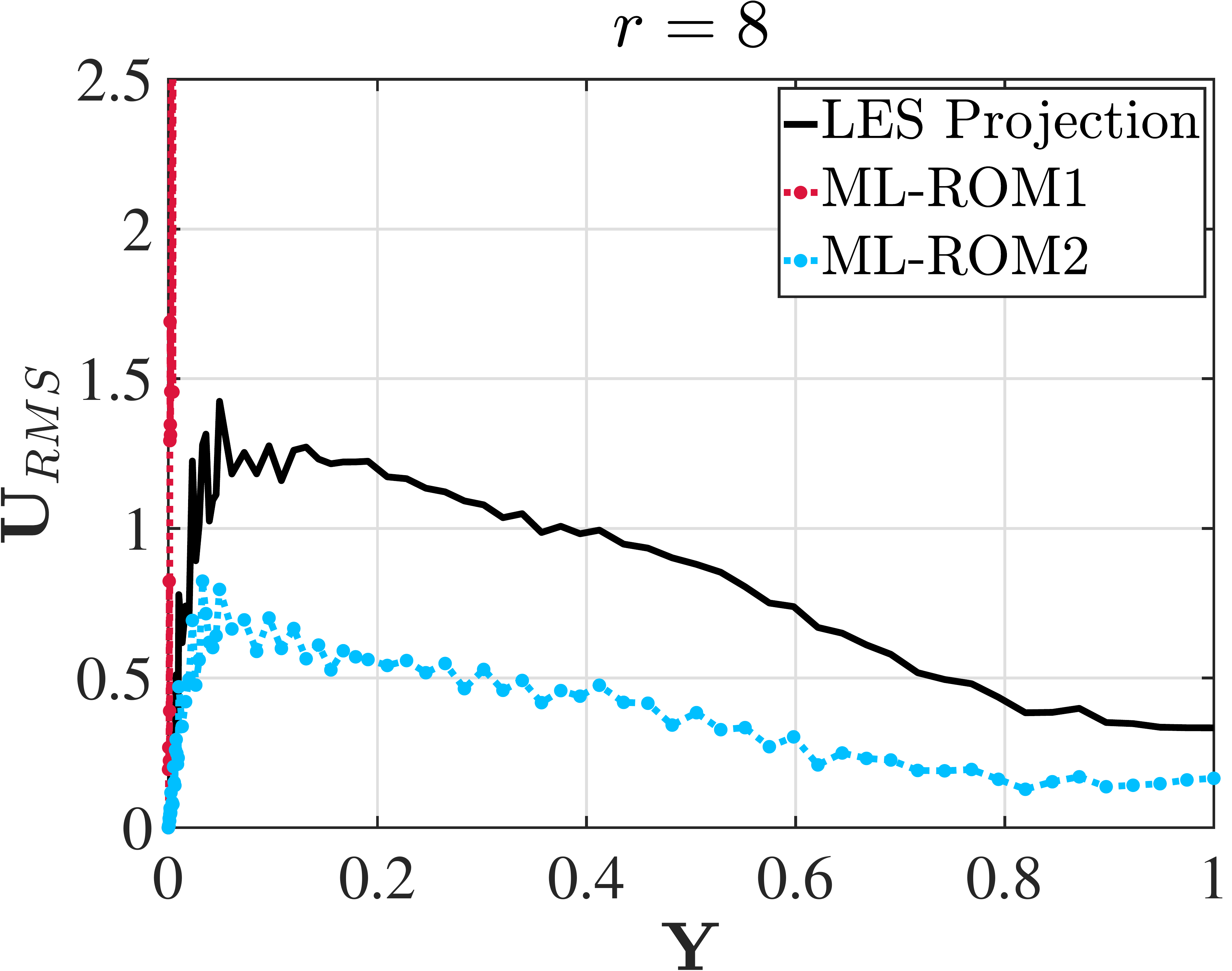}         \caption{$r=8$}
         \label{fig:stat-r-8}
     \end{subfigure}
     \begin{subfigure}[b]{0.48\textwidth}
         \centering
    \includegraphics[width=.45\textwidth]{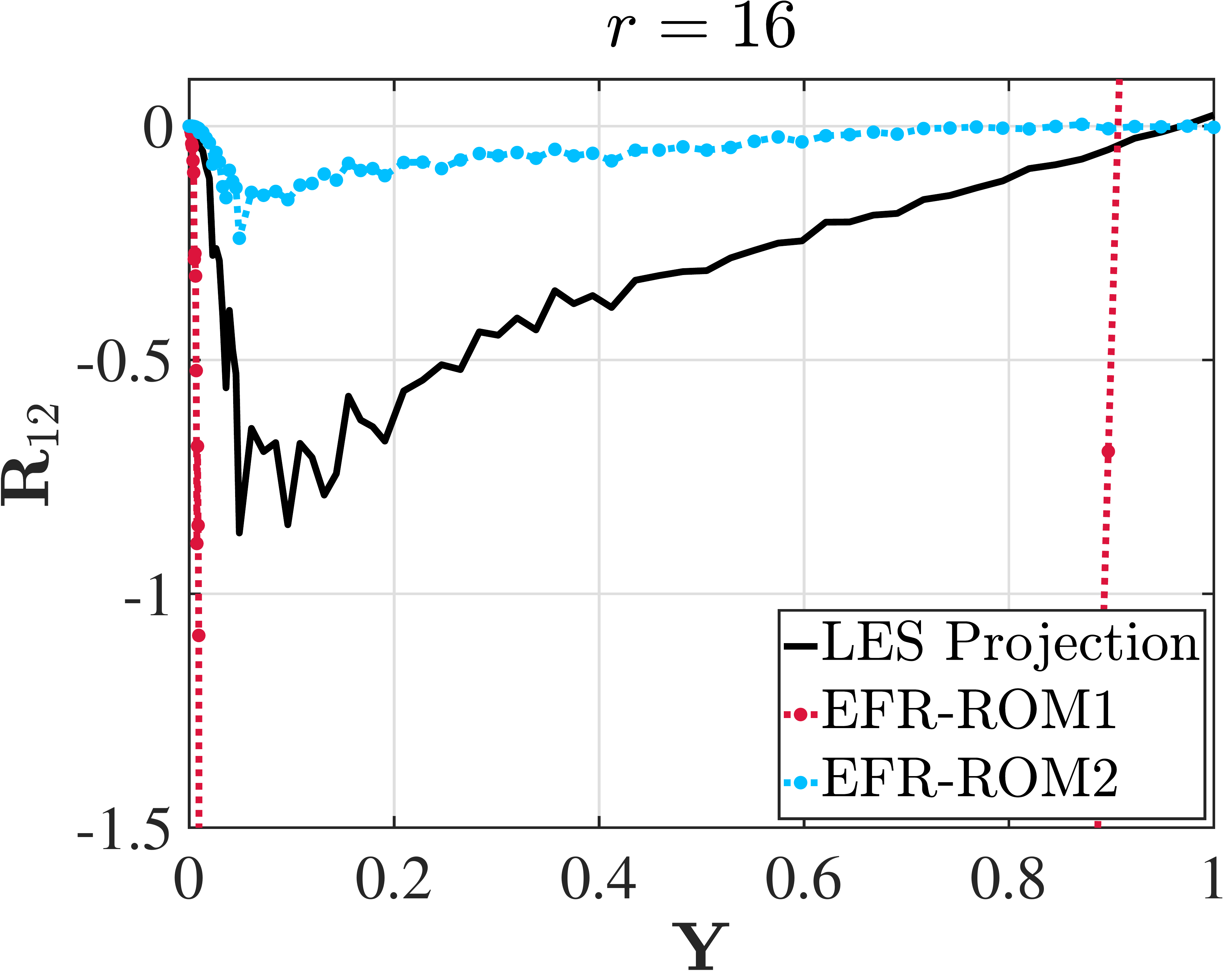}
    \includegraphics[width=.45\textwidth]{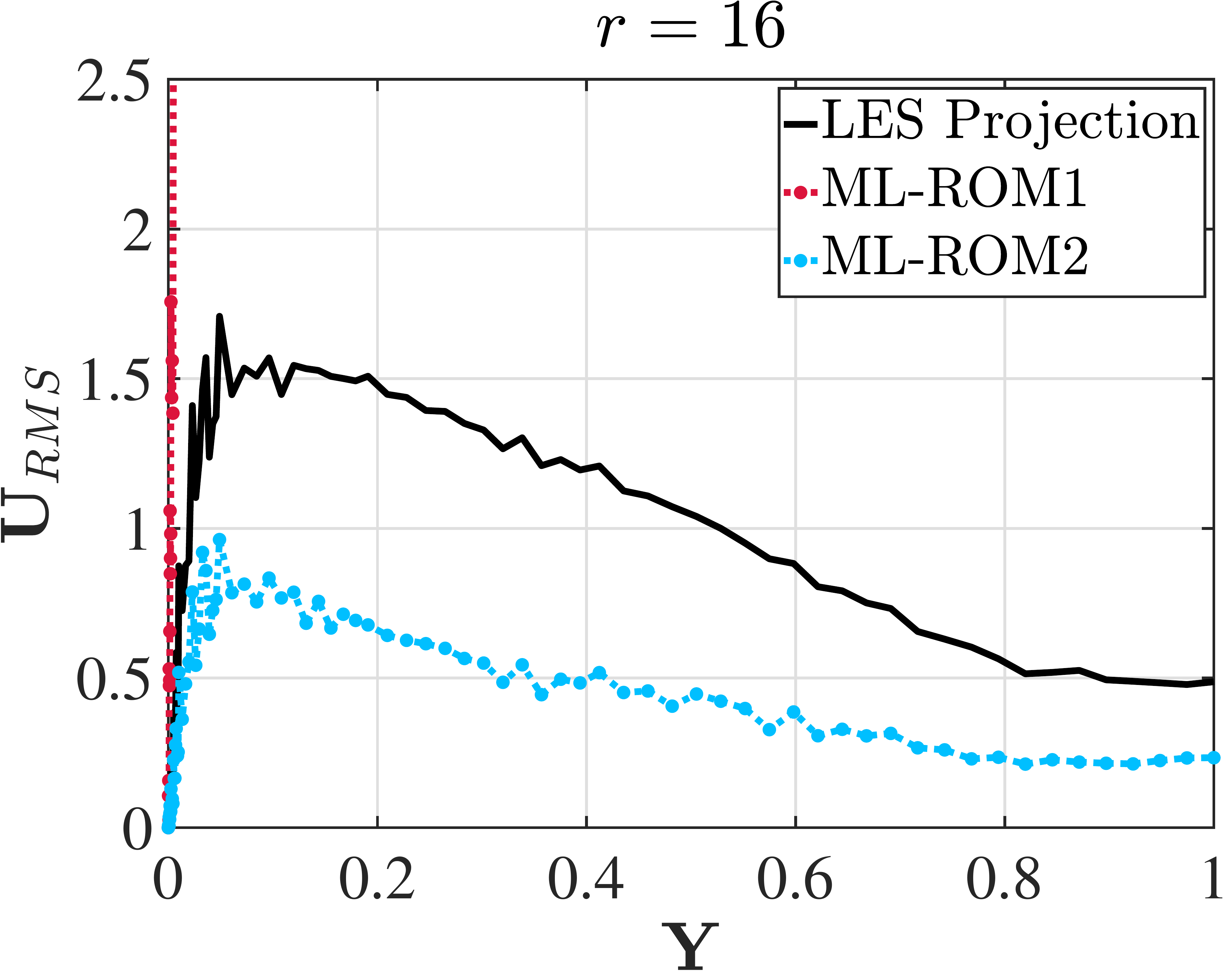}          \caption{$r=16$}
         \label{fig:stat-r-16}
     \end{subfigure}
     \begin{subfigure}[b]{0.48\textwidth}
         \centering
    \includegraphics[width=.45\textwidth]{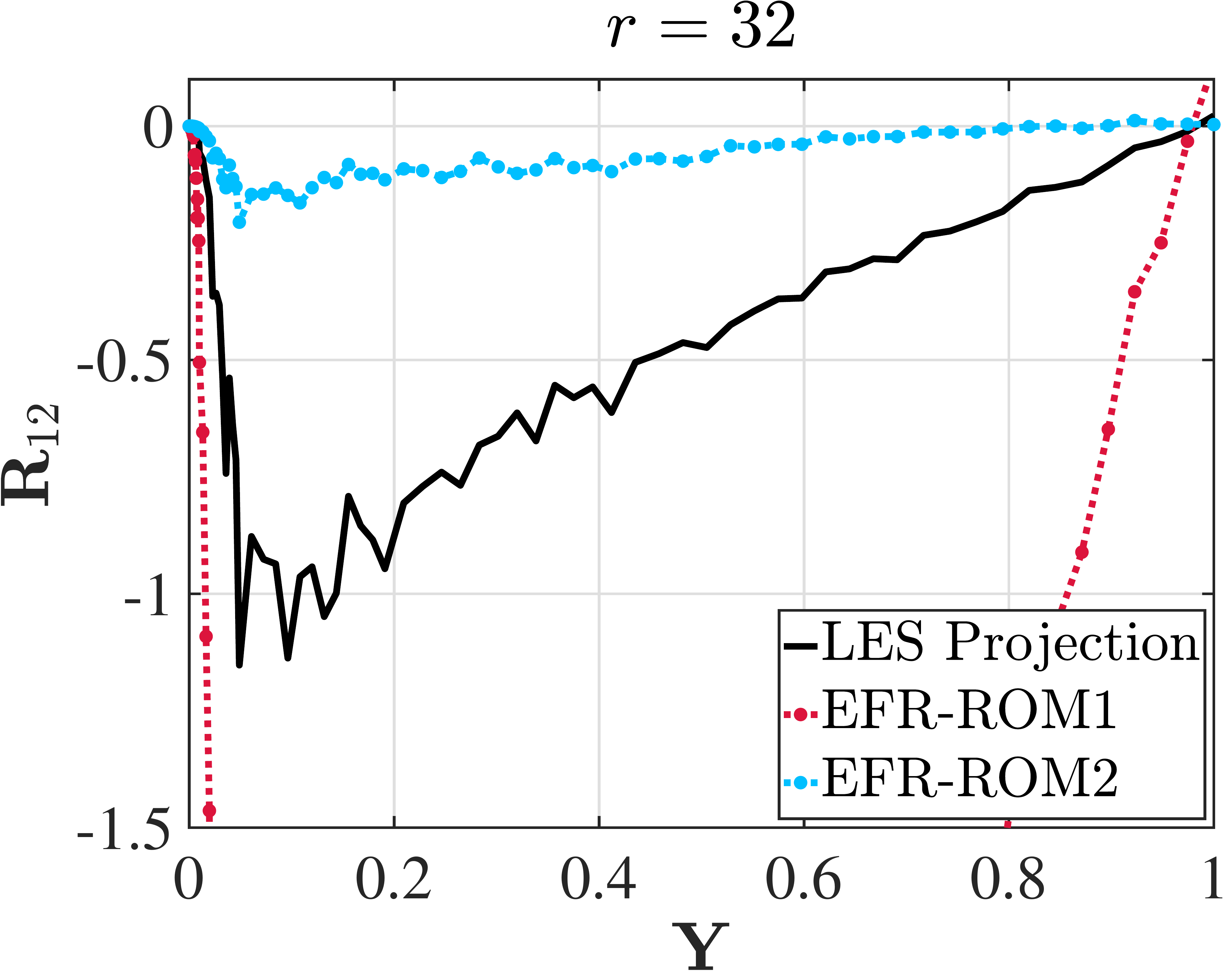}
    \includegraphics[width=.45\textwidth]{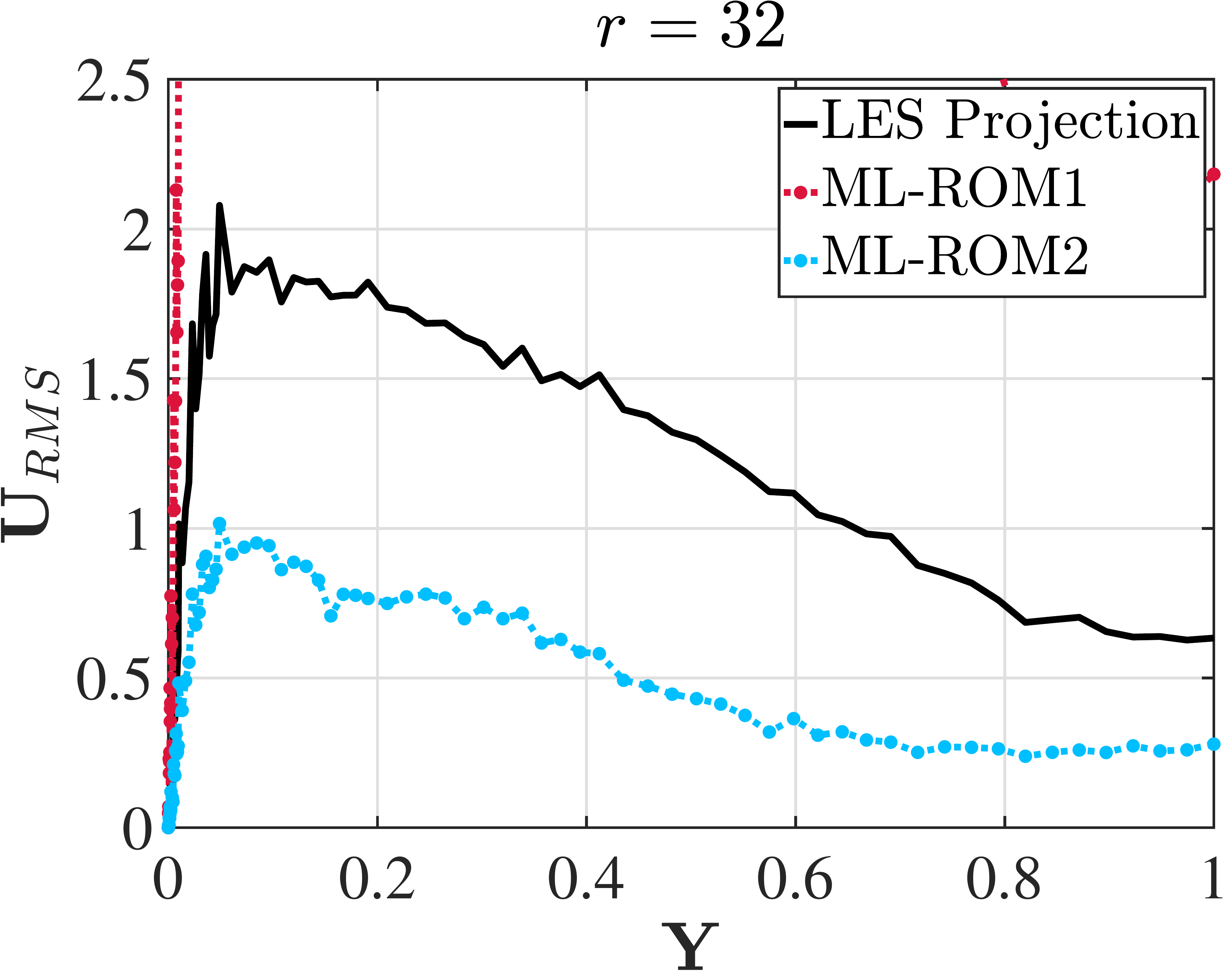}       \caption{$r=32$}
         \label{fig:stat-r-32}
    \end{subfigure}     
     \begin{subfigure}[b]{0.48\textwidth}
         \centering
    \includegraphics[width=.45\textwidth]{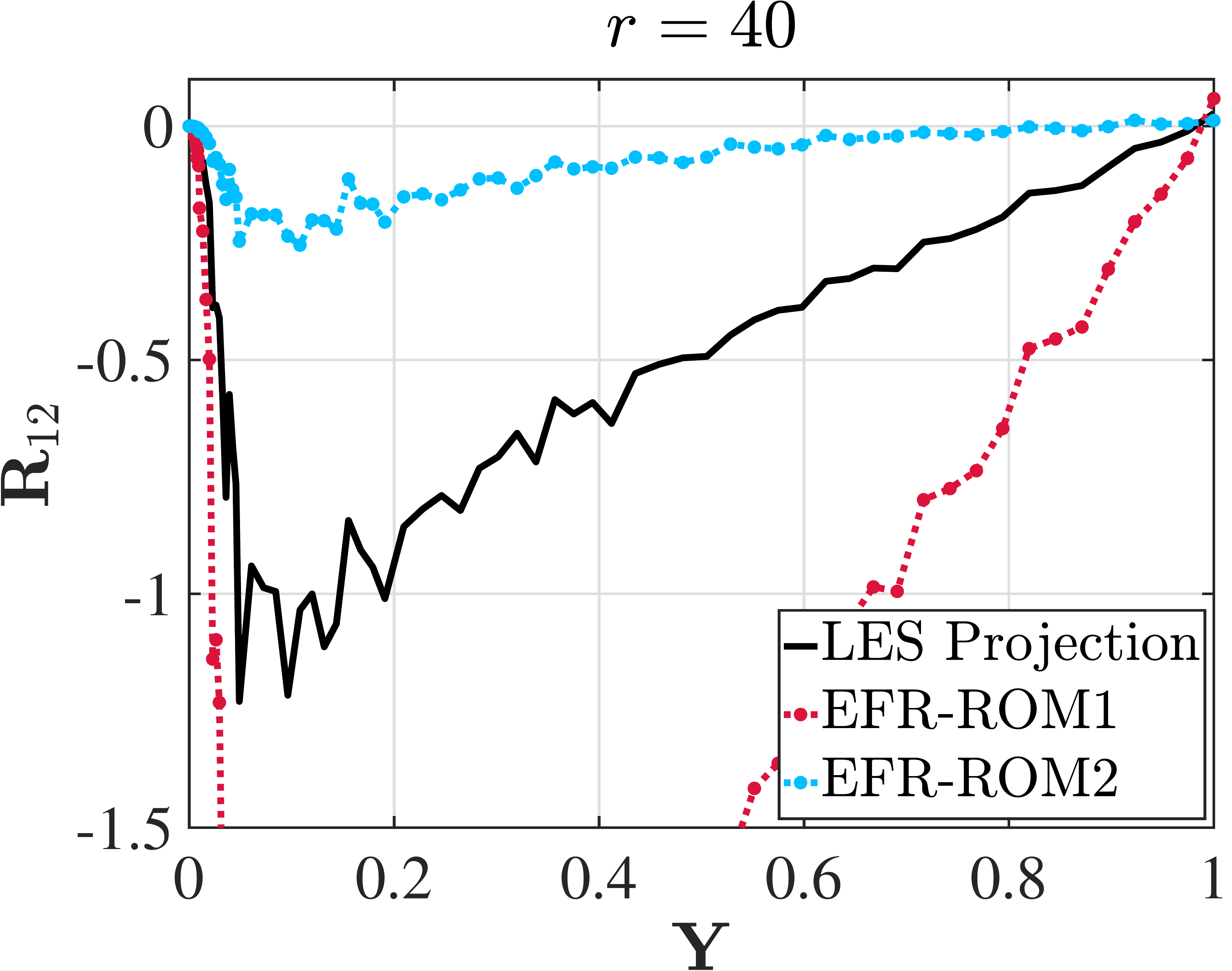}
    \includegraphics[width=.45\textwidth]{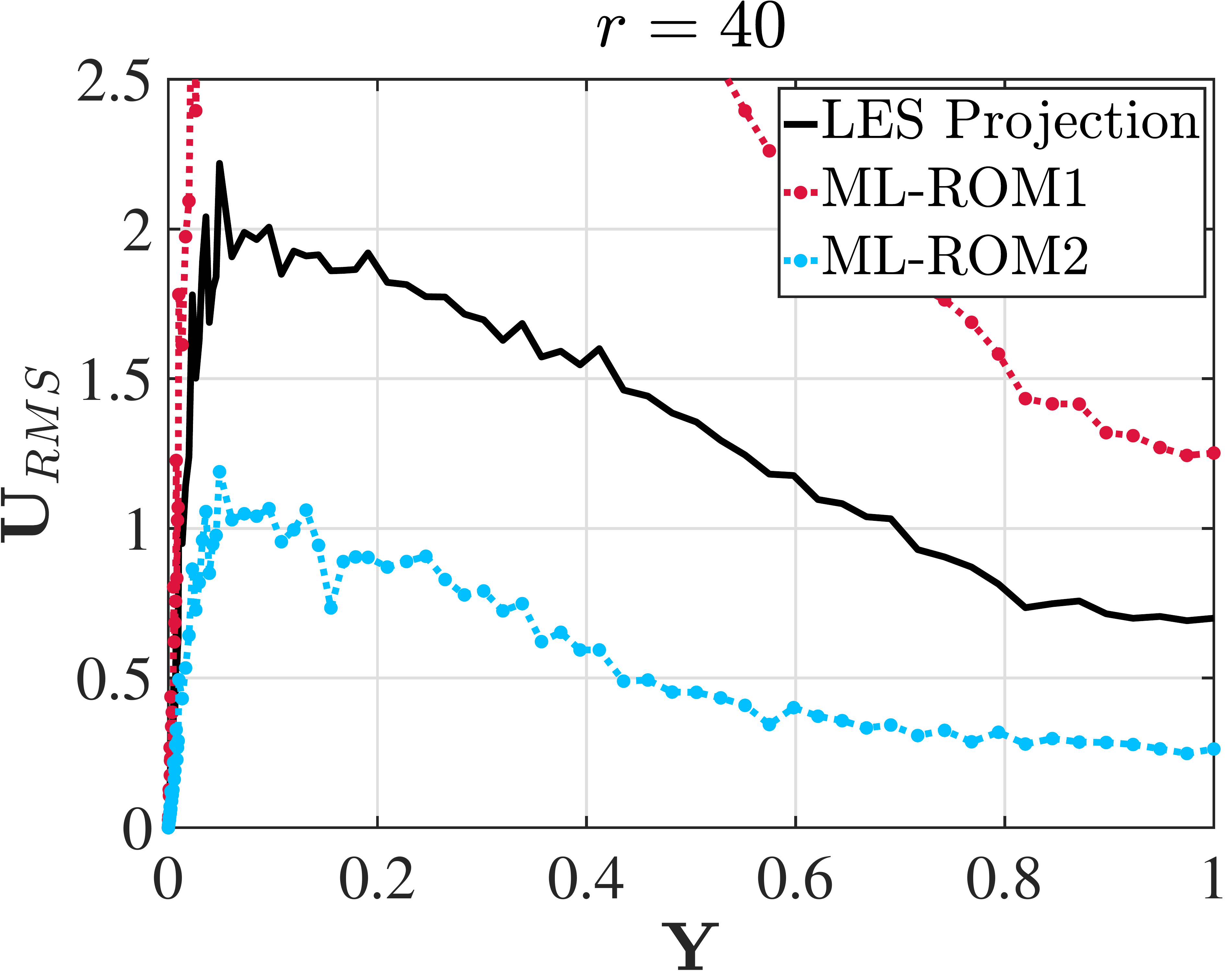}        \caption{$r=40$}
         \label{fig:stat-r-40}
     \end{subfigure} 
     \begin{subfigure}[b]{0.48\textwidth}
         \centering
     \includegraphics[width=.45\textwidth]{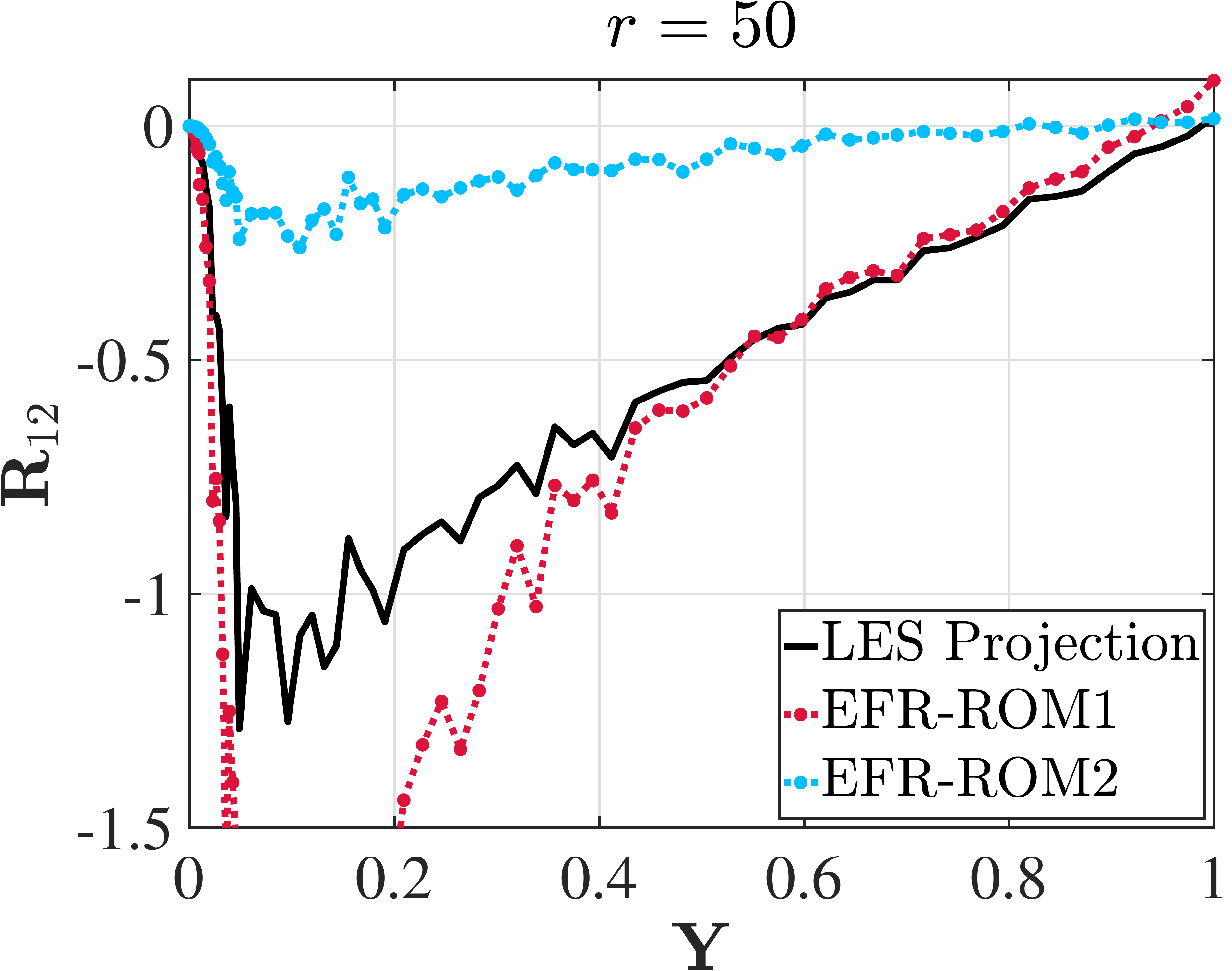}
    \includegraphics[width=.45\textwidth]{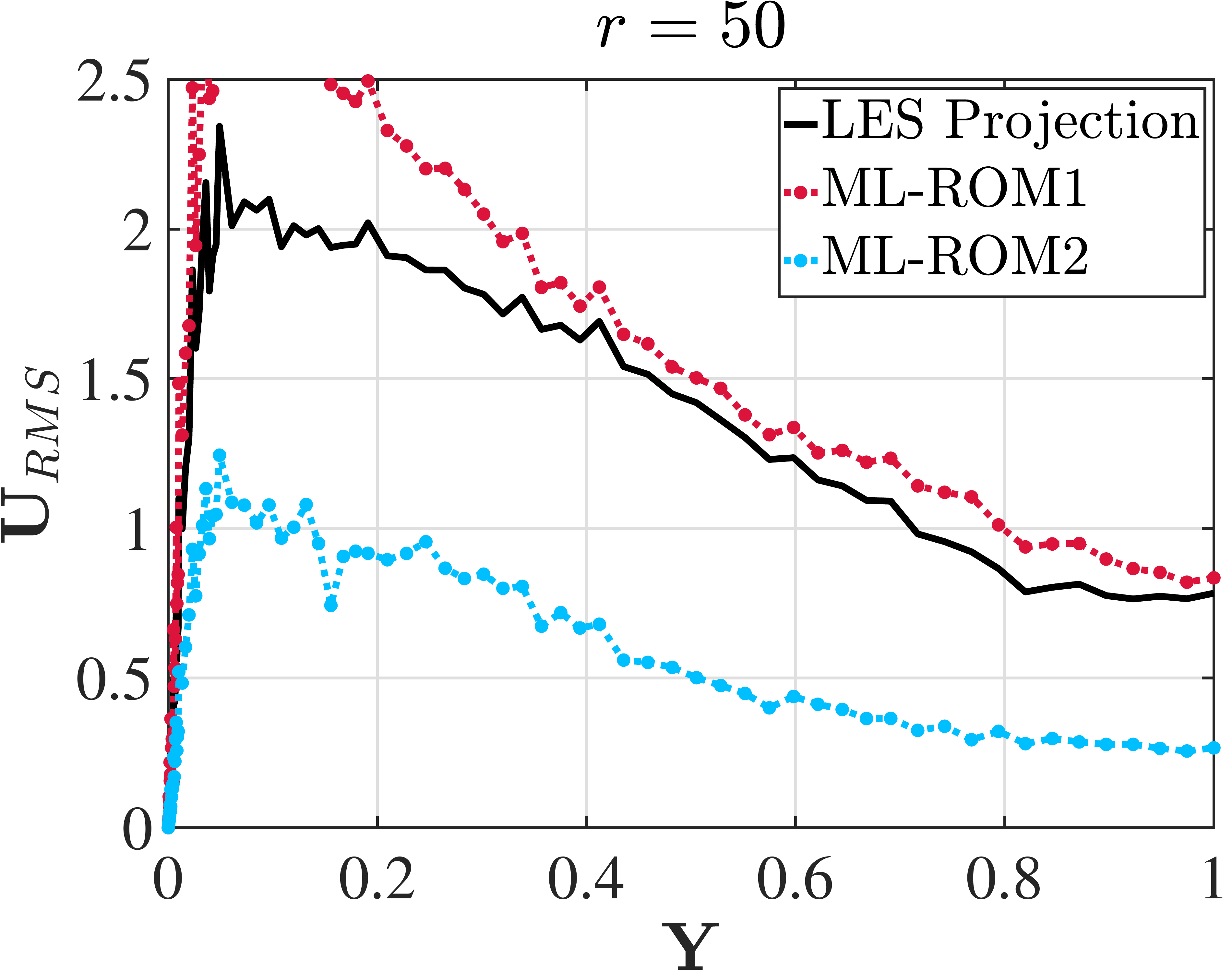}         \caption{$r=50$}
         \label{fig:stat-r-50}
     \end{subfigure} 
     \caption{
Second-order ML-ROM statistics for $\alpha=6\times 10^{-3}$}
    \label{fig:stat-alpha-large}
\end{figure}

\begin{figure}[H]
\centering
     \begin{subfigure}[b]{0.48\textwidth}
         \centering
    \includegraphics[width=.45\textwidth]{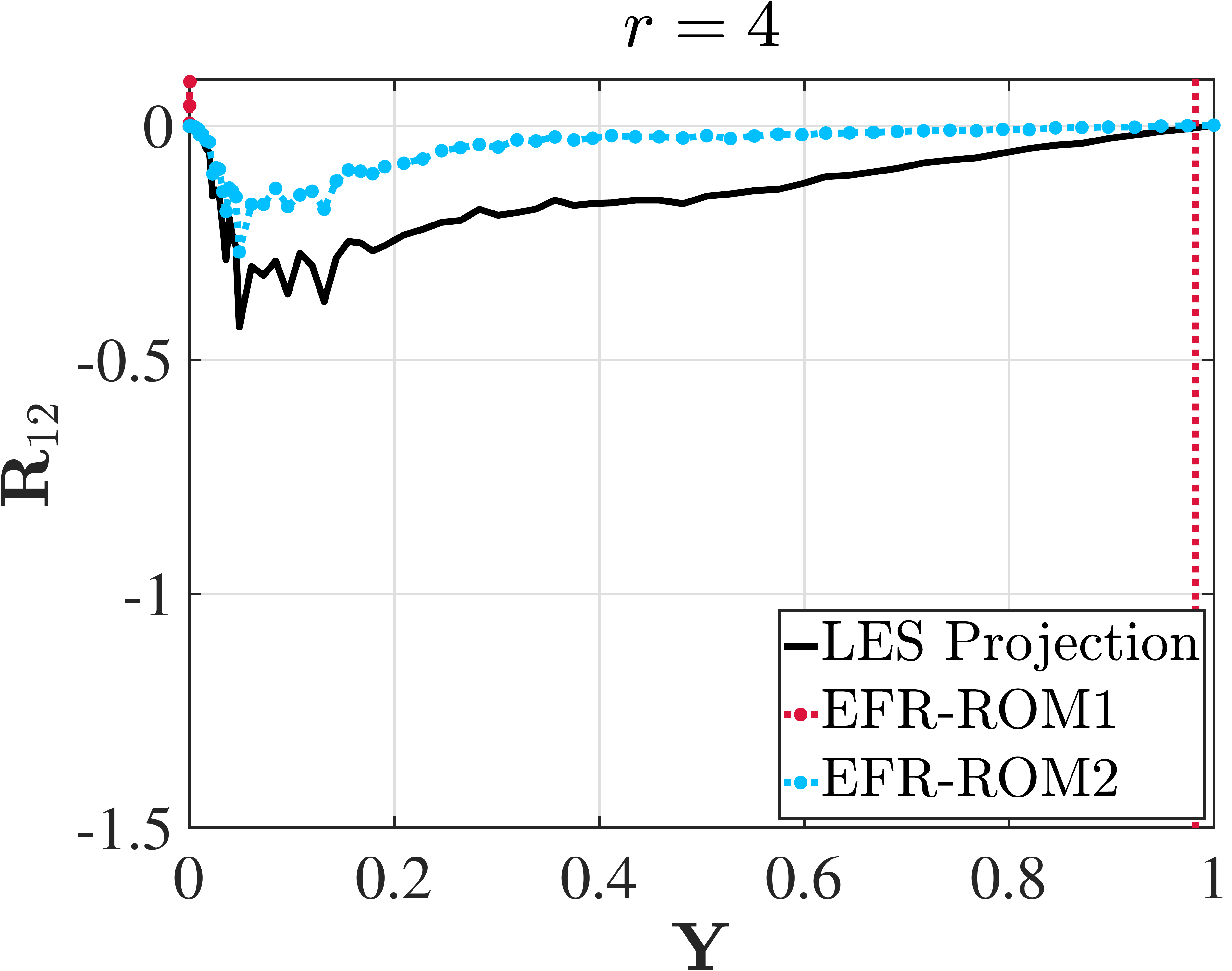}
    \includegraphics[width=.45\textwidth]{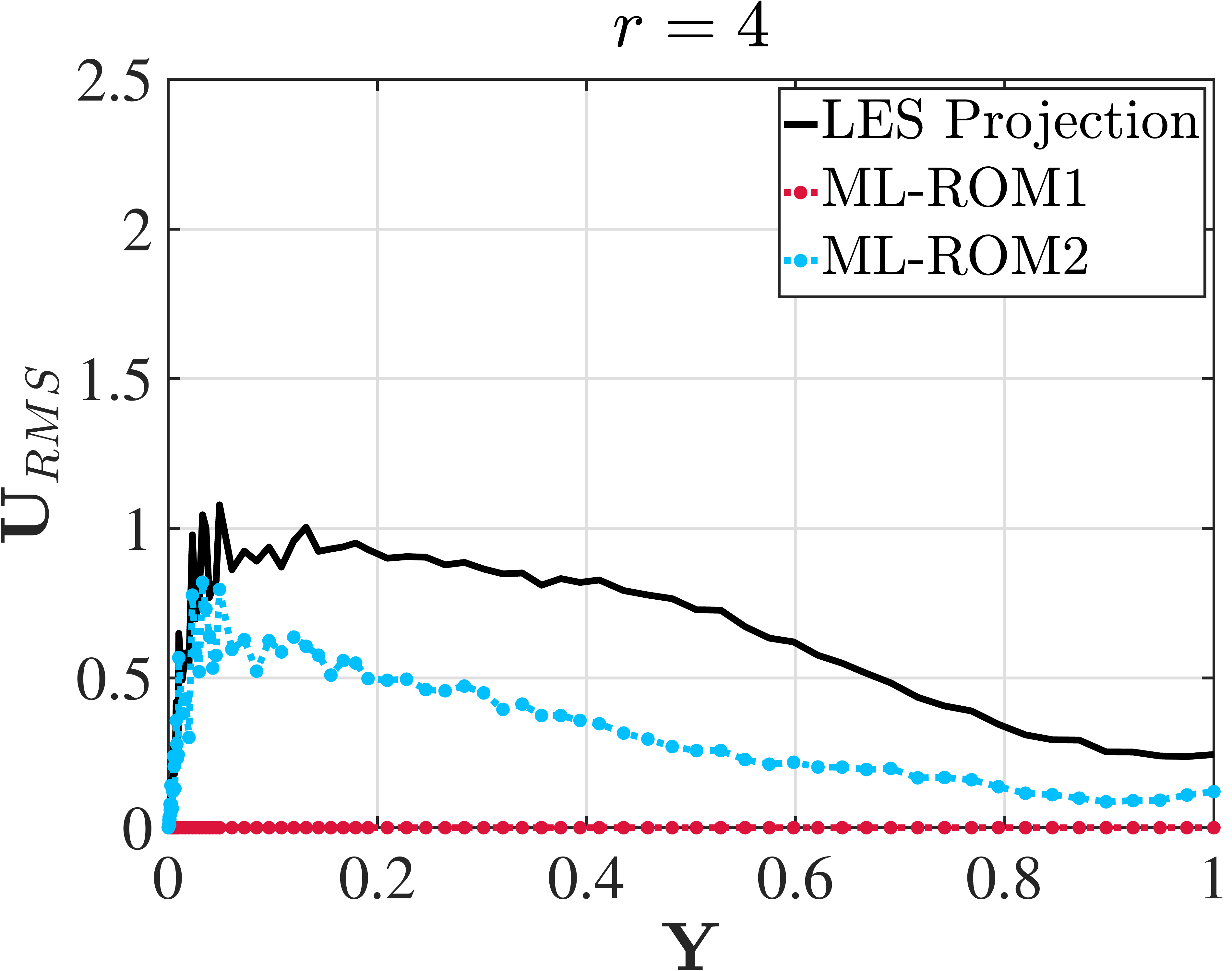}         \caption{$r=4$}
         \label{fig:stat-r-4}
     \end{subfigure}
     \begin{subfigure}[b]{0.48\textwidth}
         \centering
    \includegraphics[width=.45\textwidth]{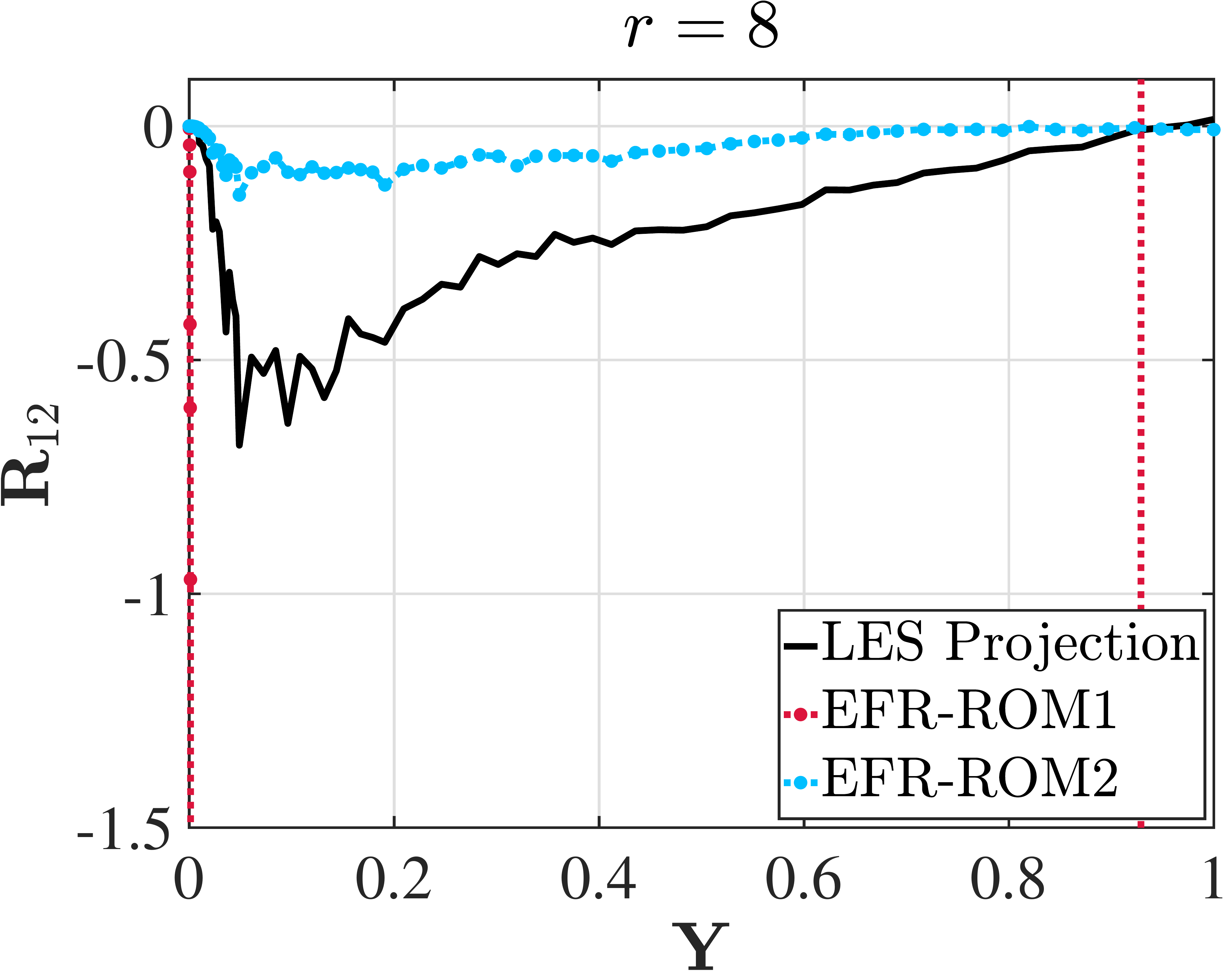}
    \includegraphics[width=.45\textwidth]{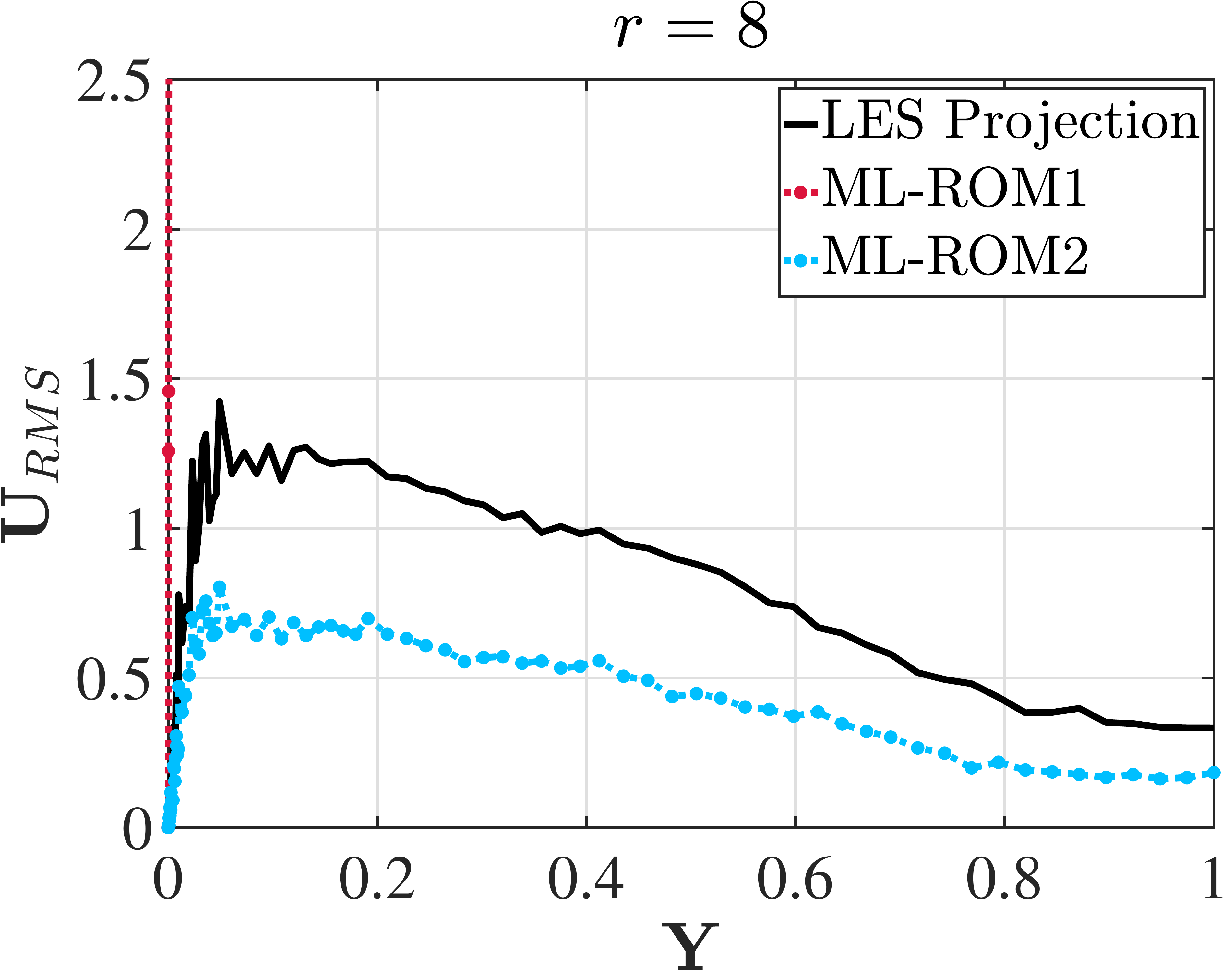}         \caption{$r=8$}
         \label{fig:stat-r-8}
     \end{subfigure}
     \begin{subfigure}[b]{0.48\textwidth}
         \centering
    \includegraphics[width=.45\textwidth]{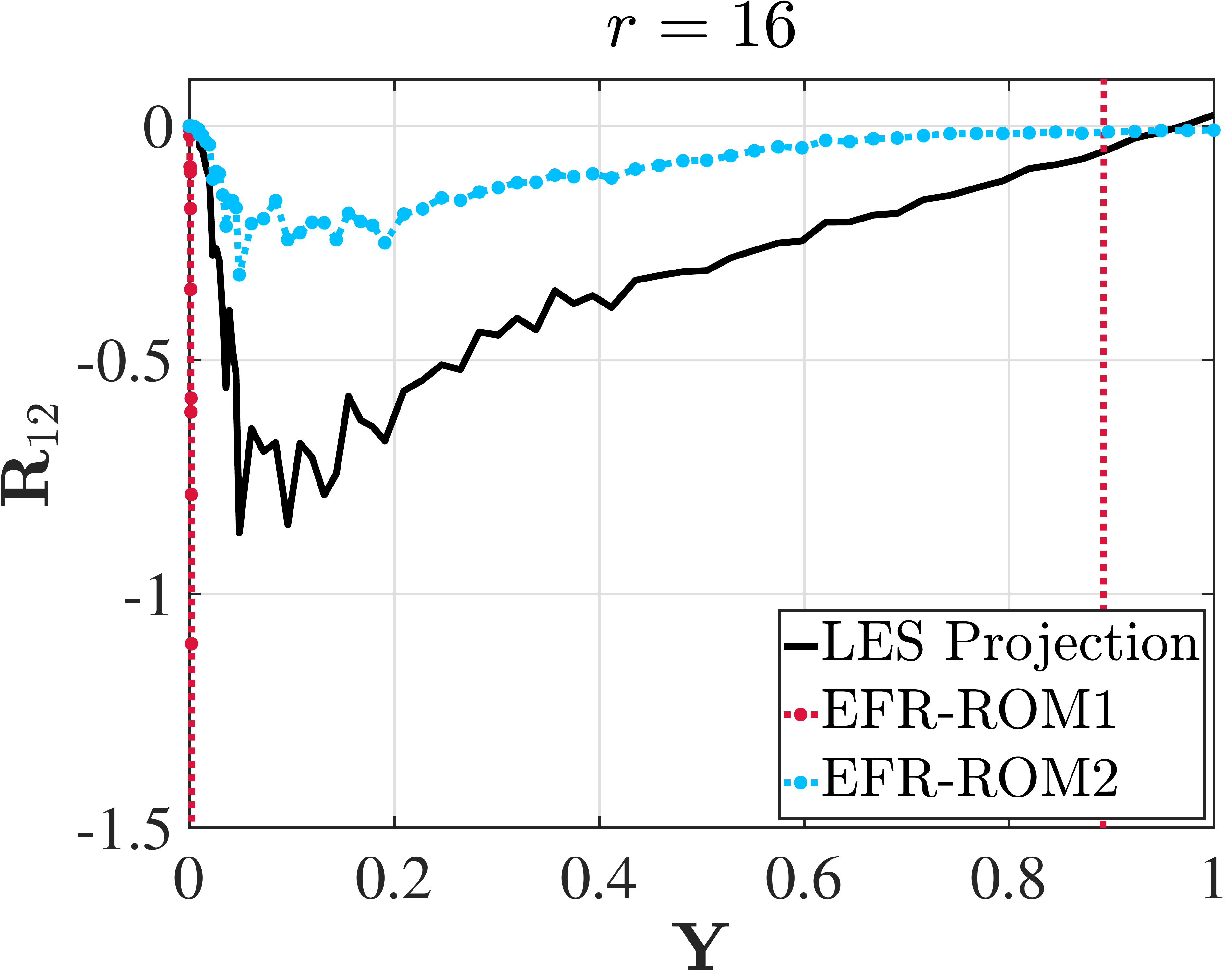}
    \includegraphics[width=.45\textwidth]{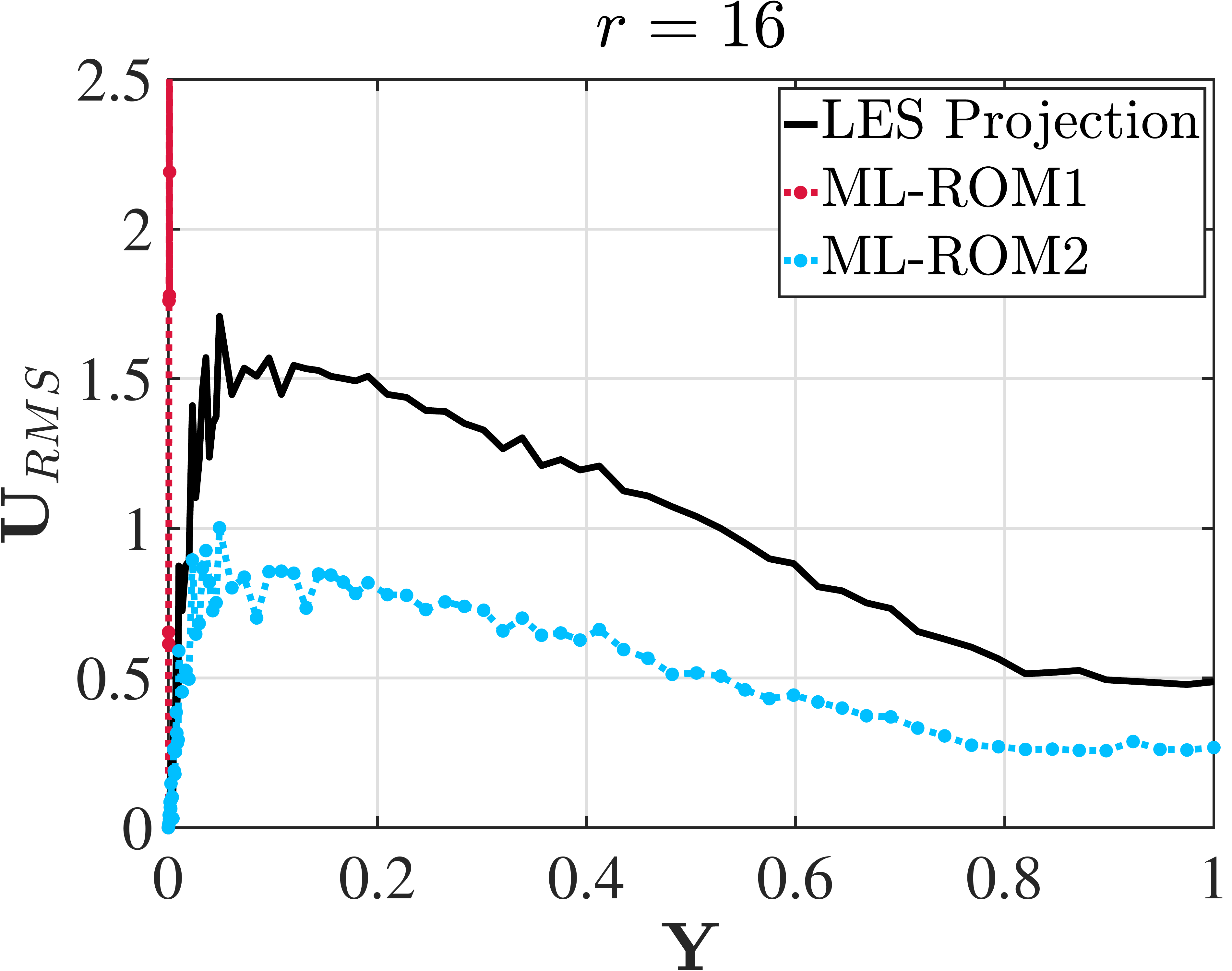}          \caption{$r=16$}
         \label{fig:stat-r-16}
     \end{subfigure}
     \begin{subfigure}[b]{0.48\textwidth}
         \centering
    \includegraphics[width=.45\textwidth]{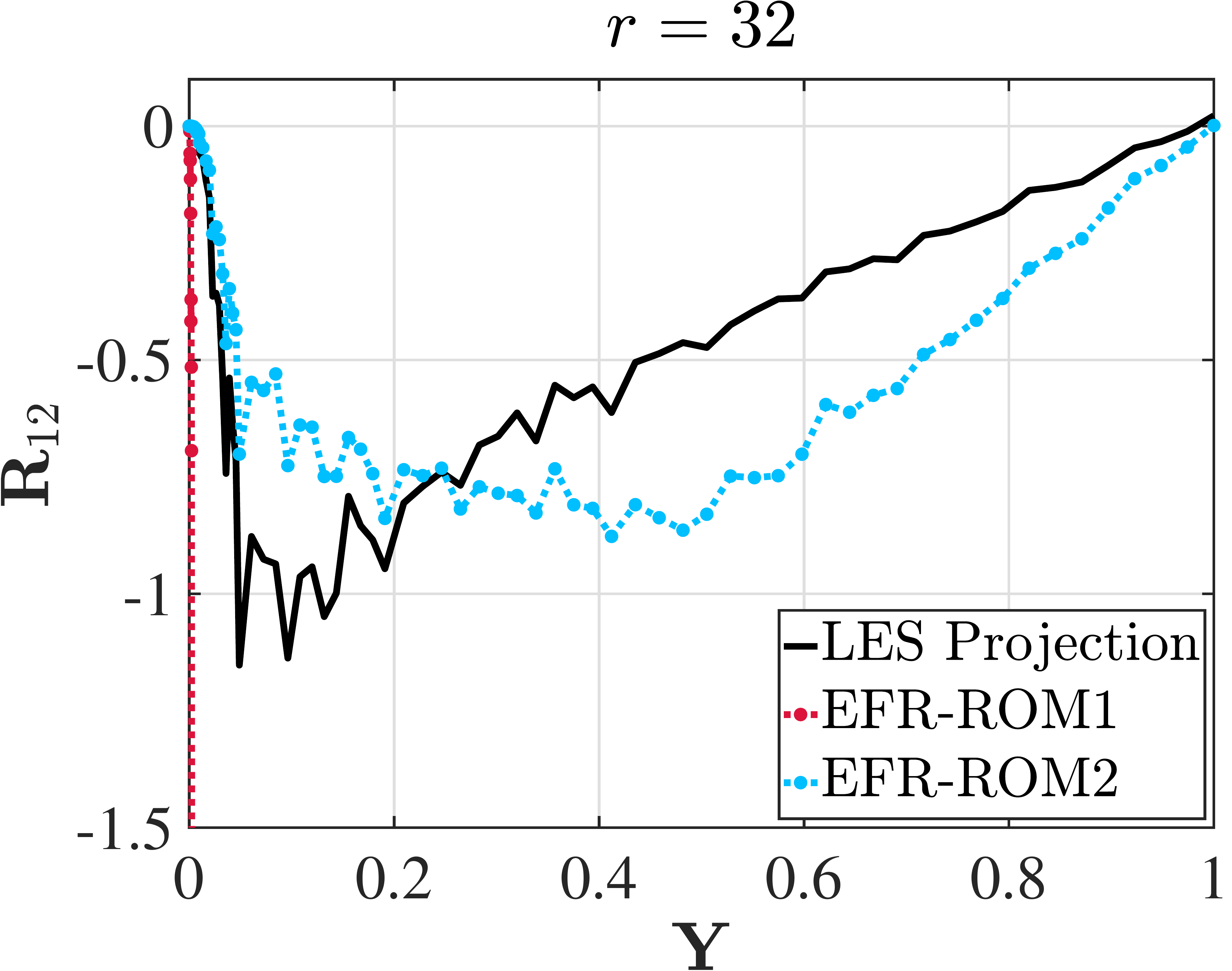}
    \includegraphics[width=.45\textwidth]{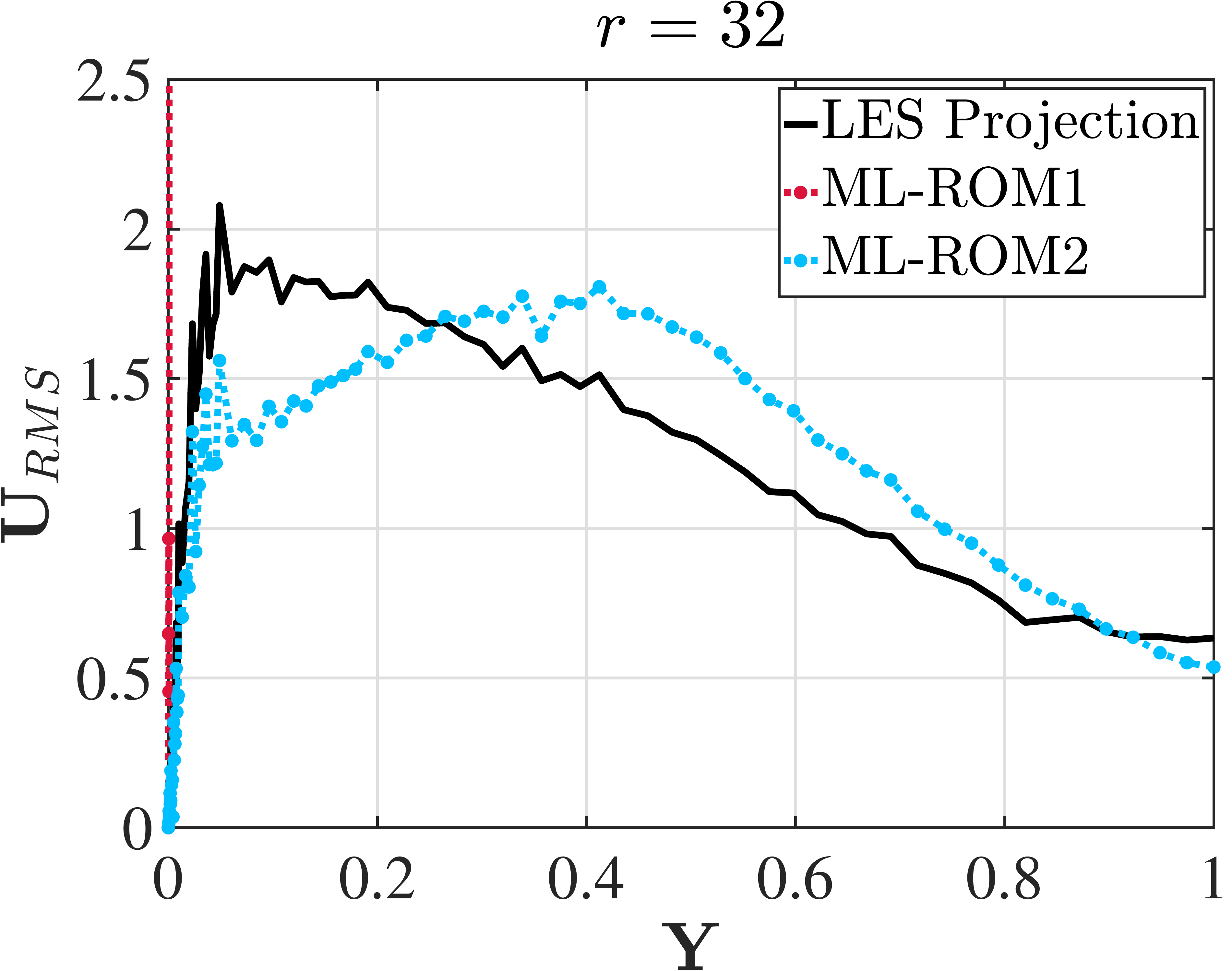}       \caption{$r=32$}
         \label{fig:stat-r-32}
    \end{subfigure}     
     \begin{subfigure}[b]{0.48\textwidth}
         \centering
    \includegraphics[width=.45\textwidth]{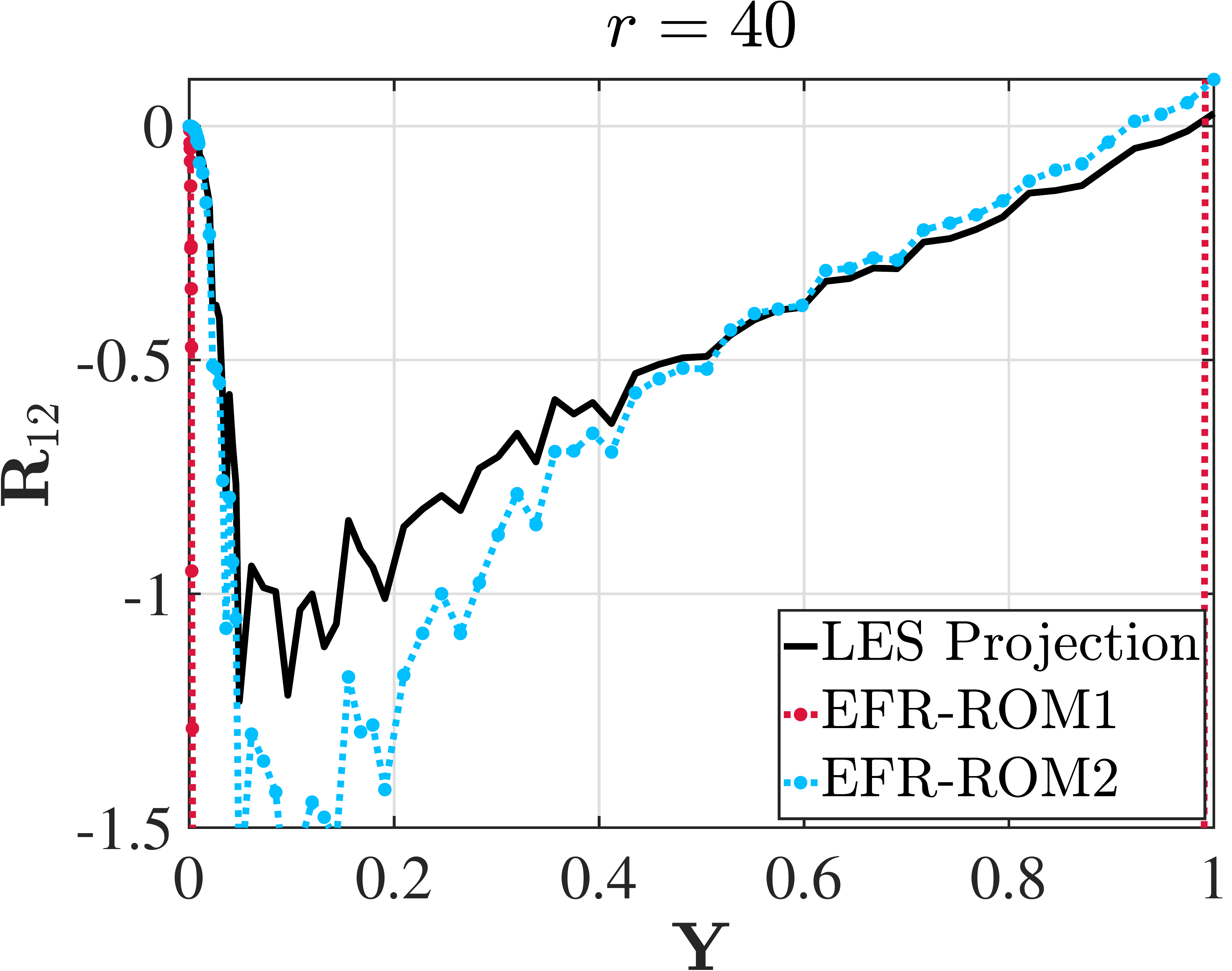}
    \includegraphics[width=.45\textwidth]{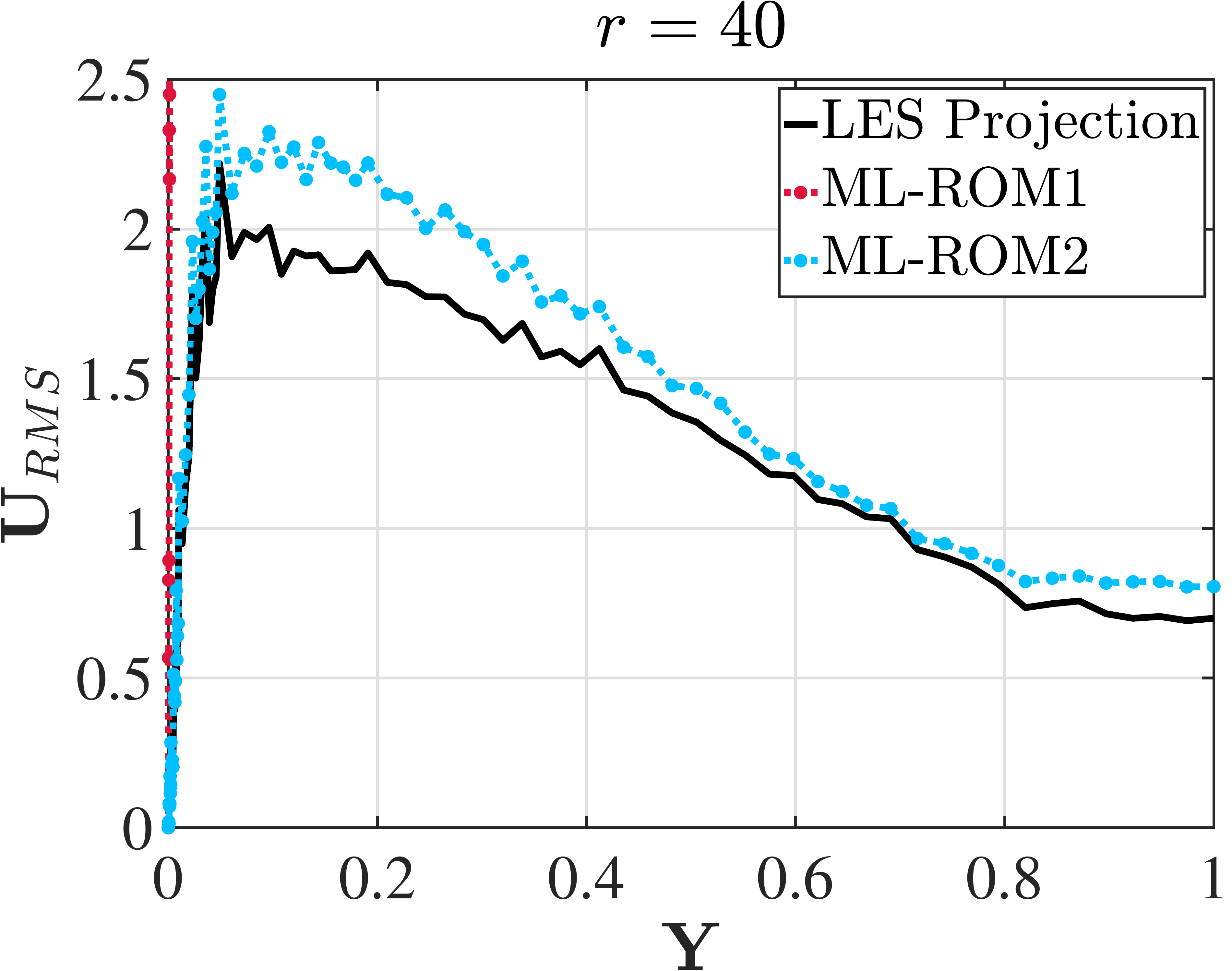}        \caption{$r=40$}
         \label{fig:stat-r-40}
     \end{subfigure} 
     \begin{subfigure}[b]{0.48\textwidth}
         \centering
     \includegraphics[width=.45\textwidth]{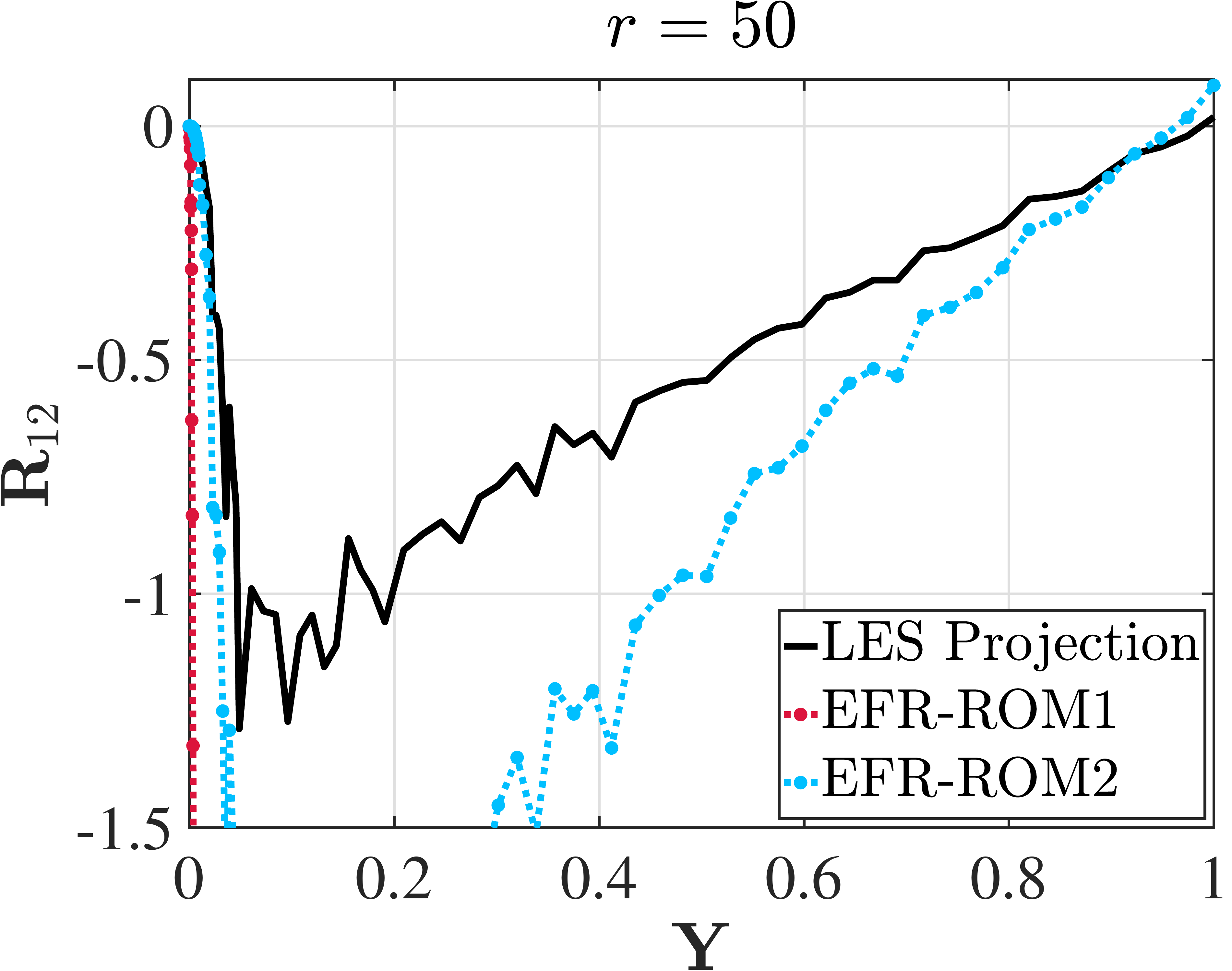}
    \includegraphics[width=.45\textwidth]{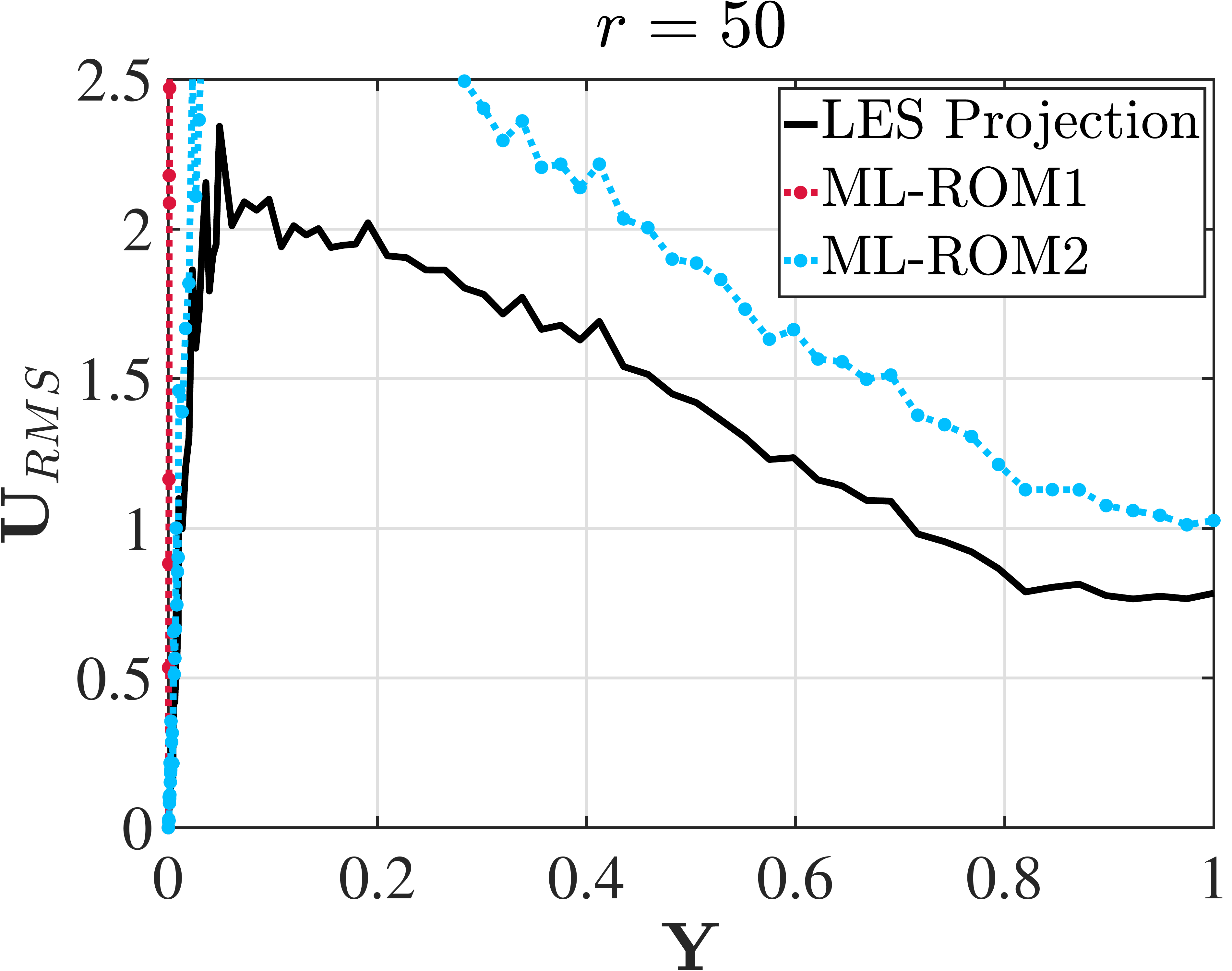}         \caption{$r=50$}
         \label{fig:stat-r-50}
     \end{subfigure} 
     \caption{
Second-order ML-ROM statistics for $\alpha=6\times 10^{-4}$}
    \label{fig:stat-alpha-small}
\end{figure}

In this section, we investigate the role played by the two ROM lengthscales, $\delta_1$ and $\delta_2$, in the ML-ROM
~\eqref{eqn:ml-v2}: 
\begin{eqnarray}
	\btau
	= - \bigl( \alpha \, U_{ML} \, L_{ML} \bigr) \, S_r \, \ba \, .
	\label{eqn:ml-v3}
\end{eqnarray}
We denote ML-ROM1 the ML-ROM in which $L_{ML}=\delta_1$ in~\eqref{eqn:ml-v3} and ML-ROM2 the ML-ROM in which $L_{ML}=\delta_2$ in~\eqref{eqn:ml-v3}.
To ensure a fair comparison of ML-ROM1 and ML-ROM2, we use the same constant $\alpha$ and the same velocity scale $U_{ML}$ (i.e., the time averaged streamwise velocity component) 
in~\eqref{eqn:ml-v3} for both models, and vary only the ROM lengthscale, i.e., $L_{ML}=\delta_1$ or $L_{ML}=\delta_2$.
To vary the ROM lengthscale, we vary the $r$ value in the definitions of $\delta_1$ and $\delta_2$.
In Figures \ref{fig:ke-alpha-1}--\ref{fig:stat-alpha-small} 
we plot the time evolution of the kinetic energy and the second-order statistics of the ML-ROM1 and ML-ROM2 for different $r$ values and two different $\alpha$ values:
$\alpha = 6 \times 10^{-3}$ (Figures~\ref{fig:ke-alpha-1} and \ref{fig:stat-alpha-large}) and 
$\alpha = 6 \times 10^{-4}$ (Figures~\ref{fig:ke-alpha-3}and \ref{fig:stat-alpha-small}).
(
Results for more $\alpha$ values are presented in the preliminary numerical investigation in~\cite{mou2022numerical}.
)
As a benchmark for the ROM results, we use the projection of the the FOM results on the ROM basis (denoted as LES-proj in these plots).

\paragraph{Stability}
Since the results in Table~\ref{table:delta} show that $\delta_2$ is between one and two orders of magnitude higher than $\delta_1$, we expect ML-ROM2 to yield more stable results than ML-ROM1.
Indeed, since we fix all the ML-ROM parameters (i.e., the constant $\alpha$ and the velocity scale $U_{ML}$) and $\delta_2$ is significantly larger than $\delta_1$, we expect the ML-ROM2 artificial viscosity to be higher than the ML-ROM1 artificial viscosity (and, thus, ML-ROM2 to be more stable than ML-ROM1). 
This is clearly shown in the plots in Figures \ref{fig:ke-alpha-1}--\ref{fig:stat-alpha-small}, in which ML-ROM2 yields stable results for {\it all} $r$ values and for {\it both} $\alpha$ values.
In contrast, for the smallest $\alpha$ value, $\alpha = 6 \times 10^{-4}$ (Figures~\ref{fig:ke-alpha-3} and \ref{fig:stat-alpha-small}), ML-ROM1 blows up for {\it all} $r$ values. 
Furthermore, for the largest $\alpha$ value, $\alpha = 6 \times 10^{-3}$ (Figures~\ref{fig:ke-alpha-1} and \ref{fig:stat-alpha-large}), ML-ROM1 blows up for the small $r$ values (i.e., $r = 4, 8, 16$, and $32$). 
To quantify the stability of the two ML-ROMs, 
in Table~\ref{tab:alpha-value-threshold},  
for different $r$ values, we list the threshold $\alpha_0$ value, i.e., the value that ensures that, if $\alpha > \alpha_0$, then the ML-ROM is stable.
These results show that, for each $r$ value, the threshold $\alpha_0$ value is 
more than an order of magnitude lower for ML-ROM2 than for ML-ROM1.
Thus, we conclude that ML-ROM2 is more stable than ML-ROM1, 
which is the same conclusion 
as that yielded by Figures~\ref{fig:ke-alpha-1}--\ref{fig:stat-alpha-small}.

\begin{table}[H]
    \centering
    \begin{tabular}{c c|c c c c c c c c c c c c }
    \hline\hline
        &$r$ & 4 & 8 & 16 & 32 &40  & 50 
    \\ \hline
        ML-ROM1  &$\alpha_0$ & 10e-3&  9.8e-3& 9.2e-3& 8.5e-3& 6.5e-3& 6.2e-3
    \\ 
        ML-ROM2   &$\alpha_0$
        &2.9e-4 & 3.4e-4 &4.1e-4 & 6.7e-4 &5.9e-4 &7.5e-4        
    \\ \hline
    \end{tabular}
    \caption{ML-ROM threshold $\alpha_0$ values for different $r$ values.}
    \label{tab:alpha-value-threshold}
\end{table}

\paragraph{Accuracy}
Since $\delta_2$ is between one and two orders of magnitude higher than $\delta_1$, we expect the ML-ROM accuracy to depend on the constant $\alpha$. 
Indeed, 
the ML-ROM1 and ML-ROM2 plots in Figures \ref{fig:ke-alpha-1}--\ref{fig:stat-alpha-small} do not show a clear winner:
For the largest $\alpha$ value (i.e., $\alpha = 6 \times 10^{-3}$), ML-ROM1 is more accurate than ML-ROM2 for large $r$ values (i.e., $r = 40$ and $50$) and less accurate for small $r$ values (i.e., $r=4, 8$, and $16$).
For the smallest $\alpha$ value (i.e., $\alpha = 6 \times 10^{-4}$), ML-ROM2 is more accurate than ML-ROM1 for {\it all} $r$ values (since ML-ROM1 simply blows up).
We also note that, for $\alpha = 6 \times 10^{-4}$, ML-ROM2 is quite accurate for larger $r$ values (i.e., $r = 32, 40$, and $50$).

\paragraph{Parameter Sensitivity}
Finally, we investigate the ML-ROM's parameter sensitivity.
Specifically, we investigate which ROM lengthscale yields ML-ROMs that are less sensitive (i.e., more robust) with respect to the ML-ROM's $\alpha$ parameter.
We emphasize that the model sensitivity with respect to model parameters is a well known issue that has 
hindered the development of closures and stabilizations in CFD over the years~\cite{sagaut2006large,BIL05,layton2012approximate}.
Thus, finding robust (i.e., less sensitive) ML-ROMs that require as little parameter tuning as possible is an important practical problem.

To quantify the ML-ROM's parameter sensitivity, 
in Table~\ref{tab:alpha-value-optimal}, for different $r$ values, we list the optimal $\alpha$ value in ML-ROM, i.e., the $\alpha$ value that ensures that the average ROM kinetic energy, $KE^{ROM}$,  
is the closest to the average FOM (LES) kinetic energy, $KE^{LES}$. 
Specifically, we 
solve the following optimization problem:
\begin{align}
    \min_{\alpha} \left| 
    \overline{KE}^{ROM} 
    - 
    \overline{KE}^{LES} 
    \right|,
    \label{eqn:min-alpha}
\end{align}
where
$ \overline{KE} = \frac{1}{M} \sum_{k=1}^M KE(t_k) $, 
and $M$ is the number of snapshots.

\begin{table}[H]
    \centering
    \begin{tabular}{c c|c c c c c c c c c c c c }
    \hline\hline
        &$r$ & 4 & 8 & 16 & 32 &40  & 50 
    \\ \hline
        ML-ROM1  &$\alpha$ 
        &1.50e-2 &1.38e-2 &1.38e-2 &1.28e-2
&7.80e-3 &6.51e-3
    \\ 
        ML-ROM2   &$\alpha$
        & 4.35e-4 &4.42e-4 &5.33e-4 &6.70e-4 &6.20e-4 &7.50e-4
    \\ \hline
    \end{tabular}
    \caption{ML-ROM optimal $\alpha$ values for different $r$ values.
    }
    \label{tab:alpha-value-optimal}
\end{table}

The results in Table~\ref{tab:alpha-value-optimal} show that the optimal ML-ROM $\alpha$ values are sensitive with respect to changes in $r$.
Indeed, as $r$ increases from $4$ to $50$, $\alpha$ decreases by a factor of $2.3$ for ML-ROM1 and increases by a factor of $1.7$ for ML-ROM2.
Overall, the results in Table~\ref{tab:alpha-value-optimal} 
show that the ML-ROM's parameter $\alpha$ sensitivity is 
higher for ML-ROM1 than for ML-ROM2.

 \subsection{Numerical Results: EFR-ROM Investigation}

To further investigate the role played by the two ROM lengthscales, $\delta_1$ and $\delta_2$, in reduced order modeling of turbulent flows, in this section we consider the EFR-ROM presented in Section~\ref{sec:efr-rom}.
We emphasize that the EFR-ROM is completely different from the ML-ROM investigated in Section~\ref{sec:numerical-results-ml}:
The EFR-ROM is a Reg-ROM based on numerical stabilization, whereas the ML-ROM is a closure model.
Thus, the EFR-ROM investigation in this section could shed new light on the novel ROM lengthscale, $\delta_2$.

To ensure a fair comparison of the two ROM lengthscales, $\delta_1$ and $\delta_2$, we fix the parameters $\gamma$ and $\chi$ in the EFR-ROM in Section~\ref{sec:efr-rom} and change only the filter radius, $\delta$, of the differential filter used in Step (II) of the EFR-ROM algorithm.
We denote the resulting models as EFR-ROM1 (when $\delta  = \delta_1$) and EFR-ROM2 (when $\delta  = \delta_2$).
In our numerical investigation, we fix the $\chi$ value for both EFR-ROM1 and EFR-ROM2 to $\chi = 6 \times 10^{-3}$, which is the minimum value that yields a stable solution for $r=4$ EFR-ROM1. 
We also fix the $\gamma$ value for both EFR-ROM1 and EFR-ROM2.
In our numerical investigation, we consider two $\gamma$ values.
To vary the ROM lengthscale, we vary the $r$ value in the definitions of $\delta_1$ and $\delta_2$.

\begin{figure}[H]
\centering
    \includegraphics[width=.45\textwidth]{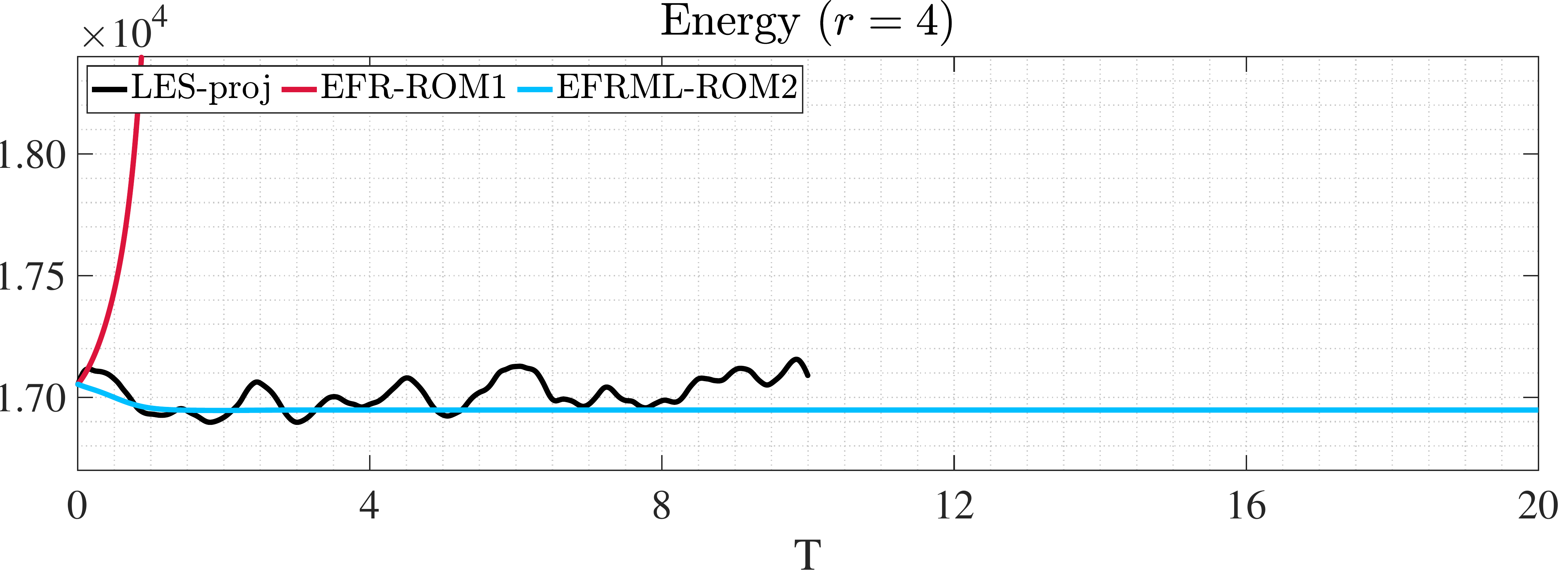}
    \includegraphics[width=.45\textwidth]{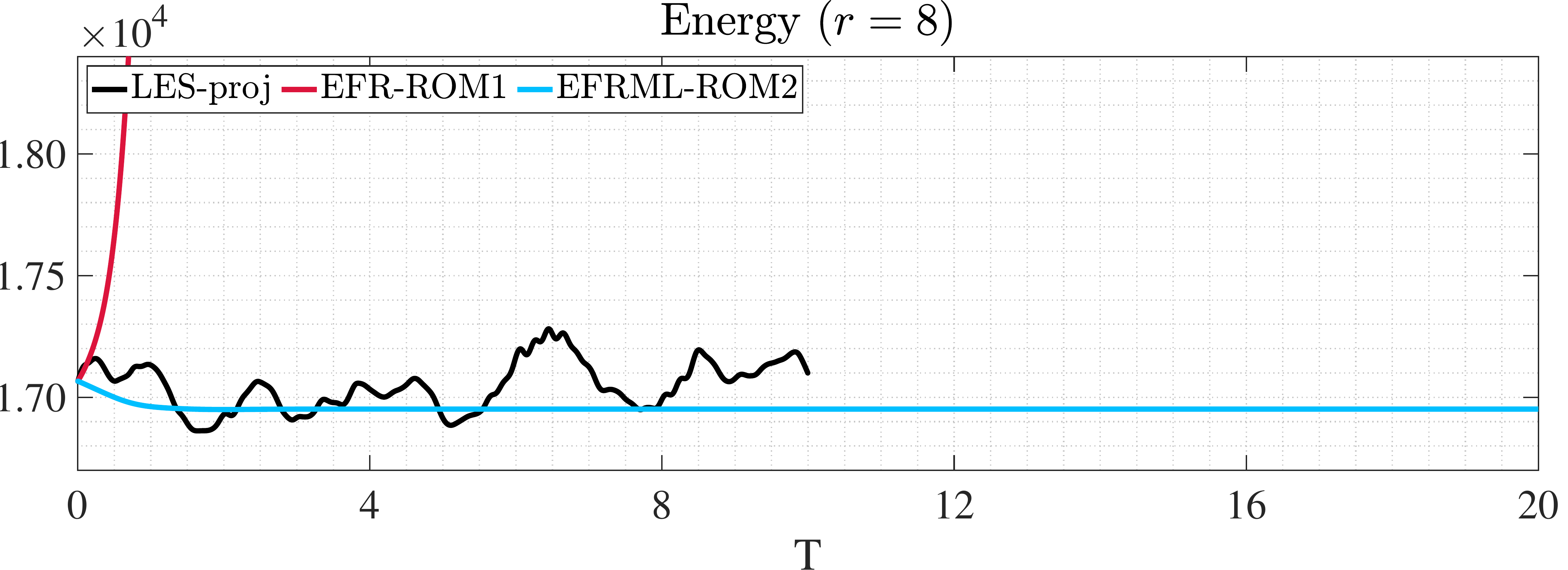}
    \includegraphics[width=.45\textwidth]{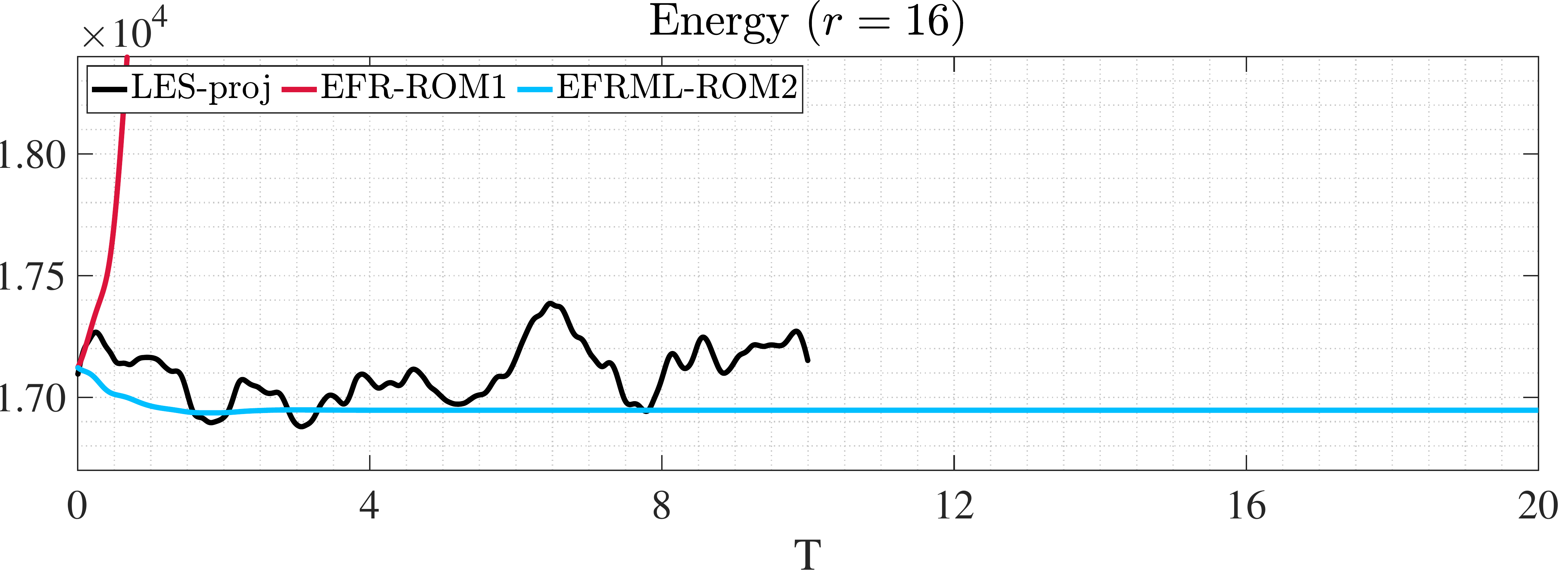}
    \includegraphics[width=.45\textwidth]{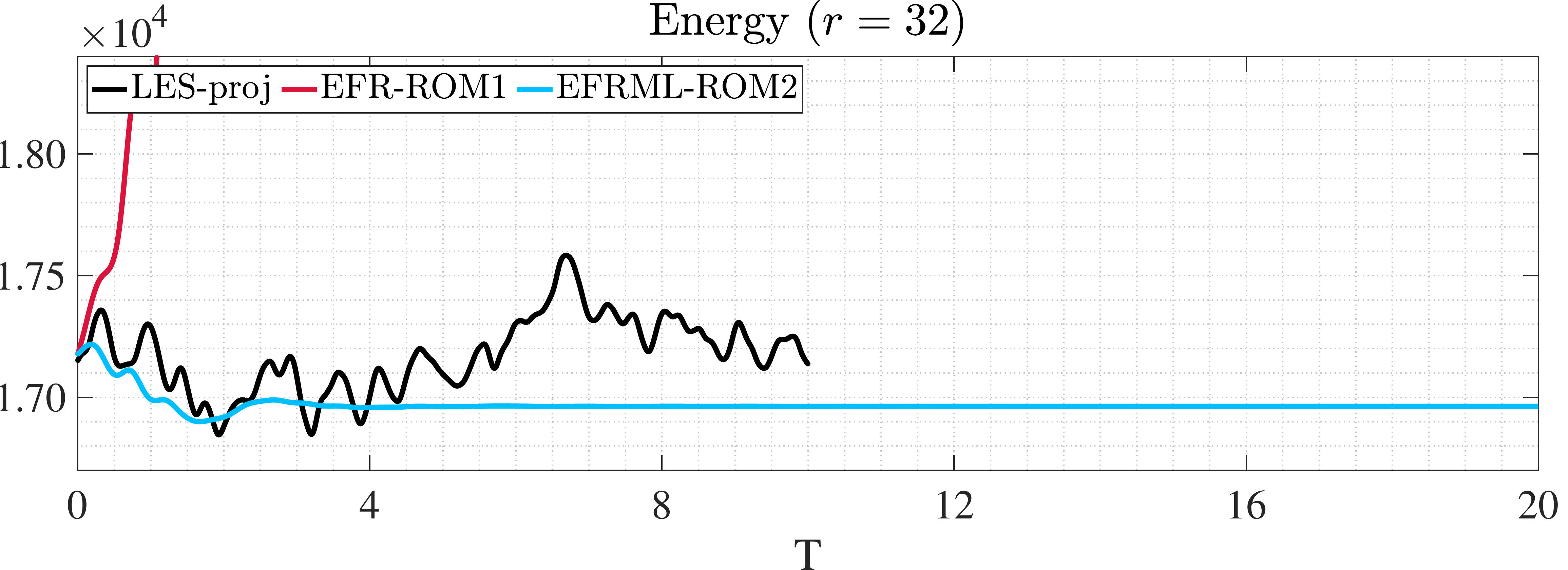}
    \includegraphics[width=.45\textwidth]{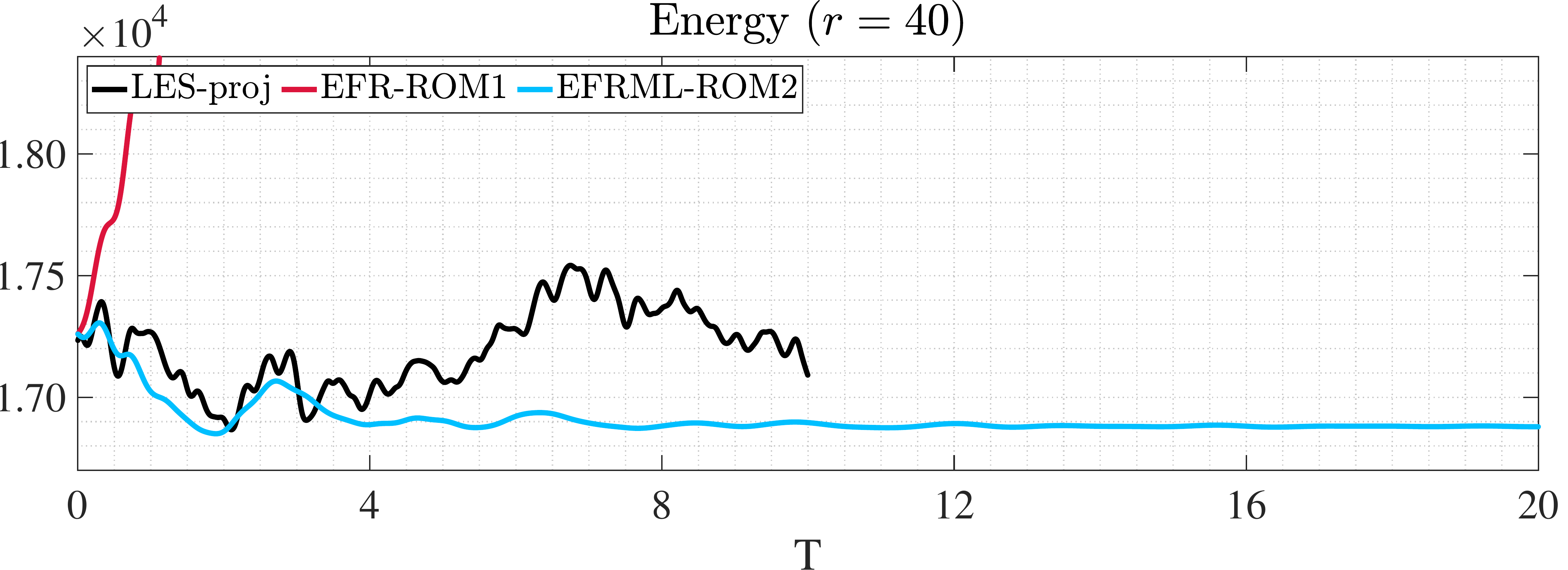}
    \includegraphics[width=.45\textwidth]{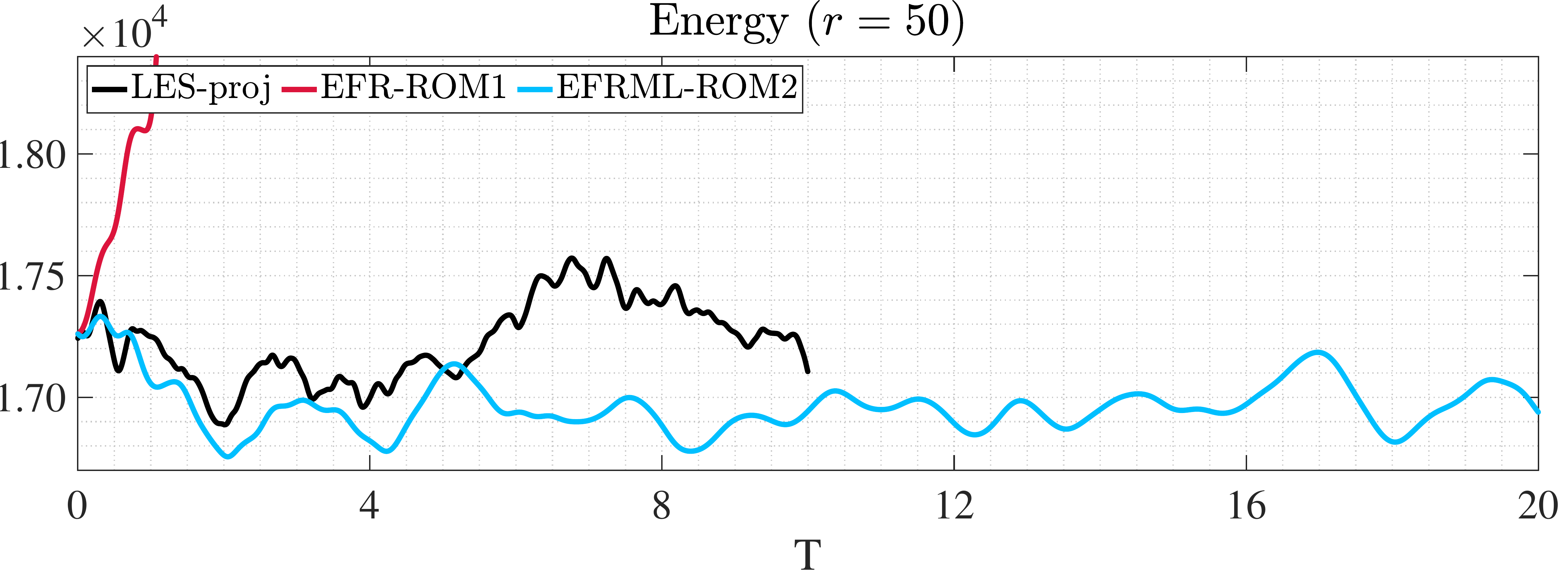}
    \caption{
    Time evolution of the EFR-ROM kinetic energy for $\gamma=8\times 10^{-2}$
    }    
    \label{fig:ke-gamma-1-efr}
\end{figure}

\begin{figure}[H]
\centering
    \includegraphics[width=.45\textwidth]{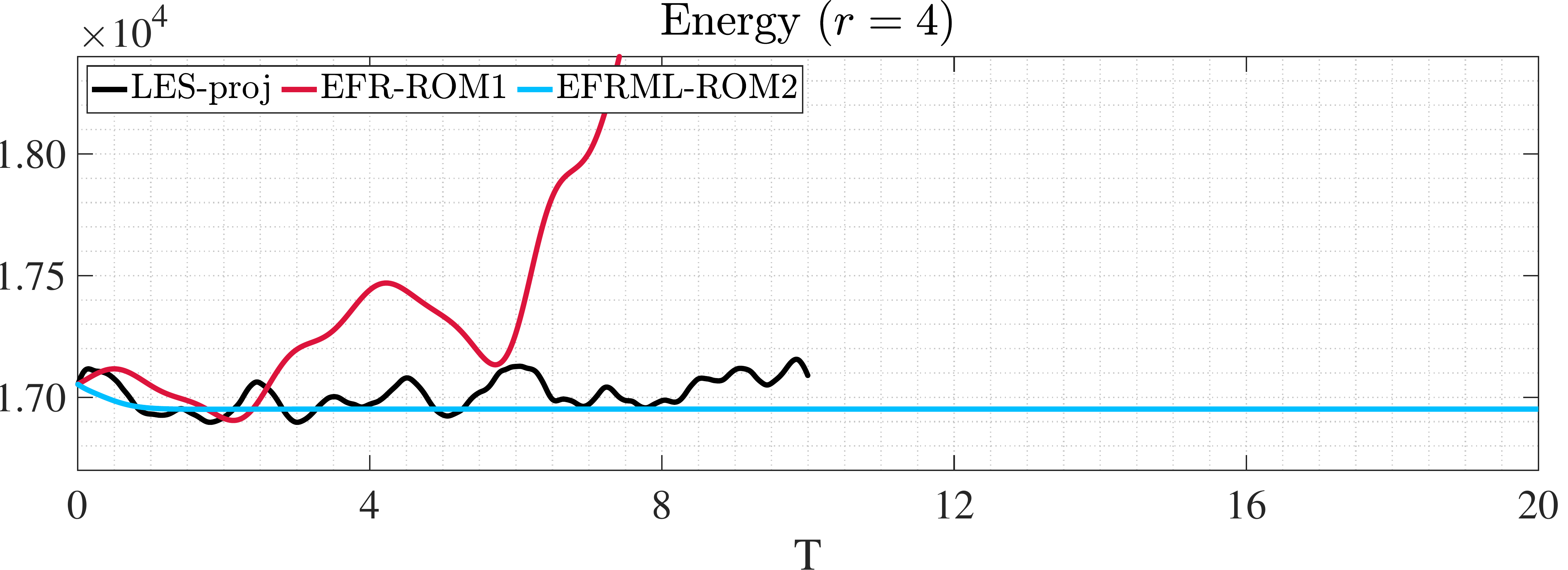}
    \includegraphics[width=.45\textwidth]{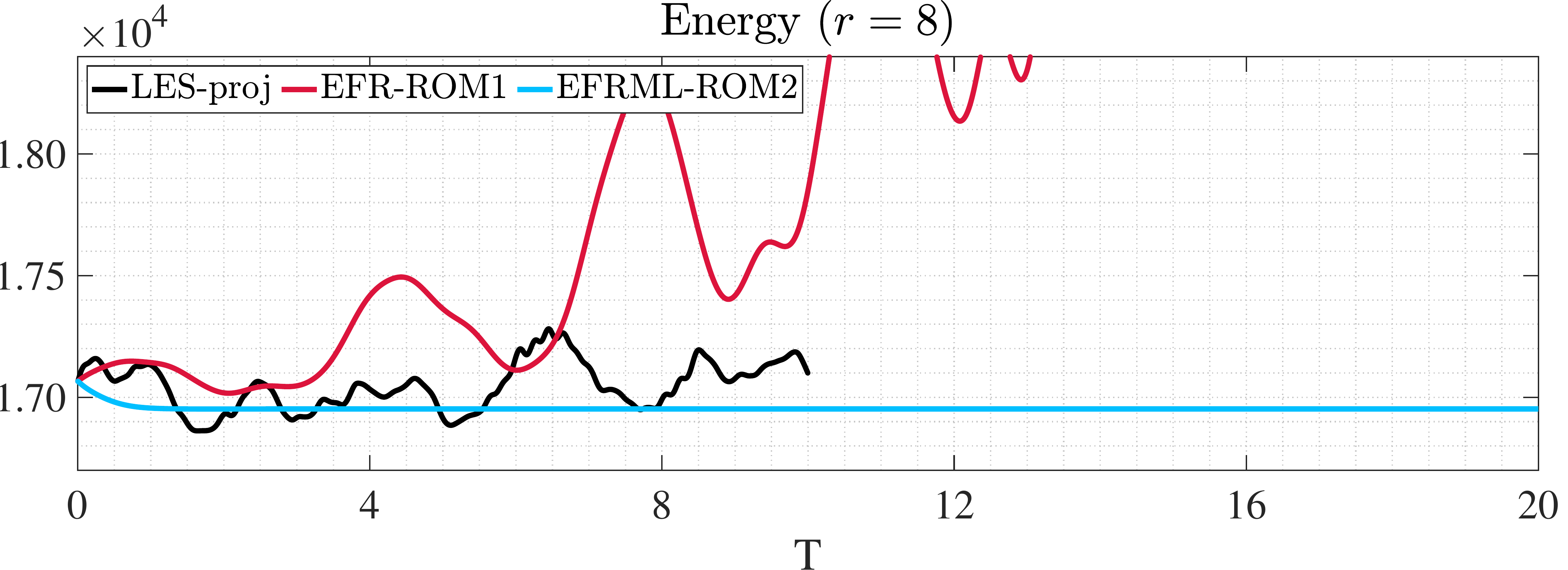}
    \includegraphics[width=.45\textwidth]{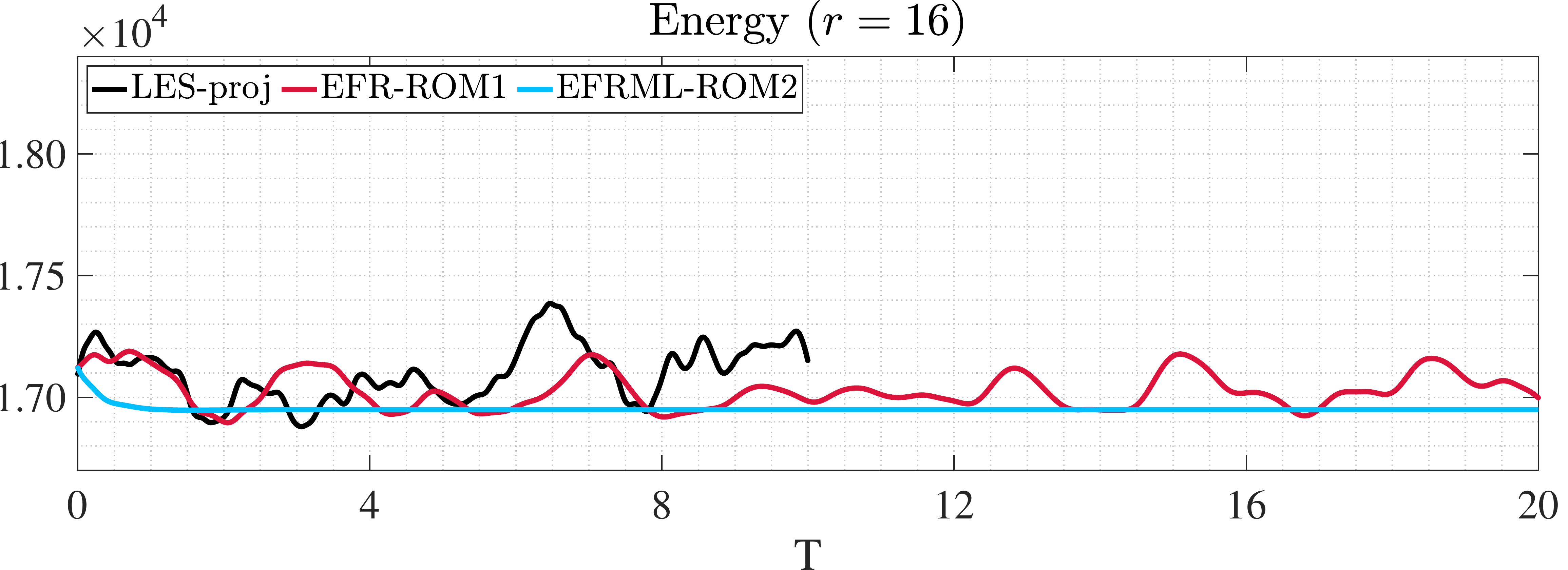}
    \includegraphics[width=.45\textwidth]{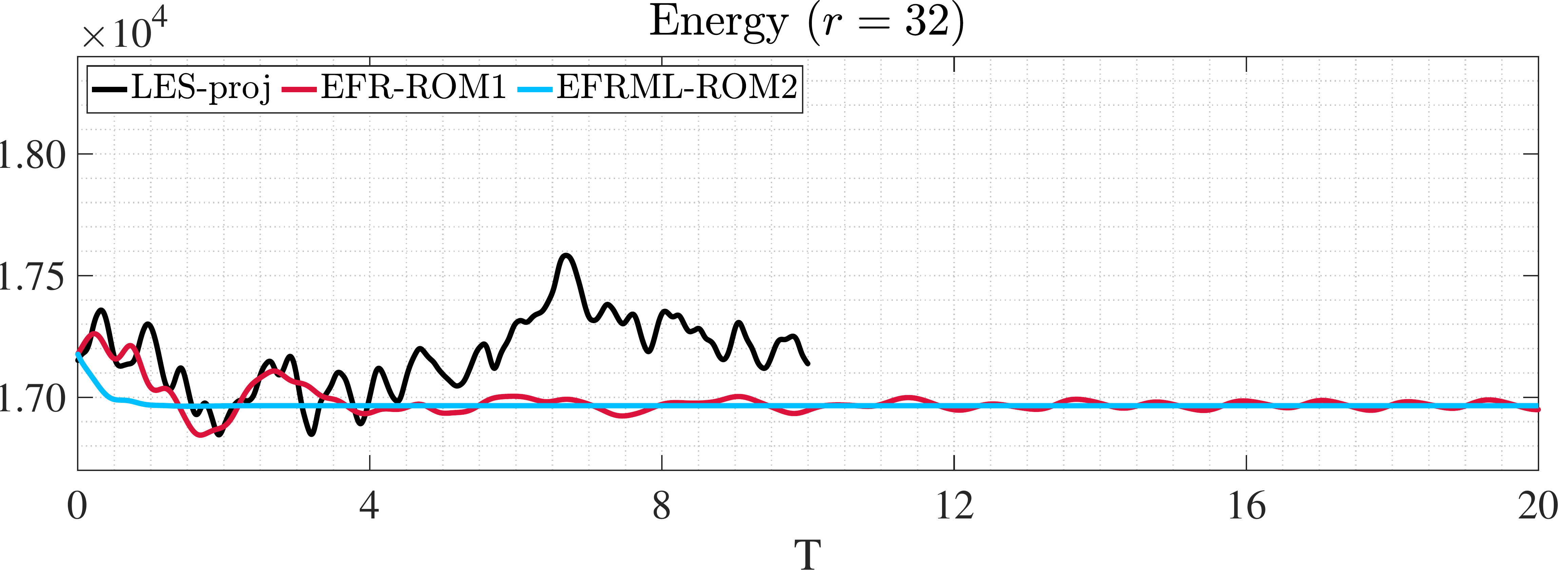}
    \includegraphics[width=.45\textwidth]{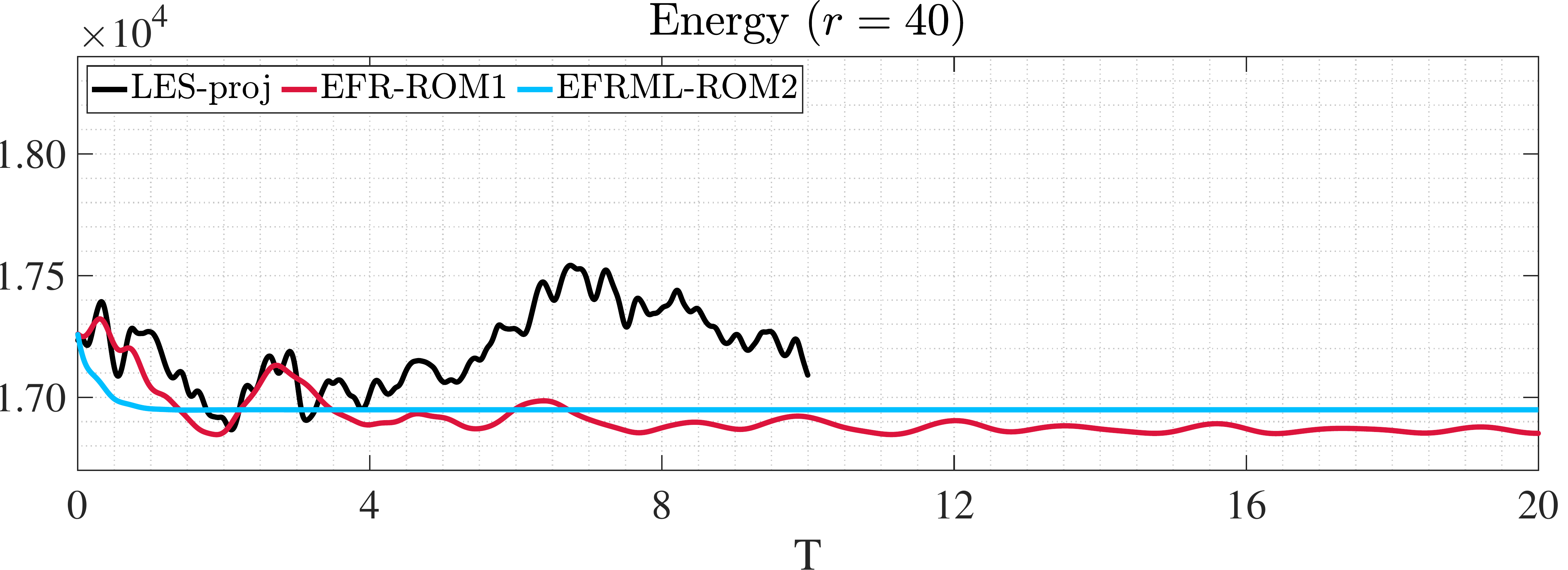}
    \includegraphics[width=.45\textwidth]{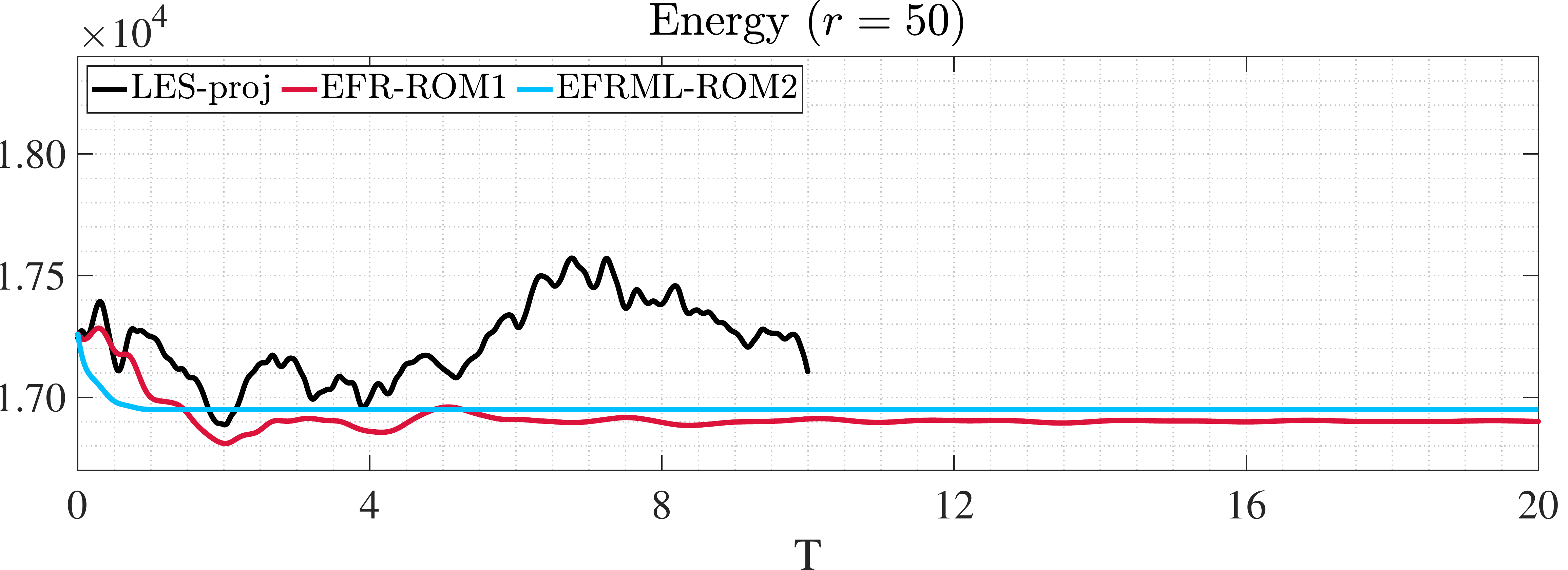}
    \caption{
    Time evolution of the EFR-ROM kinetic energy for $\gamma=9\times 10^{-1}$
    }    
    \label{fig:ke-gamma-2-efr}
\end{figure}

\begin{figure}[H]
\centering
     \begin{subfigure}[b]{0.48\textwidth}
         \centering
    \includegraphics[width=.45\textwidth]{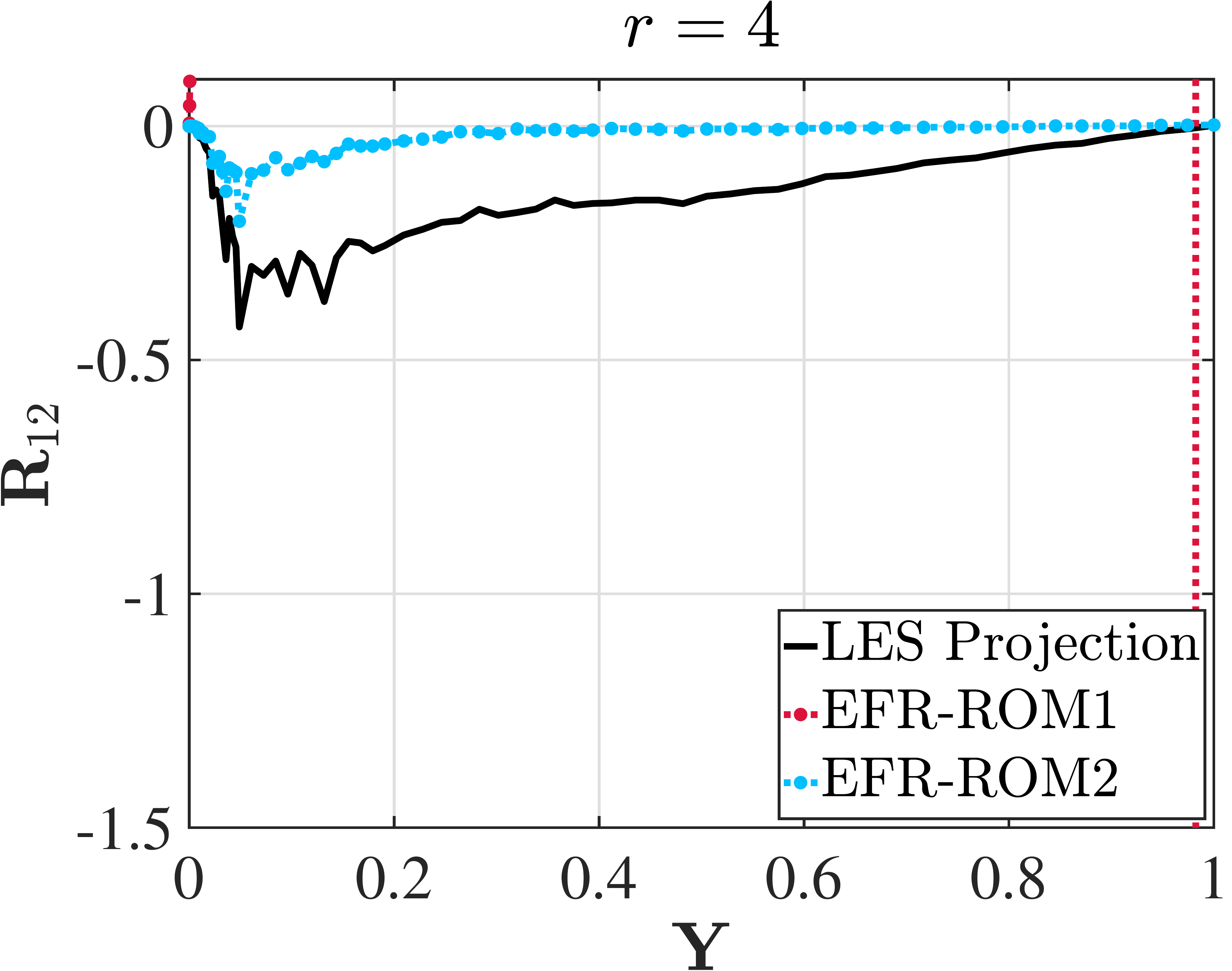}
    \includegraphics[width=.45\textwidth]{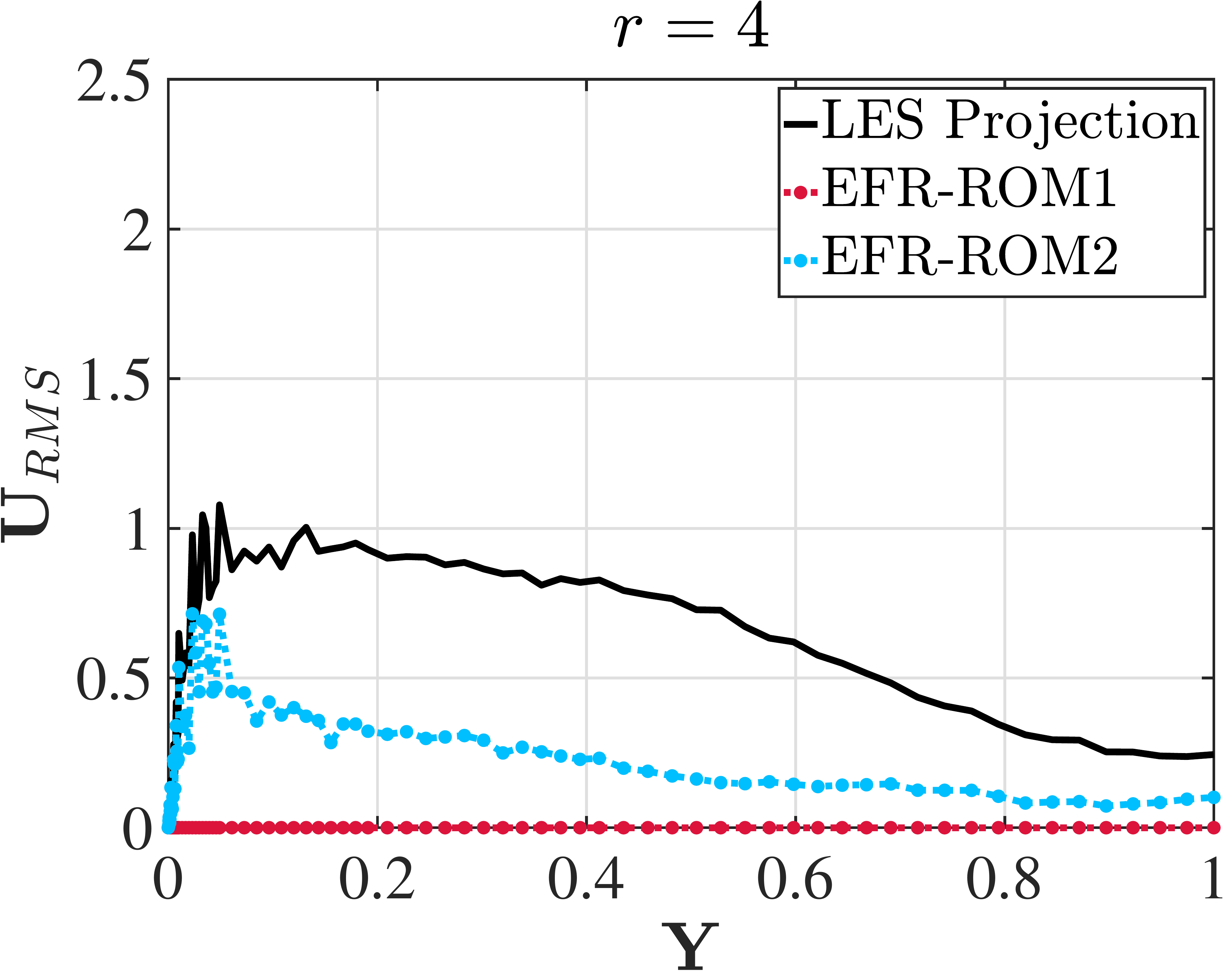}         \caption{$r=4$}
         \label{fig:efr-stat-r-4}
     \end{subfigure}
     \begin{subfigure}[b]{0.48\textwidth}
         \centering
    \includegraphics[width=.45\textwidth]{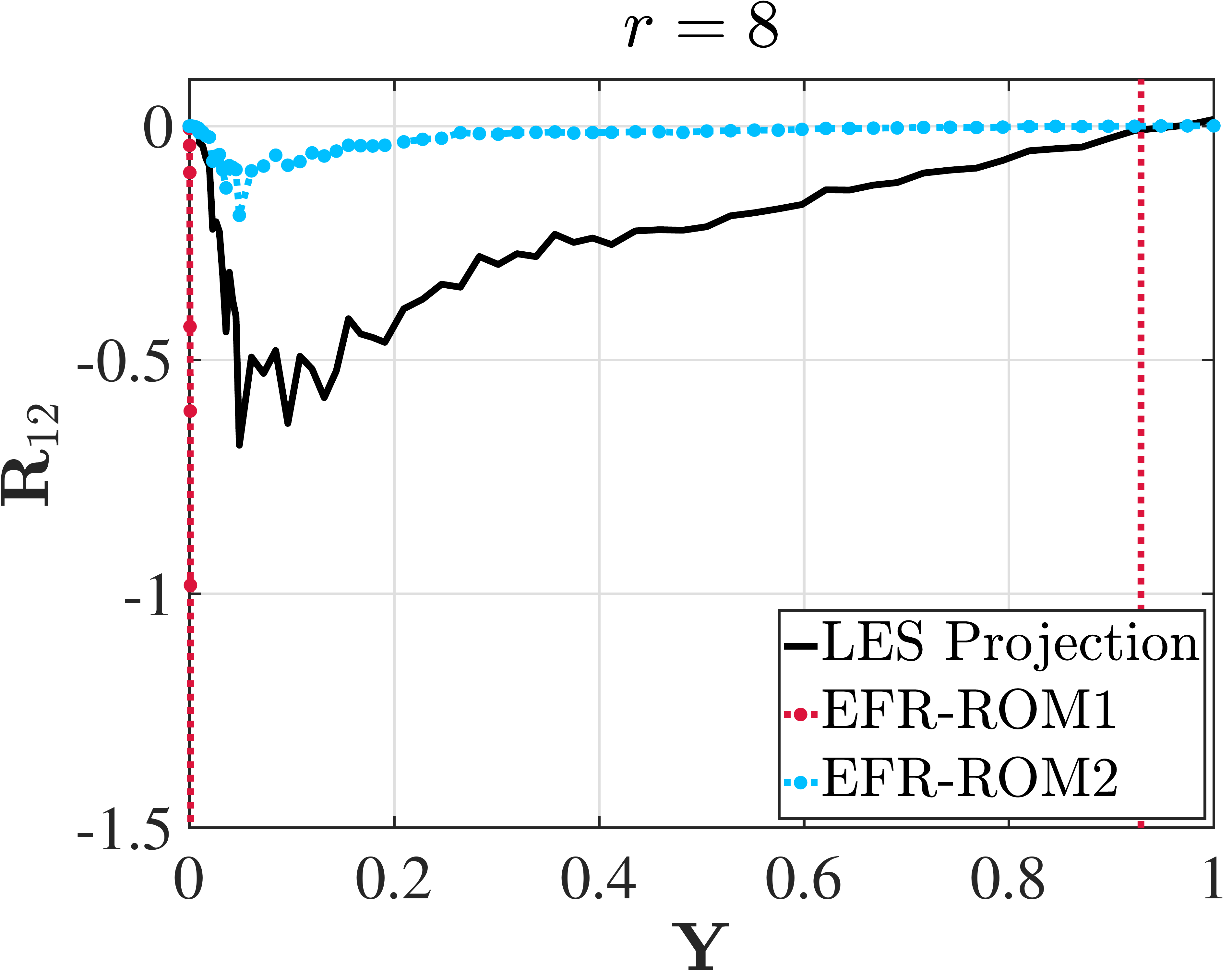}
    \includegraphics[width=.45\textwidth]{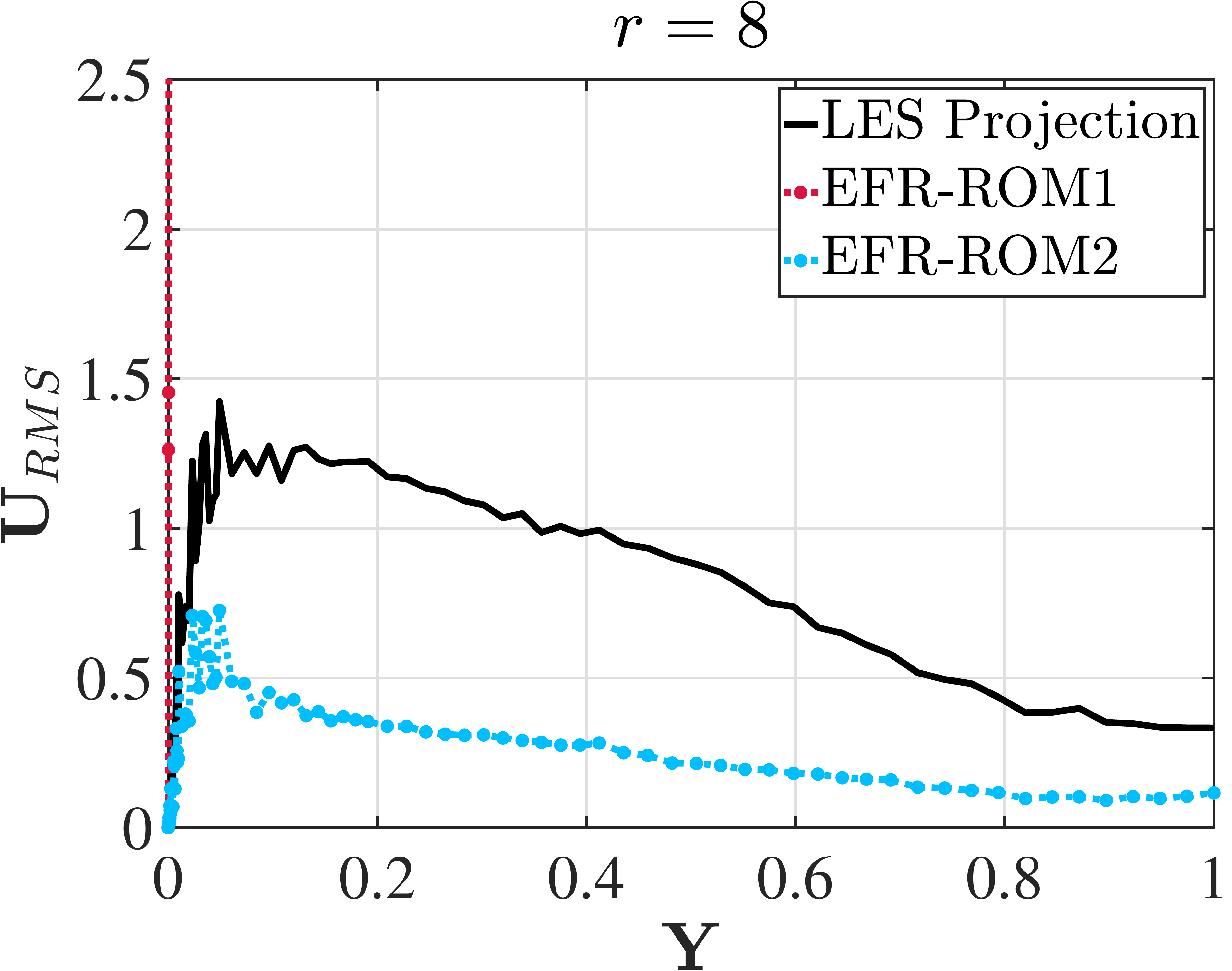}         \caption{$r=8$}
         \label{fig:efr-stat-r-8}
     \end{subfigure}
     \begin{subfigure}[b]{0.48\textwidth}
         \centering
    \includegraphics[width=.45\textwidth]{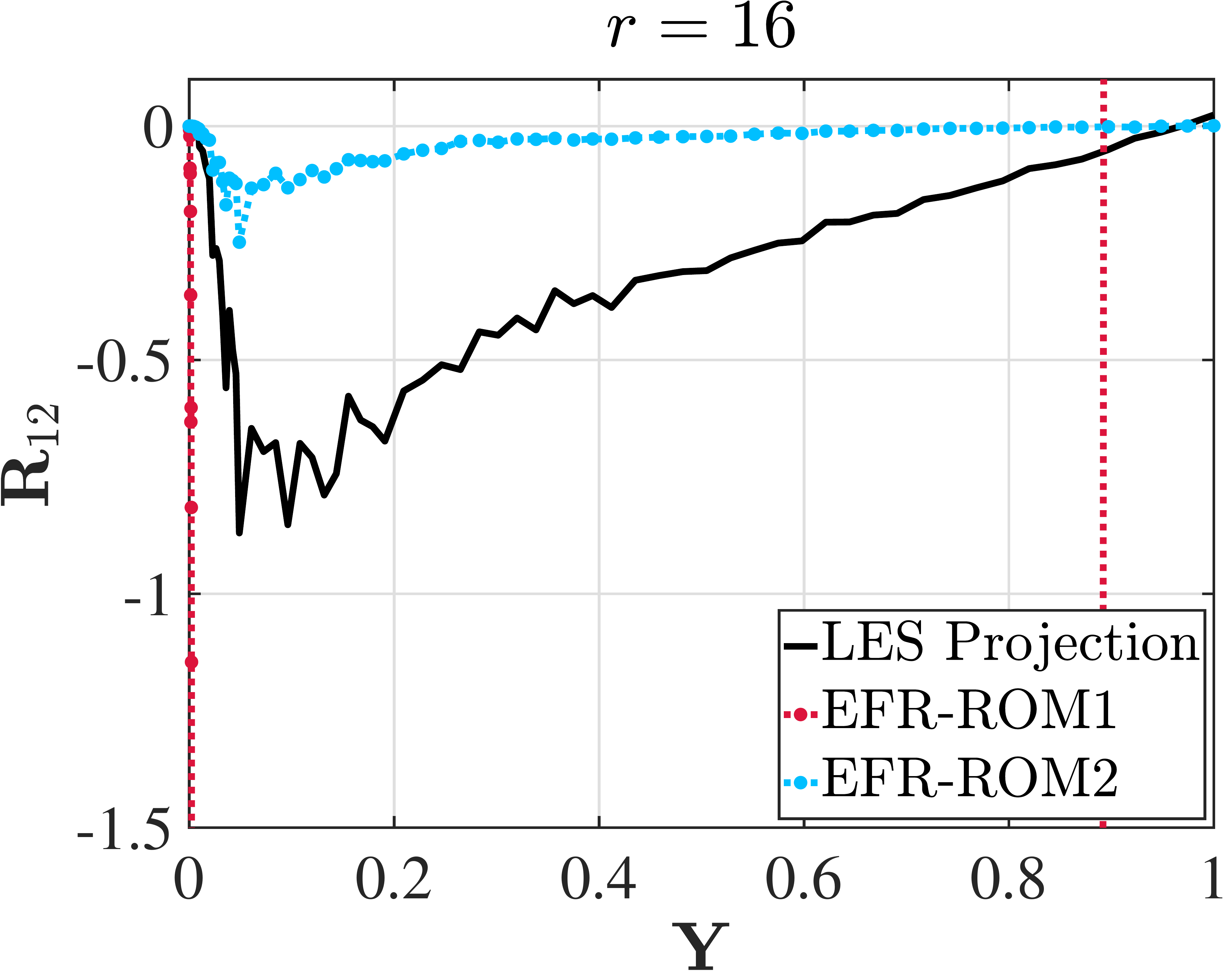}
    \includegraphics[width=.45\textwidth]{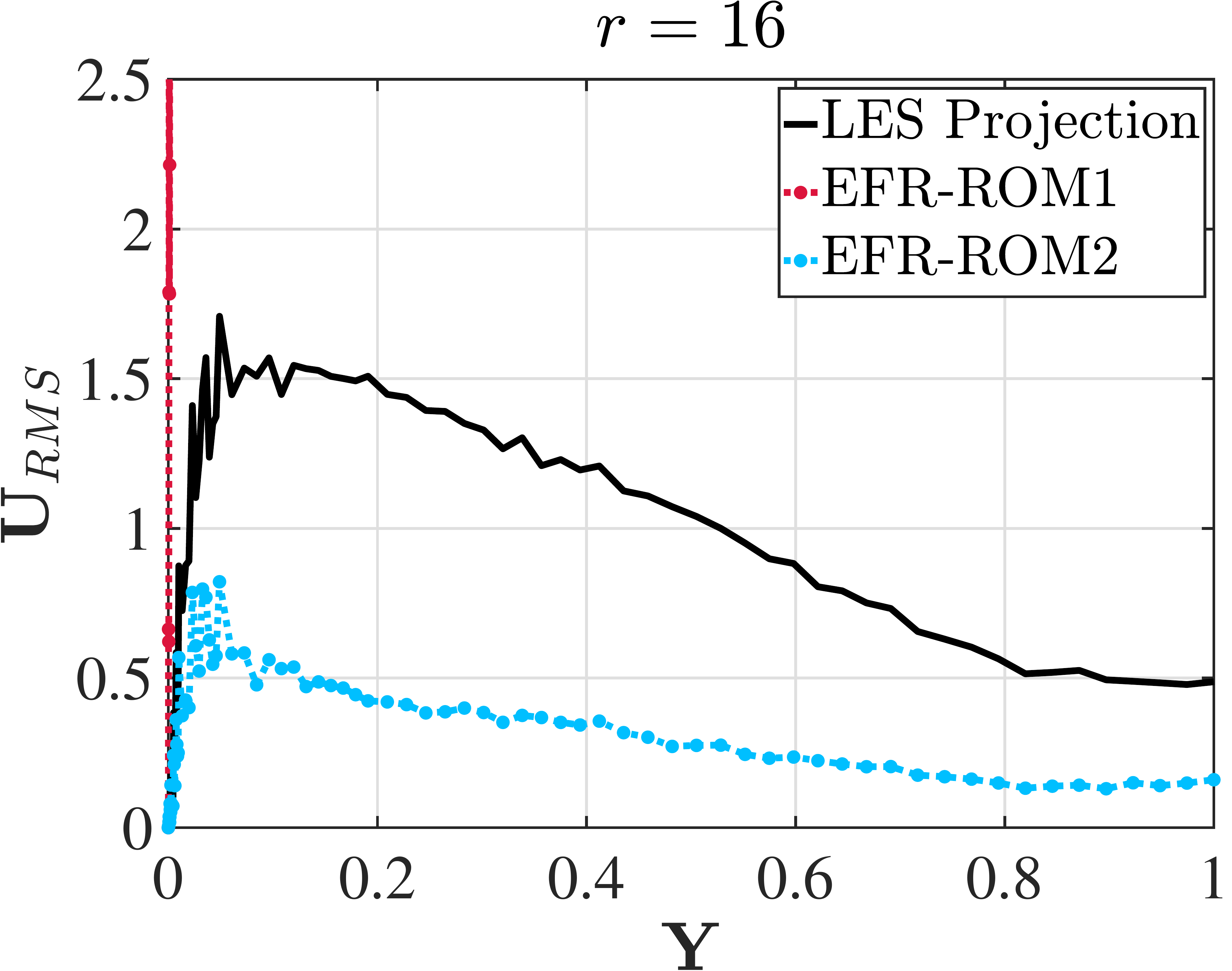}            \caption{$r=16$}
         \label{fig:efr-stat-r-16}
     \end{subfigure}
     \begin{subfigure}[b]{0.48\textwidth}
         \centering
    \includegraphics[width=.45\textwidth]{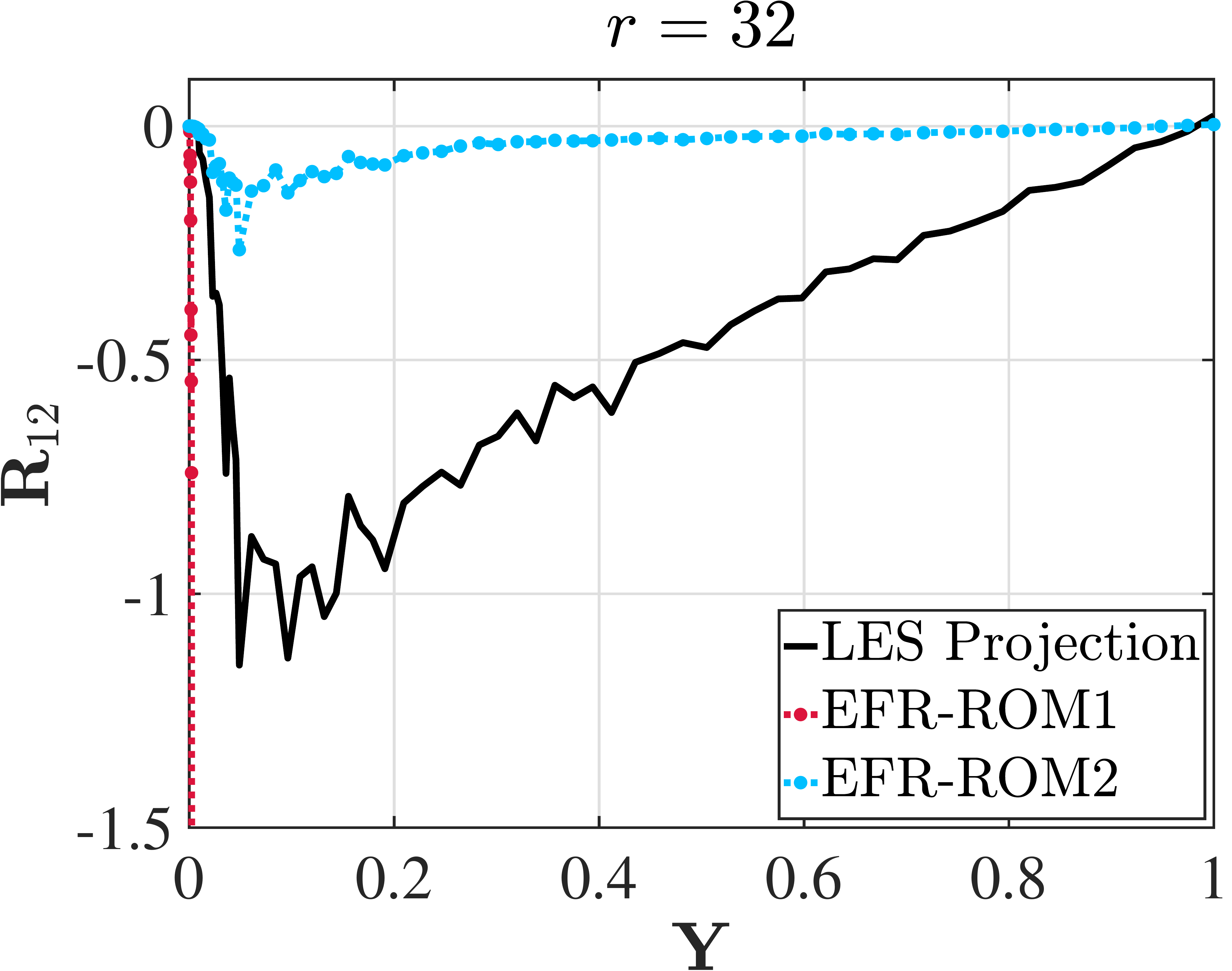}
    \includegraphics[width=.45\textwidth]{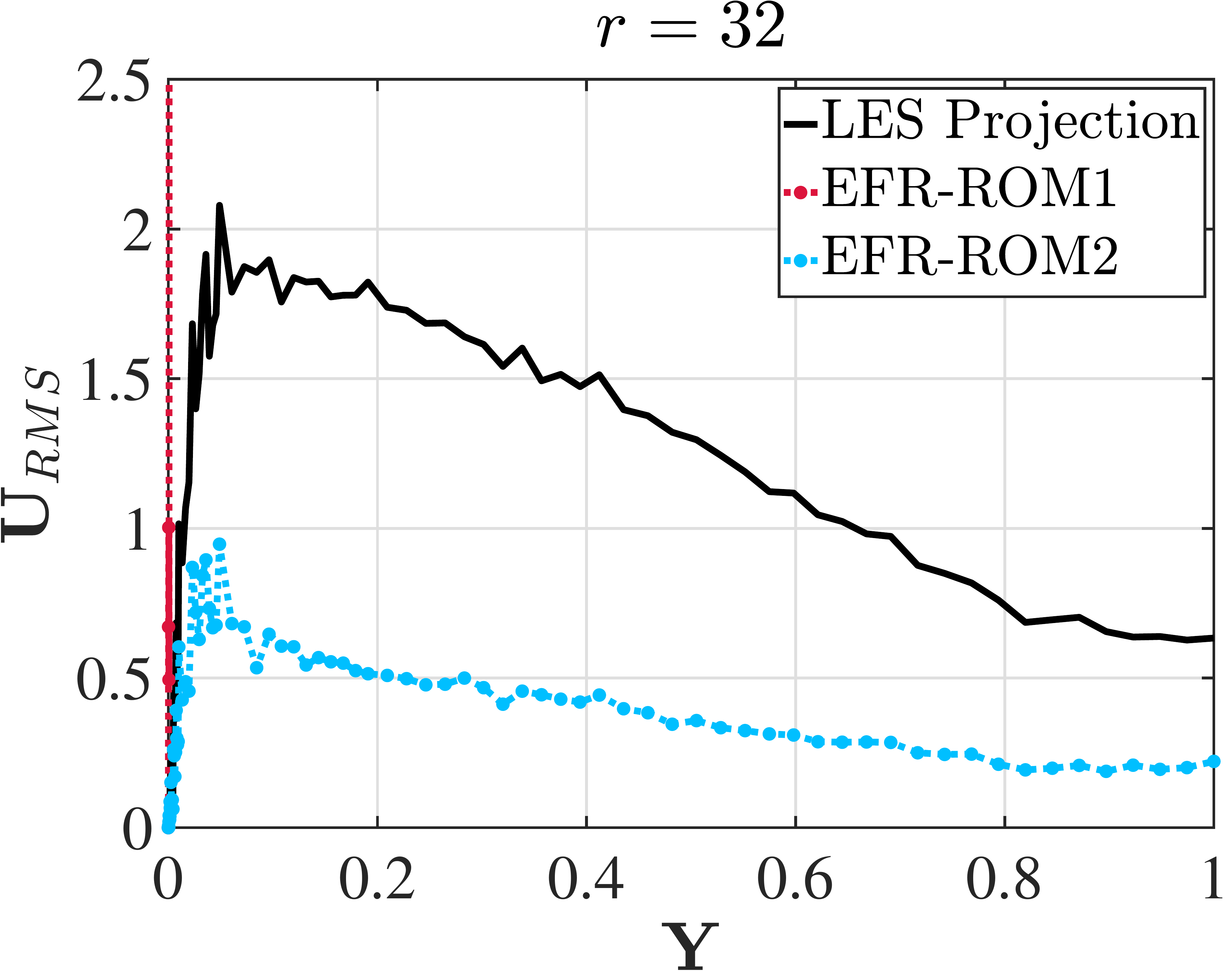}        \caption{$r=32$}
         \label{fig:efr-stat-r-32}
    \end{subfigure}     
     \begin{subfigure}[b]{0.48\textwidth}
         \centering
    \includegraphics[width=.45\textwidth]{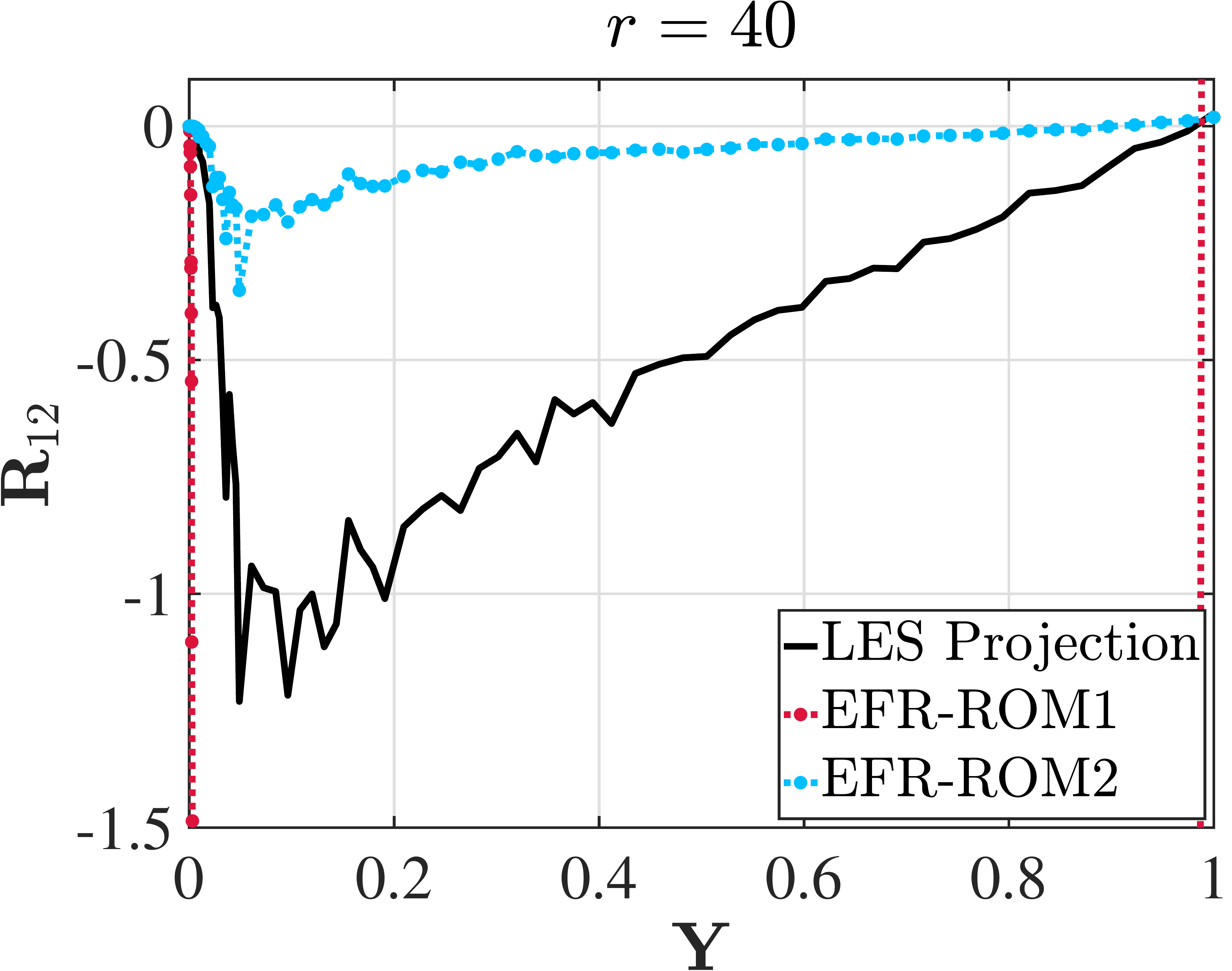}
    \includegraphics[width=.45\textwidth]{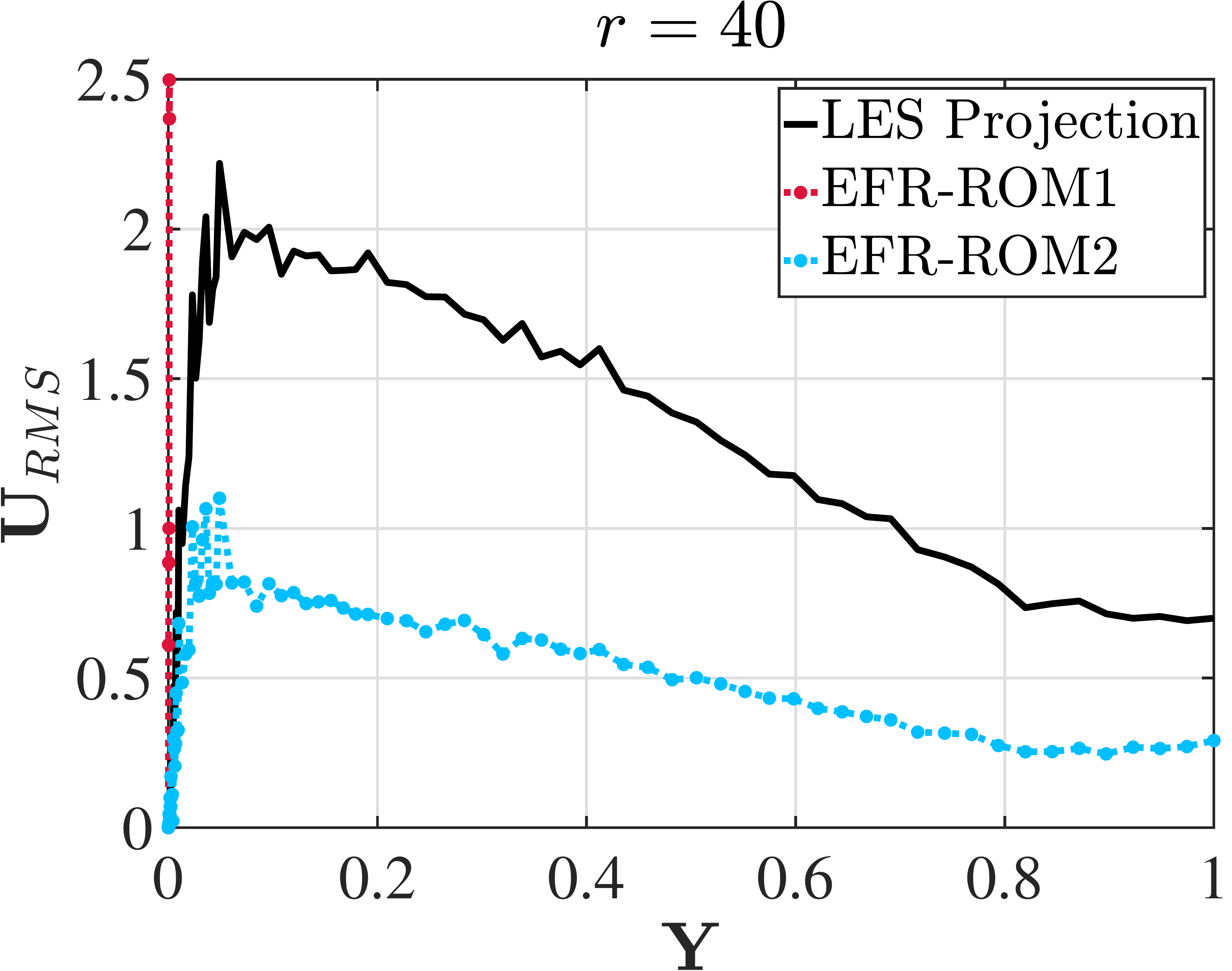}           \caption{$r=40$}
         \label{fig:efr-stat-r-40}
     \end{subfigure} 
     \begin{subfigure}[b]{0.48\textwidth}
         \centering
    \includegraphics[width=.45\textwidth]{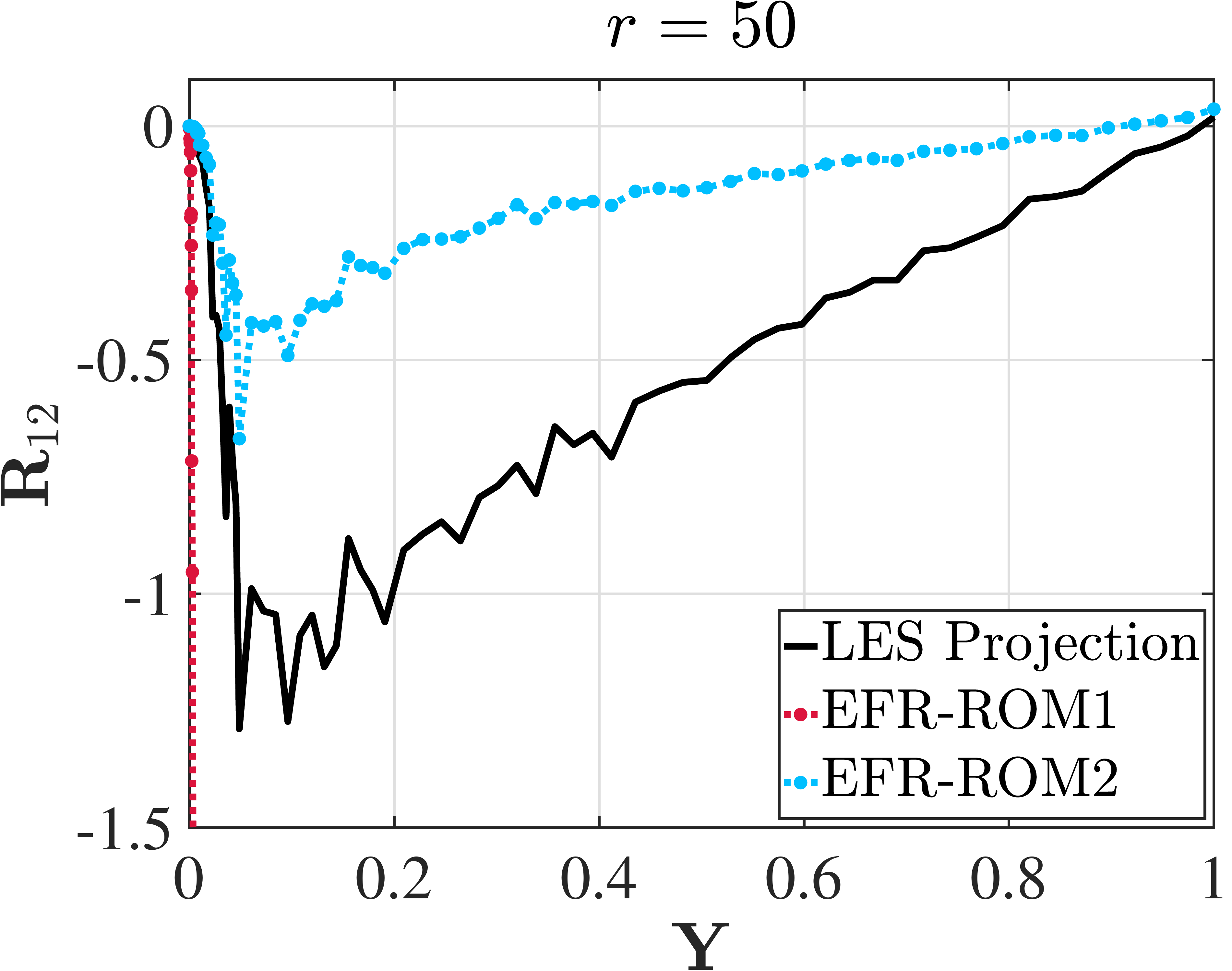}
    \includegraphics[width=.45\textwidth]{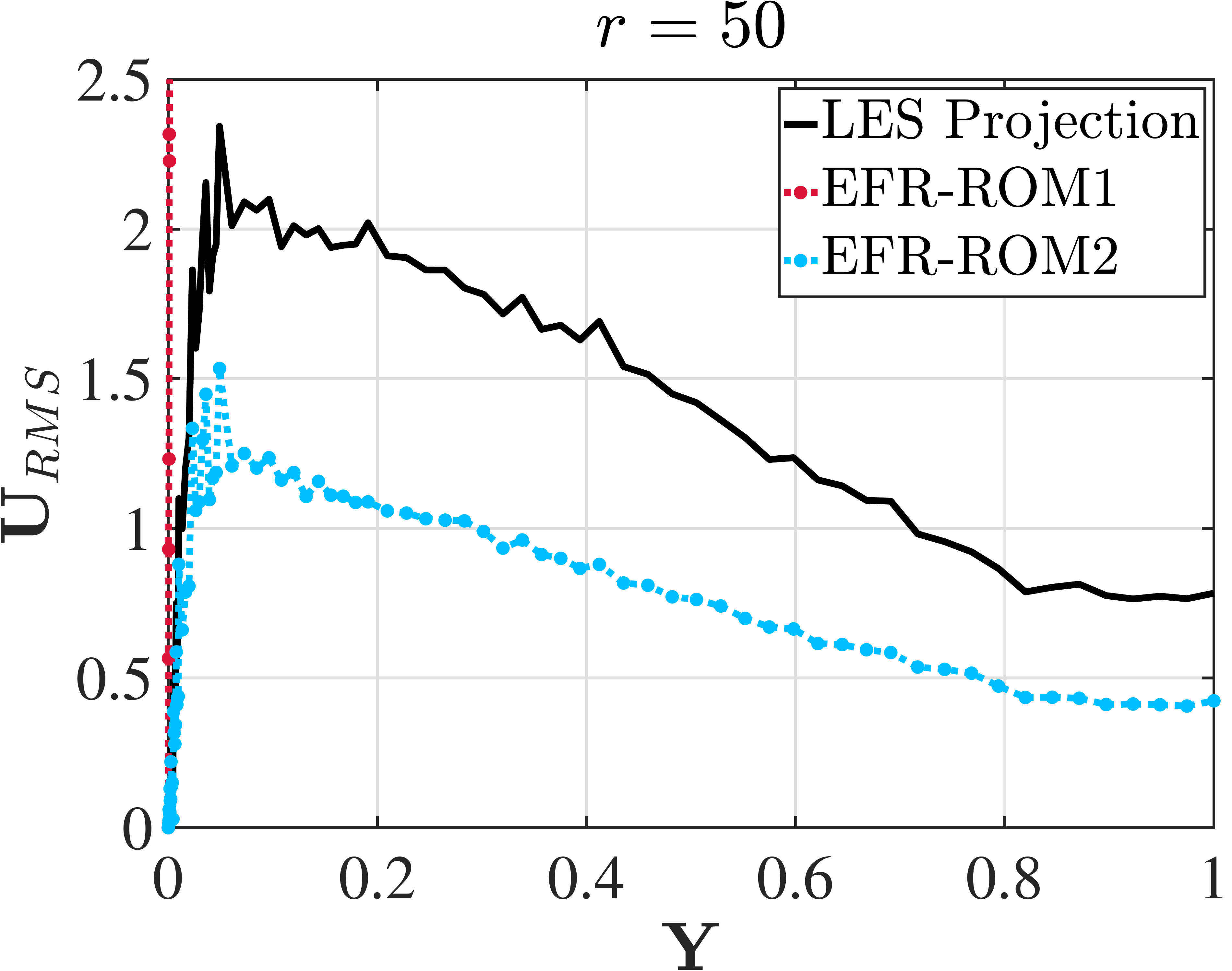}           \caption{$r=50$}
         \label{fig:efr-stat-r-50}
     \end{subfigure} 
     \caption{
    Second-order EFR-ROM statistics for $\gamma=8\times 10^{-2}$
      }   
    \label{fig:stat-gamma-1-efr}
\end{figure}

\begin{figure}[H]
\centering
     \begin{subfigure}[b]{0.48\textwidth}
         \centering
    \includegraphics[width=.45\textwidth]{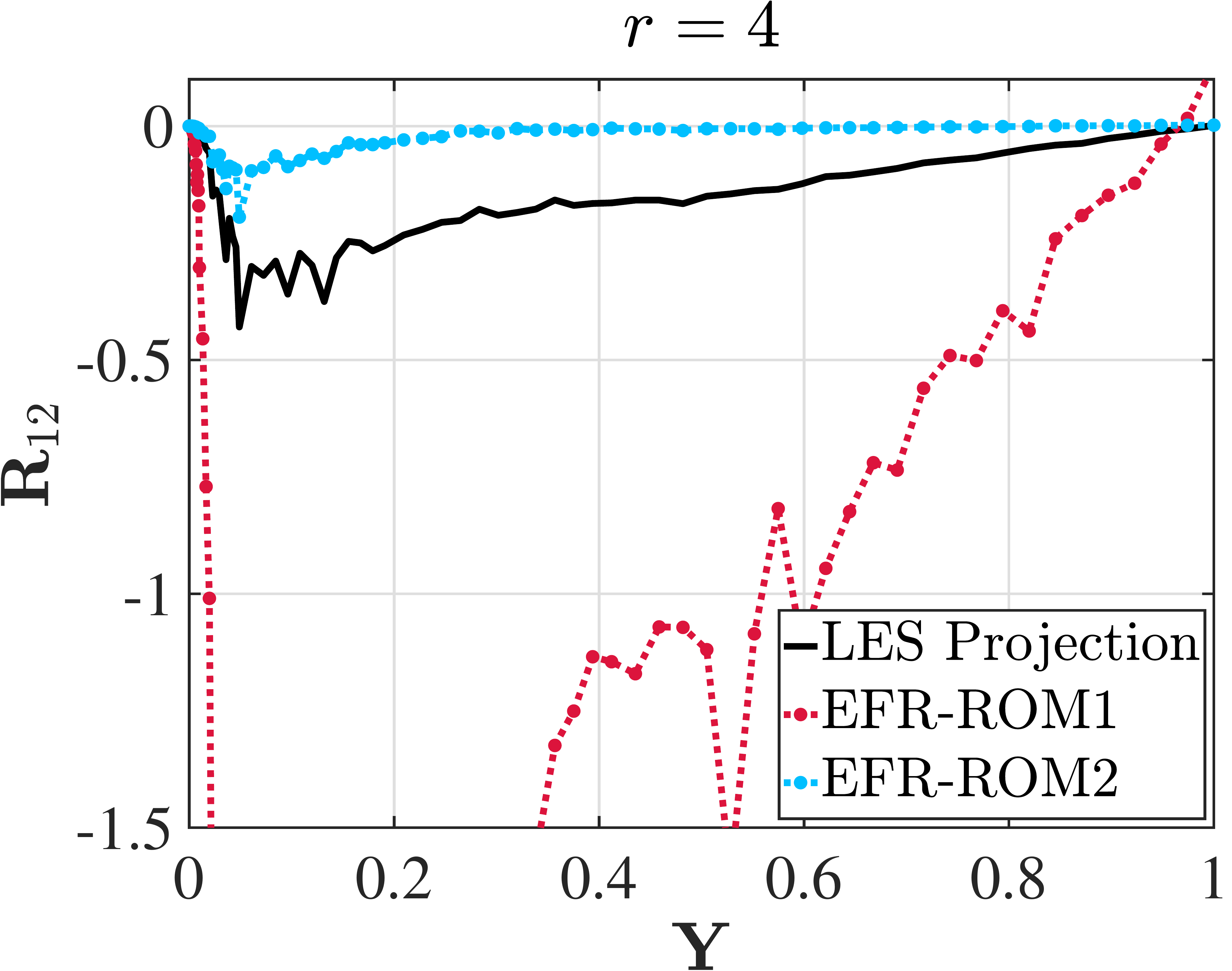}
    \includegraphics[width=.45\textwidth]{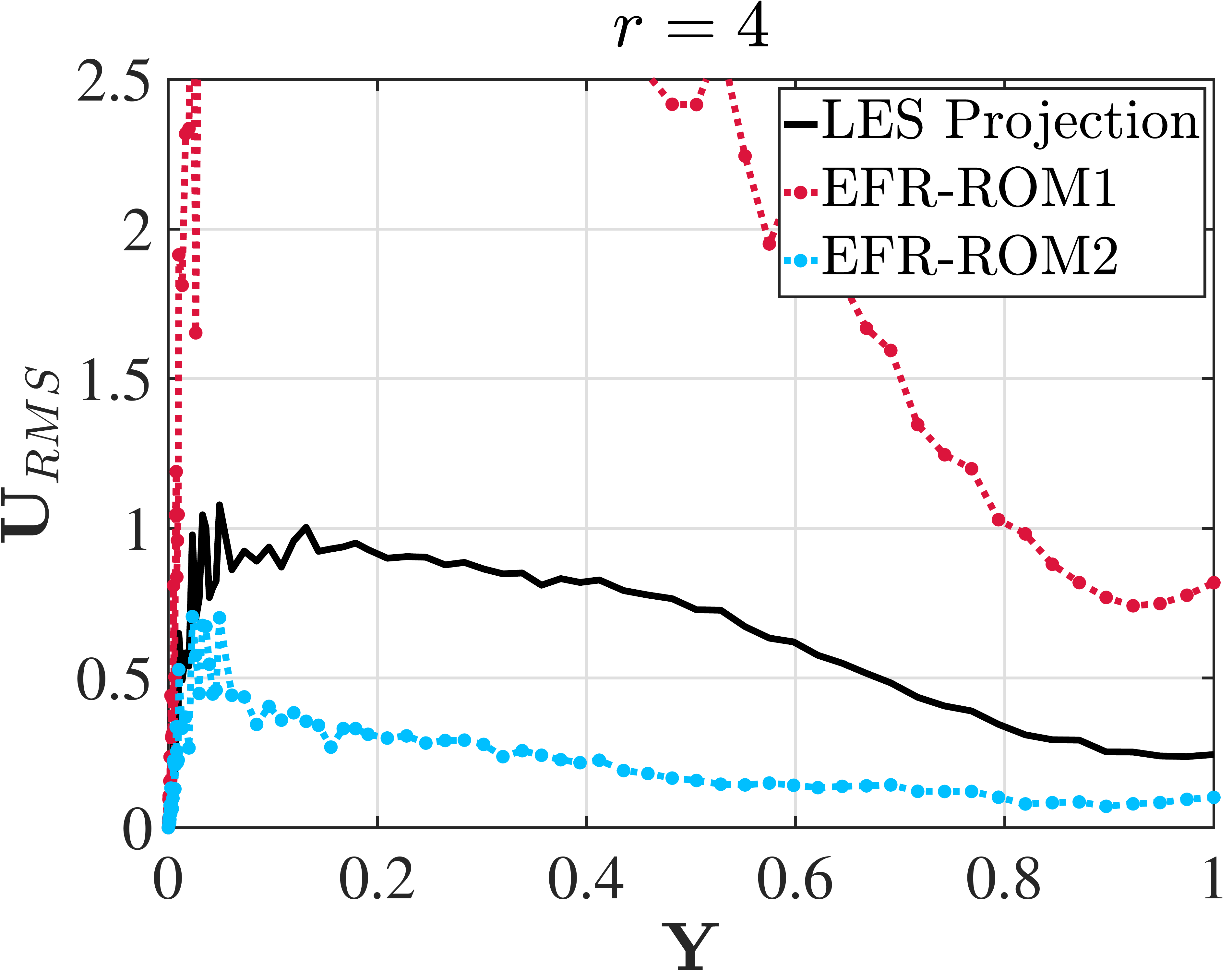}         \caption{$r=4$}
         \label{fig:efr-stat-r-4-gamma2}
     \end{subfigure}
     \begin{subfigure}[b]{0.48\textwidth}
         \centering
    \includegraphics[width=.45\textwidth]{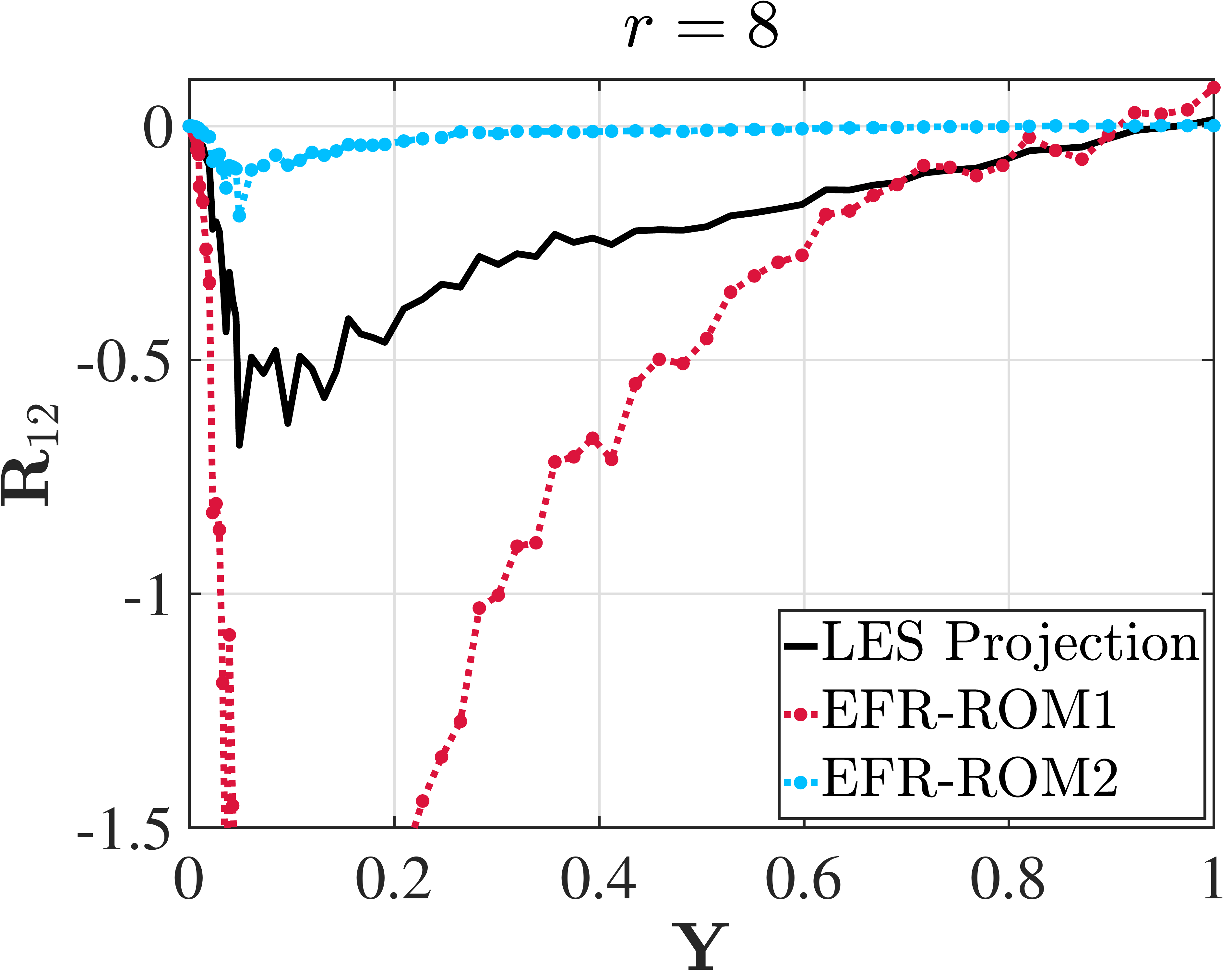}
    \includegraphics[width=.45\textwidth]{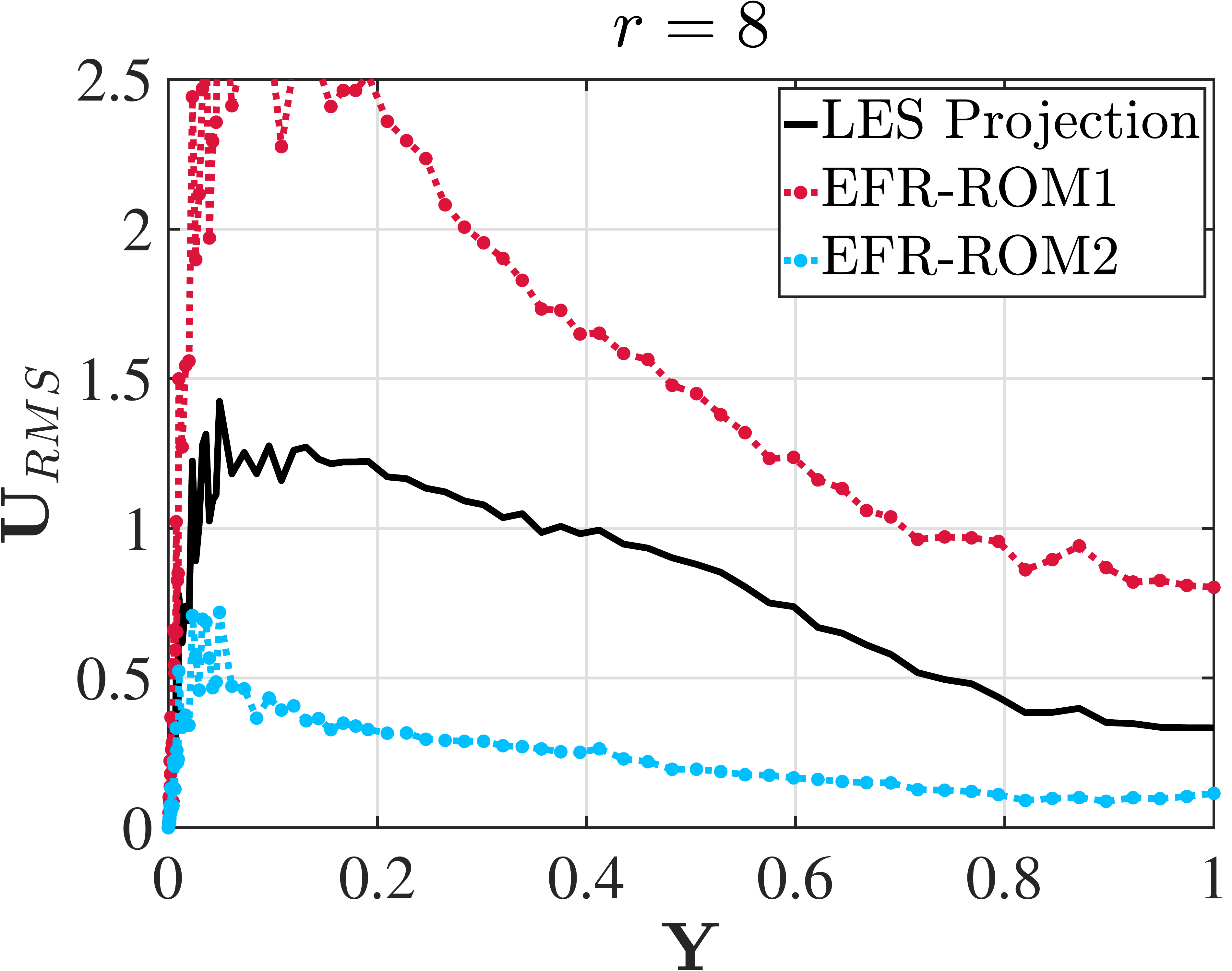}         \caption{$r=8$}
         \label{fig:efr-stat-r-8-gamma2}
     \end{subfigure}
     \begin{subfigure}[b]{0.48\textwidth}
         \centering
    \includegraphics[width=.45\textwidth]{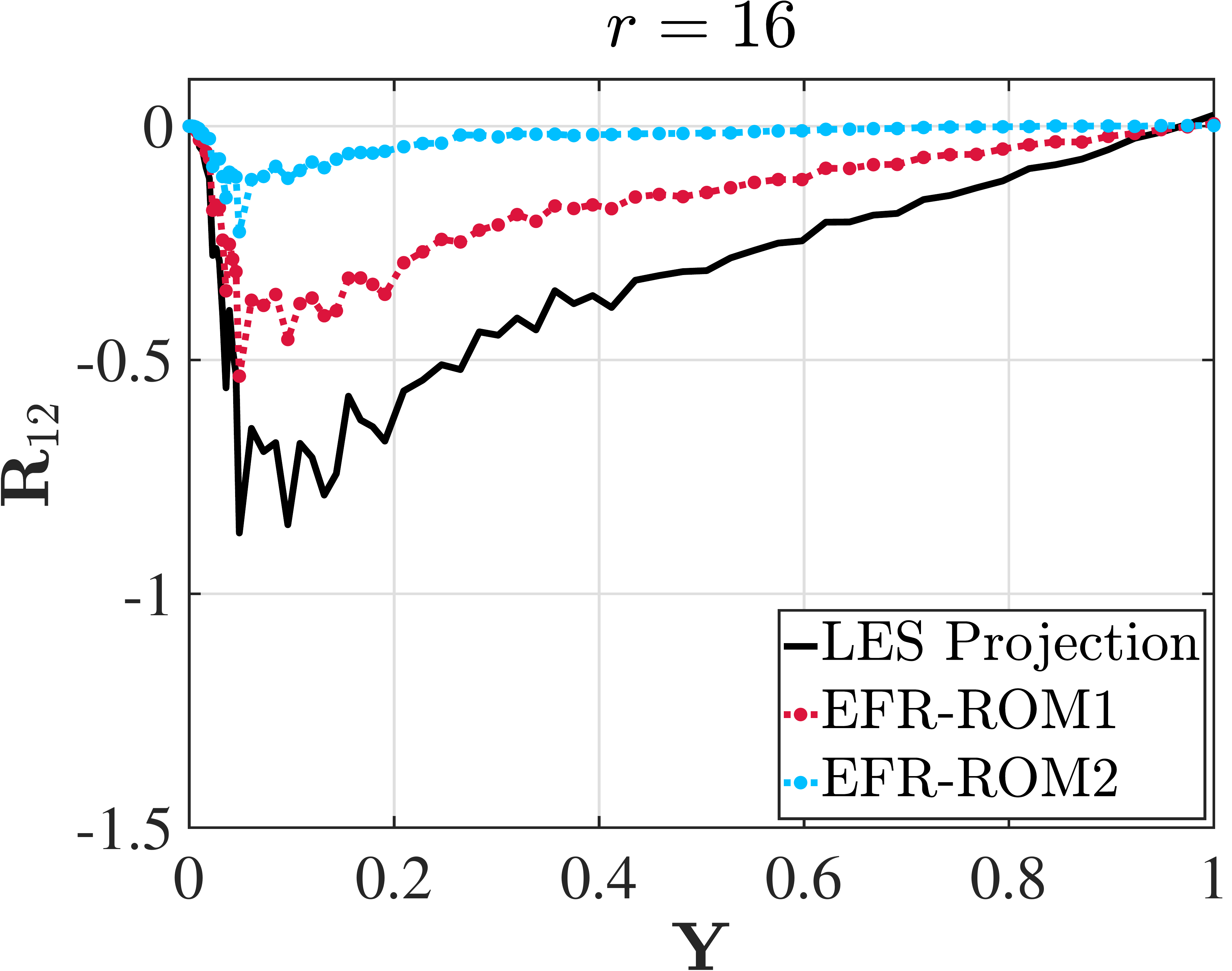}
    \includegraphics[width=.45\textwidth]{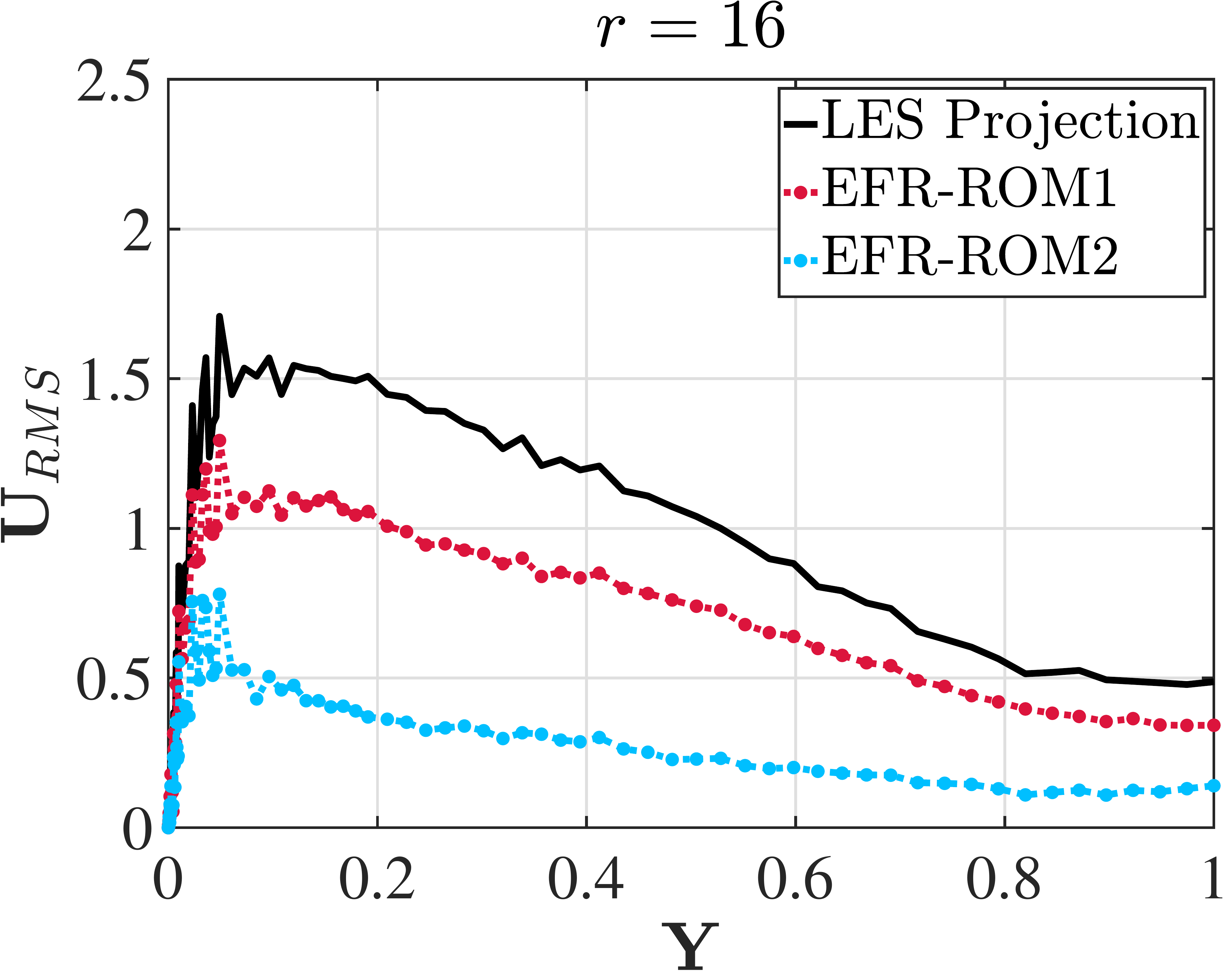}            \caption{$r=16$}
         \label{fig:efr-stat-r-16-gamma2}
     \end{subfigure}
     \begin{subfigure}[b]{0.48\textwidth}
         \centering
    \includegraphics[width=.45\textwidth]{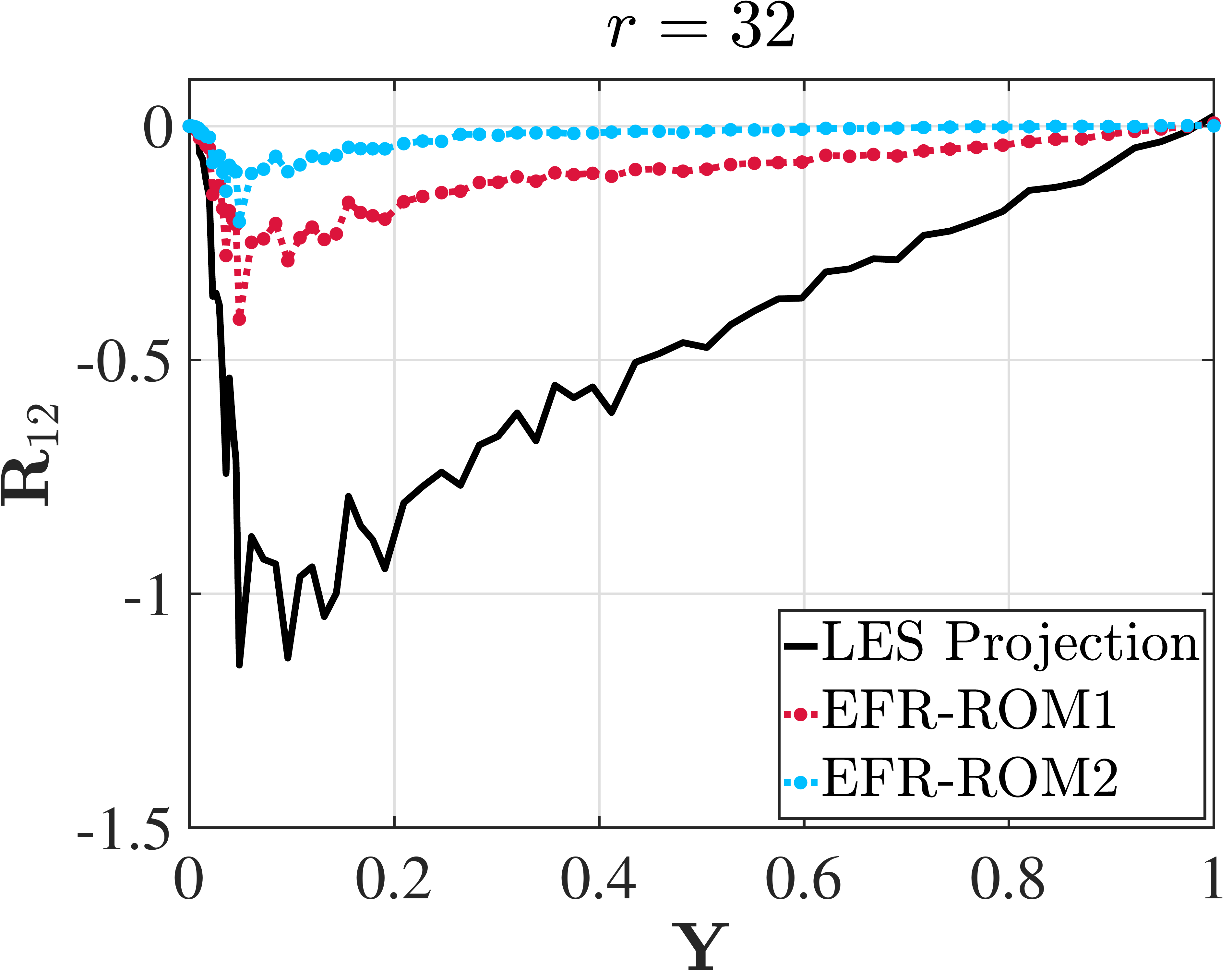}
    \includegraphics[width=.45\textwidth]{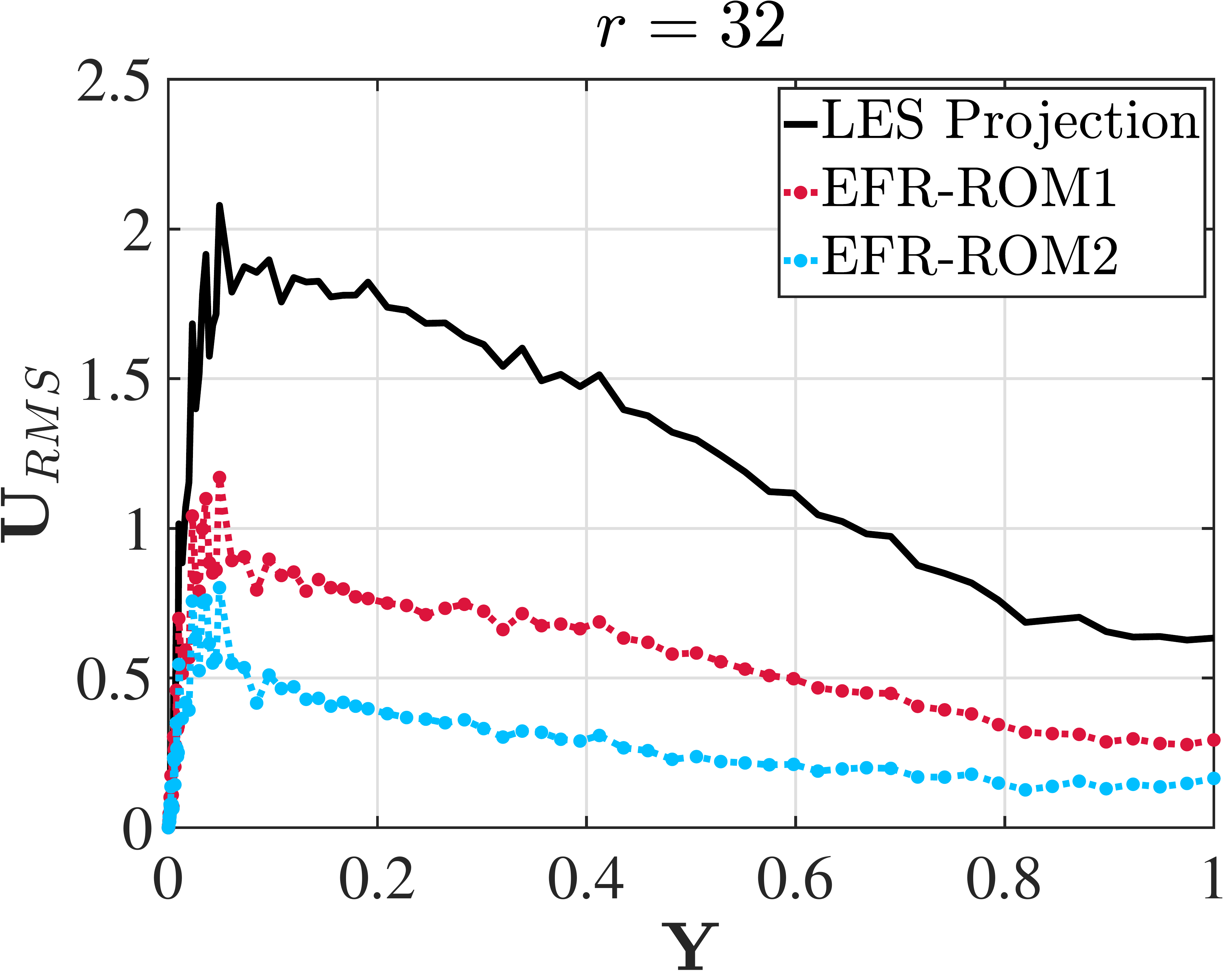}        \caption{$r=32$}
         \label{fig:efr-stat-r-32-gamma2}
    \end{subfigure}     
     \begin{subfigure}[b]{0.48\textwidth}
         \centering
    \includegraphics[width=.45\textwidth]{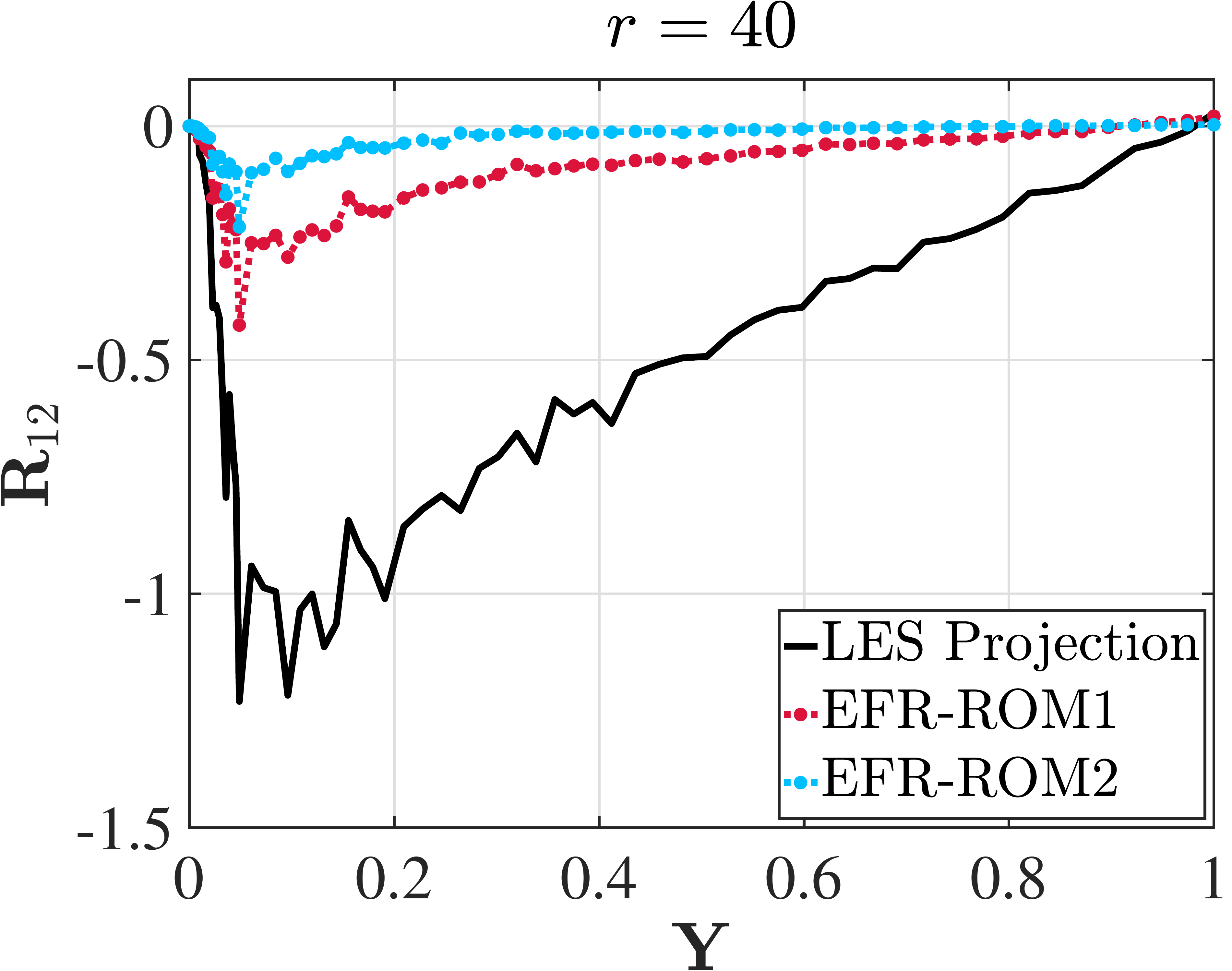}
    \includegraphics[width=.45\textwidth]{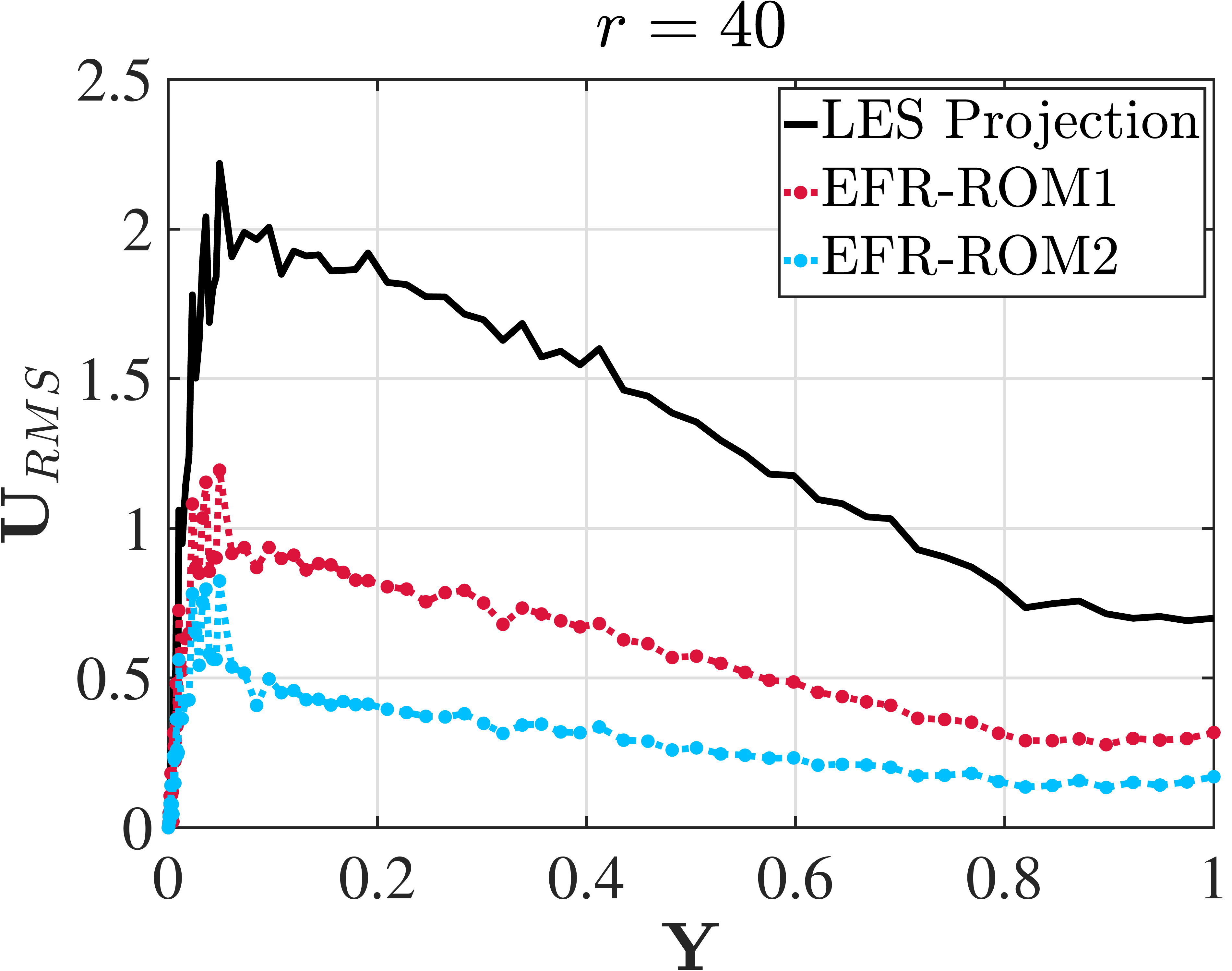}           \caption{$r=40$}
         \label{fig:efr-stat-r-40-gamma2}
     \end{subfigure} 
     \begin{subfigure}[b]{0.48\textwidth}
         \centering
    \includegraphics[width=.45\textwidth]{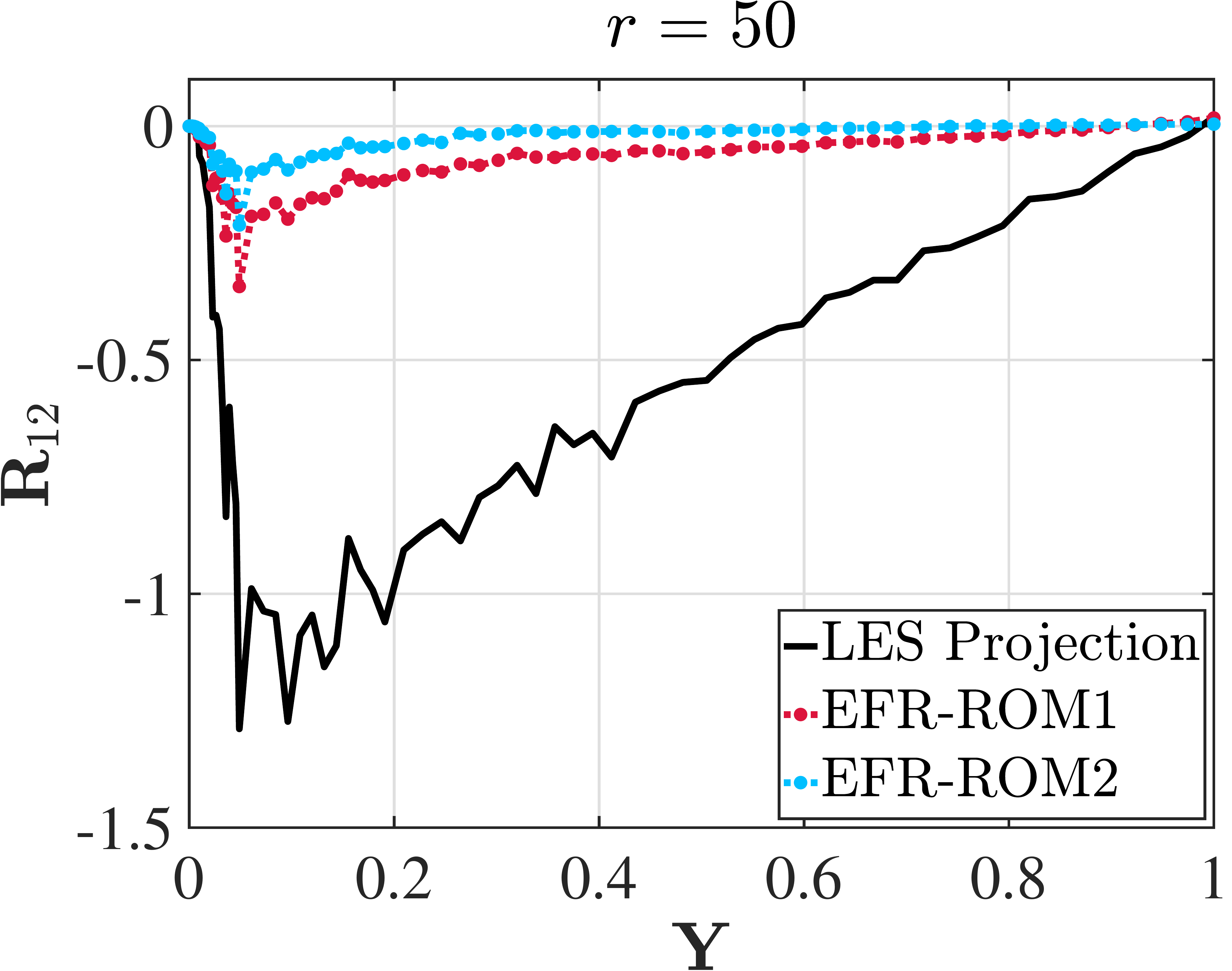}
    \includegraphics[width=.45\textwidth]{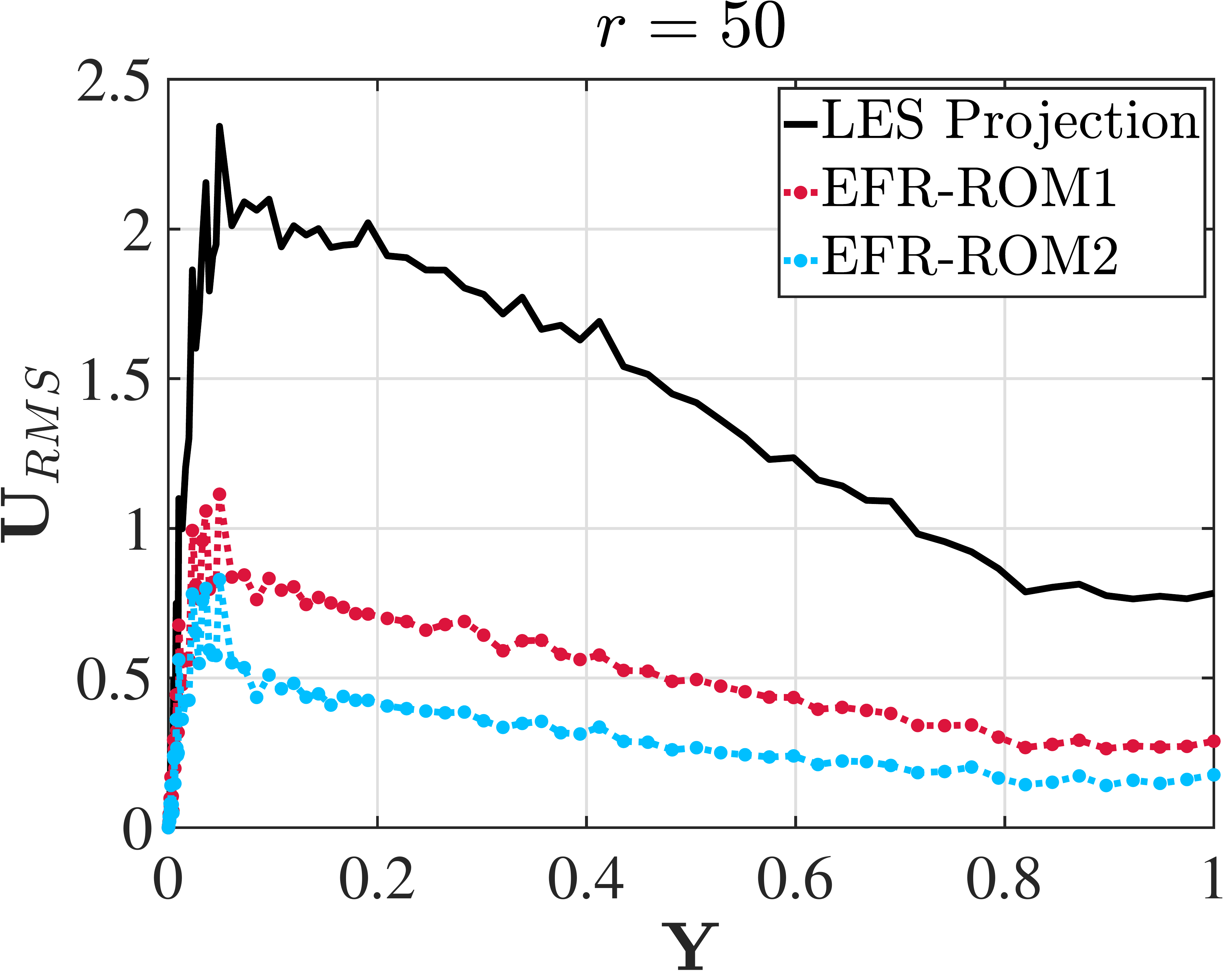}           \caption{$r=50$}
         \label{fig:efr-stat-r-50-gamma2}
     \end{subfigure} 
     \caption{
    Second-order EFR-ROM statistics for $\gamma=9\times 10^{-1}$
      }   
    \label{fig:stat-gamma-2-efr}
\end{figure}

In Figures \ref{fig:ke-gamma-1-efr}--\ref{fig:stat-gamma-2-efr}, we plot the time evolution of the kinetic energy and the second-order statistics of the EFR-ROM1 and EFR-ROM2 for different $r$ values and two different $\gamma$ values:
$\gamma = 8 \times 10^{-2}$  (Figures~\ref{fig:ke-gamma-1-efr} and \ref{fig:stat-gamma-1-efr}) and 
$\gamma = 9 \times 10^{-1}$ (Figures~\ref{fig:ke-gamma-2-efr}and \ref{fig:stat-gamma-2-efr}).
As a benchmark for the ROM results, we use the projection of the the FOM results on the ROM basis (denoted as LES-proj in these plots).

\paragraph{Stability}
Since the results in Table~\ref{table:delta} show that $\delta_2$ is between one and two orders of magnitude higher than $\delta_1$, we expect EFR-ROM2 to yield more stable results than EFR-ROM1.
Indeed, since we fix all the EFR-ROM parameters (i.e., $\chi$ and $\gamma$) and $\delta_2$ is significantly larger than $\delta_1$, we expect the EFR-ROM2 filtering level to be higher than the EFR-ROM1 filtering level (and, thus, EFR-ROM2 to be more stable than EFR-ROM1). 
This is clearly shown in the plots in Figures \ref{fig:ke-gamma-1-efr}--\ref{fig:stat-gamma-2-efr}, in which EFR-ROM2 yields stable results for {\it all} $r$ values and for {\it both} $\gamma$ values.
In contrast, for the smallest $\gamma$ value, $\gamma = 8 \times 10^{-2}$ (Figures~\ref{fig:ke-gamma-1-efr} and \ref{fig:stat-gamma-1-efr}), EFR-ROM1 blows up for {\it all} $r$ values. 
Furthermore, for the largest $\gamma$ value, $\gamma = 9 \times 10^{-1}$ (Figures~\ref{fig:ke-gamma-2-efr} and \ref{fig:stat-gamma-2-efr}), EFR-ROM1 blows up for the small $r$ values (i.e., $r = 4$ and $r = 8$).
To quantify the stability of the two EFR-ROMs, in Table~\ref{tab:gamma-value-threshold}, for different $r$ values, we list the threshold $\gamma_0$ value, i.e., the value that ensures that, if $\gamma > \gamma_0$, then the EFR-ROM is stable.
These results show that, for each $r$ value, the threshold $\gamma_0$ value is more than an order of magnitude lower for EFR-ROM2 than for EFR-ROM1.
Thus, we conclude that EFR-ROM2 is more stable than EFR-ROM1, which is the same conclusion as that yielded by 
Figures \ref{fig:ke-gamma-1-efr}--\ref{fig:stat-gamma-2-efr}.

\begin{table}[H]
    \centering
    \begin{tabular}{c c|c c c c c c c c c c c c }
    \hline\hline
        &$r$ & 4 & 8 & 16 & 32 &40  & 50 
    \\ \hline
        EFR-ROM1  &$\gamma_0$   & 1.0e0&  1.0e0& 9.0e-1& 8.0e-1& 7.0e-1& 7.0e-1
\\
     EFR-ROM2  &$\gamma_0$   & 2.9e-2&  3.3e-2& 3.9e-2& 5.5e-2& 5.9e-2& 7.5e-2
    \\ \hline
    \end{tabular}
        \caption{EFR-ROM threshold $\gamma_0$ values for different $r$ values.}
    \label{tab:gamma-value-threshold}
\end{table}

\paragraph{Accuracy}
Since $\delta_2$ is between one and two orders of magnitude higher than $\delta_1$, we expect the EFR-ROM accuracy to depend on the constant $\gamma$. 
Indeed, 
the EFR-ROM1 and EFR-ROM2 plots in Figures \ref{fig:ke-gamma-1-efr}--\ref{fig:stat-gamma-2-efr} do not display a clear winner:
For the largest $\gamma$ value (i.e., $\gamma = 9 \times 10^{-1}$), EFR-ROM2 is more accurate than EFR-ROM1 for the small $r$ values (i.e., $r=4$ and $8$) since EFR-ROM1 simply blows up.
For the remaining $r$ values, EFR-ROM1 and EFR-ROM2 display the same accuracy level.
For the smallest $\gamma$ value (i.e., $\gamma = 8 \times 10^{-2}$), EFR-ROM2 is more accurate than EFR-ROM1 for {\it all} $r$ values (since EFR-ROM1 simply blows up).
We also note that EFR-ROM2 is relatively accurate for the largest $r$ value (i.e., $r = 50$).

\paragraph{Parameter Sensitivity}
To study the EFR-ROM's parameter sensitivity, we investigate which ROM lengthscale yields EFR-ROMs that are less sensitive (i.e., more robust) with respect to the EFR-ROM's $\gamma$ parameter.
To quantify the ML-ROM's parameter sensitivity, in Table~\ref{tab:gam-value-optimal-2}, for different $r$ values, we list the optimal $\gamma$ value in EFR-ROM, i.e., the $\gamma$ value that ensures that the average ROM kinetic energy, $KE^{ROM}$, is the closest to the everage FOM (LES) kinetic energy, $KE^{LES}$.  
Specifically, we solve the following optimization problem:
\begin{align}
    \min_{\gamma} \left| \overline{KE}^{ROM} - \overline{KE}^{LES} \right|.
    \label{eqn:min-gamma}
\end{align}

\begin{table}[H]
    \centering
    \begin{tabular}{c c|c c c c c c c c c c c c }
    \hline\hline
        &$r$ & 4 & 8 & 16 & 32 &40  & 50 
    \\ \hline
        EFR-ROM1  &$\gamma$ 
        & 1.01e0&   9.98e0&  9.03e-1& 8.05e-1& 8.44e-1& 7.00e-1
    \\ 
        EFR-ROM2   &$\gamma$
        &2.90e-2 & 3.30e-2 &3.90e-2 & 5.50e-2 &7.08e-2 &7.50e-2        
    \\ \hline
    \end{tabular}
    \caption{
    EFR-ROM optimal $\gamma$ values for $\chi = 6 \times 10^{-3}$ and different $r$ values. 
    }
    \label{tab:gam-value-optimal-2}
\end{table}

The results in Table~\ref{tab:gam-value-optimal-2} display a relatively low sensitivity of the optimal EFR-ROM parameter $\gamma$ with respect to changes in $r$.
Indeed, as $r$ varies, the order of magnitude of the optimal $\gamma$ remains the same for both EFR-ROM1 and EFR-ROM2, although the latter is more sensitive than the former.

One possible explanation for the relatively low sensitivity of the optimal EFR-ROM parameter $\gamma$ is that the EFR-ROM parameter $\chi$, which controls the amount of filtering in the relaxation step of the EFR-ROM algorithm, is low ($\chi = 6 \times 10^{-3}$).
Thus, only $0.6 \%$ filtering is applied at each time step of the EFR-ROM algorithm.
Since only a low amount of filtering is used, the effect of the filter radius (i.e., the ROM lengthscales $\delta_1$ and $\delta_2$) is not as important in the EFR-ROM case as in the ML-ROM case.

\begin{table}[H]
    \centering
    \begin{tabular}{c c|c c c c c c c c c c c c }
    \hline\hline
        &$r$ & 4 & 8 & 16 & 32 &40  & 50 
    \\ \hline
        EFR-ROM1  &$\gamma$ 
        & 2.92e-1 &   6.25e-2  & 6.00e-2 &  4.84e-2 &   4.41e-2&   3.61e-2
    \\ 
        EFR-ROM2   &$\gamma$
        & 2.36e-4&   6.83e-5 &  1.12e-04&   2.26e-4 &  3.10e-4 &  4.13e-4       
    \\ \hline
    \end{tabular}
    \caption{
    EFR-ROM optimal $\gamma$ values for $\chi = 6 \times 10^{-2}$ and different $r$ values. 
    }
    \label{tab:gam-value-optimal-3}
\end{table}

To investigate whether a higher percentage of filtering yields a higher sensitivity of the optimal EFR-ROM parameter $\gamma$, we increase the $\chi$ value.  
Specifically, we choose $\chi = 6 \times 10^{-2}$, i.e., $6 \%$ filtering at each EFR-ROM time step  (Table~\ref{tab:gam-value-optimal-3}). 
The results in Table~\ref{tab:gam-value-optimal-3} show that the optimal EFR-ROM $\gamma$ value is very sensitive with respect to changes in $r$.
Indeed, as $r$ increases from $4$ to $50$, $\gamma$ decreases by almost one order of magnitude.
The optimal EFR-ROM2 $\gamma$ values are less sensitive with respect to changes in $r$ than the optimal EFR-ROM1 $\gamma$ values.
Indeed, although the EFR-ROM2 $\gamma$ values vary, their  order of magnitude generally stays constant.

Overall, the results in Tables~\ref{tab:gam-value-optimal-2}--\ref{tab:gam-value-optimal-3} show that the EFR-ROM's parameter $\gamma$ sensitivity is higher for EFR-ROM1 than for EFR-ROM2 when a high percentage of filtering is used in the EFR-ROM algorithm.
As expected, when a low percentage of filtering is used, both EFR-ROM1 and EFR-ROM2 display a relatively low $\gamma$ sensitivity.

Based on the results in Figures \ref{fig:ke-gamma-1-efr}--\ref{fig:stat-gamma-2-efr} and in Tables~\ref{tab:gam-value-optimal-2}--\ref{tab:gam-value-optimal-3}, we conclude that the EFR-ROM investigation in this section yields qualitative results that are similar to those yielded by the ML-ROM investigation in Section~\ref{sec:numerical-results-ml}.

\section{Conclusions}
    \label{sec:conclusions}

In this paper, we proposed a novel ROM lengthscale definition. 
This new ROM lengthscale, denoted $\delta_2$, was constructed by using energy distribution arguments.
Specifically, we balanced the ROM and FOM energy content with the energy content in the $\delta_2$ and $h$ scales, respectively, where $h$ is the FOM mesh size. 
We emphasize that the novel ROM lengthscale, $\delta_2$, is fundamentally different from the current ROM lengthscales, which are built by using dimensional arguments.

We compared the new ROM lengthscale, $\delta_2$, with a standard dimensional based ROM lengthscale, denoted $\delta_1$.
To this end, we used these two ROM lengthscales to build two  mixing-length ROMs (ML-ROMs) and two evolve-filter-relax ROMs (EFR-ROMs) in which all the other parameters were the same.
We investigated the four resulting ML-ROMs and EFR-ROMs in the numerical simulation of the turbulent channel flow at $Re_{\tau} = 395$.

The numerical results 
of our investigation yielded the following conclusions:
\begin{enumerate}
    \item The new energy-based ROM lengthscale, $\delta_2$, was signficantly (two orders of magnitude) larger than the standard ROM lengthscale, $\delta_1$.
    As a result, the ML-ROMs and EFR-ROMs based on the new ROM lengthscale were  significantly more stable than the ML-ROMs and EFR-ROMs based on the standard ROM lengthscale.
    \item 
    The new energy-based ROM lengthscale displayed the correct asymptotic behavior with respect to the ROM dimension, whereas the standard dimensionality-based ROM lengthscale did not.
    \item 
    The ML-ROM parameters based on the new energy-based ROM lengthscale were less sensitive (i.e., more robust) with respect to changes in the ROM dimension ($r$) than the ML-ROMs parameters based on the standard dimensionality-based ROM lengthscale.
    The EFR-ROM parameters based on the new lengthscale were less sensitive with respect to changes in $r$ than the EFR-ROMs parameters based on the standard lengthscale when a significant percentage of filtering was performed.
    For a low filtering percentage, as expected, the EFR-ROM parameters based on the two lengthscales displayed a relatively low sensitivity. 
    In this setting, the standard lengthscale yielded less sensitive parameters than the new lengthscale.
\end{enumerate}

The numerical assessment of the new energy-based ROM lengthscale yielded encouraging results.
We plan to further investigate this  ROM lengthscale in the construction of other types of ROMs, e.g., large eddy simulation ROMs~\cite{wang2012proper,xie2017approximate}. 
We also plan to leverage the new energy based lengthscale to develop scale-aware ROM strategies that are better suited for flow-specific applications.

\section*{Acknowledgments}

The work of the first and fourth authors was supported by NSF through grant DMS-2012253 and CDS\&E-MSS-1953113.
The third author gratefully acknowledges the U.S. DOE Early Career Research Program support through grant DE-SC0019290 and the NSF support through grant DMS-2012255.
Part of this work was funded under the nuclear energy advanced modeling and simulation program.


\bibliography{traian}

\end{document}

%% file: notation.tex
\newcommand{\EX}{{\Bbb{E}}}
\newcommand{\PX}{{\Bbb{P}}}

\newcommand{\lp}{\left(}
\newcommand{\rp}{\right)}
\newcommand{\lb}{\left[}
\newcommand{\rb}{\right]}
\newcommand{\lbr}{\left\{}
\newcommand{\rbr}{\right\}}
\newcommand{\lnorm}{\left\|}
\newcommand{\rnorm}{\right\|}

\newtheorem{remark}{Remark}[section]
\newtheorem{lemma}{Lemma}[section]
\newtheorem{theorem}{Theorem}[section]
\newtheorem{corollary}{Corollary}[section]
\newtheorem{proposition}{Proposition}[section]
\newtheorem{definition}{Definition}[section]
\newtheorem{assumption}{Assumption}[section]

\def\PP{{{\rm l}\kern - .15em {\rm P} }}
\def\PN2{{\PP_{N}-\PP_{N-2}}}

\newcommand{\erf}[1]{\mbox{erf}\left(#1\right)}

\newcommand{\D}{\mathbbm{D}}
\newcommand{\I}{\mathbbm{I}}
\newcommand{\N}{\mathbbm{N}}
\newcommand{\R}{\mathbbm{R}}
\newcommand{\Z}{\mathbbm{Z}}

\newcommand{\cD}{\mathcal{D}}
\newcommand{\cE}{\mathcal{E}}
\newcommand{\cF}{\mathcal{F}}
\newcommand{\cH}{\mathcal{H}}
\newcommand{\cO}{\mathcal{O}}
\newcommand{\cP}{\mathcal{P}}

\newcommand{\bfeta}{\boldsymbol{\eta}}
\newcommand{\bLambdar}{\boldsymbol{\Lambda}_r}
\newcommand{\bLambdarL}{\boldsymbol{\Lambda}_r^{L^2}}
\newcommand{\bLambdarH}{\boldsymbol{\Lambda}_r^{H^1}}
\newcommand{\bmu}{\boldsymbol{\mu}}
\newcommand{\bPhi}{\boldsymbol{\Phi}}
\newcommand{\bPhir}{\boldsymbol{\Phi}_r}
\newcommand{\bphi}{\boldsymbol{\varphi}}
\newcommand{\bphir}{\boldsymbol{\varphi}_r}
\newcommand{\bPsi}{\boldsymbol{\Psi}}
\newcommand{\btau}{\boldsymbol{\tau}}

\newcommand{\ba}{\boldsymbol{a}}
\newcommand{\bas}{{\boldsymbol a}^{snap}}
\newcommand{\bA}{\boldsymbol{A}}
\newcommand{\bb}{\boldsymbol{b}}
\newcommand{\bB}{\boldsymbol{B}}
\newcommand{\bc}{\boldsymbol{c}}
\newcommand{\bd}{\boldsymbol{d}}
\newcommand{\be}{\boldsymbol{e}}
\newcommand{\bff}{\boldsymbol{f}}
\newcommand{\bFF}{{\boldsymbol F}}
\newcommand{\bG}{{\boldsymbol G}}
\newcommand{\bGs}{{\boldsymbol G}^{snap}}
\newcommand{\bh}{\boldsymbol{h}}
\newcommand{\bH}{\boldsymbol{H}}
\newcommand{\bk}{\boldsymbol{k}}
\newcommand{\bL}{\boldsymbol{L}}
\newcommand{\bM}{\boldsymbol{M}}
\newcommand{\bq}{{\boldsymbol q}}
\newcommand{\bqs}{{\boldsymbol q}^{snap}}
\newcommand{\br}{\boldsymbol{r}}
\newcommand{\bS}{\boldsymbol{S}}
\newcommand{\bu}{\boldsymbol{u}}
\newcommand{\bU}{\boldsymbol{U}}
\newcommand{\bur}{{\boldsymbol{u}}_r}
\newcommand{\bUr}{{\boldsymbol{U}}_r}
\newcommand{\bv}{\boldsymbol{v}}
\newcommand{\bV}{\boldsymbol{V}}
\newcommand{\bvr}{{\boldsymbol{v}}_r}
\newcommand{\bw}{\boldsymbol{w}}
\newcommand{\bW}{\boldsymbol{W}}
\newcommand{\bwr}{{\boldsymbol{w}}_r}
\newcommand{\bWr}{{\boldsymbol{W}}_r}
\newcommand{\bx}{\boldsymbol{x}}
\newcommand{\bX}{\boldsymbol{X}}
\newcommand{\bXh}{{\bf X}^h}
\newcommand{\bXr}{{\bf X}^r}
\newcommand{\bY}{\boldsymbol{Y}}
\newcommand{\bz}{\boldsymbol{z}}

\newcommand{\tA}{\tilde{A}}
\newcommand{\tB}{\widetilde{B}}
\newcommand{\tC}{\widetilde{C}}

\newcommand{\oa}{\overline{a}}
\newcommand{\oA}{\overline{A}}
\newcommand{\oc}{\overline{c}}
\newcommand{\obc}{\overline{\boldsymbol c}}
\newcommand{\op}{\overline{p}}
\newcommand{\oU}{\overline{U}}
\newcommand{\obu}{\overline{\boldsymbol u}}
\newcommand{\obU}{\overline{\boldsymbol U}}
\newcommand{\obur}{\overline{\boldsymbol{u}_r}}
\newcommand{\obUr}{\overline{\boldsymbol{U}_r}}
\newcommand{\obv}{\overline{\boldsymbol v}}
\newcommand{\obx}{\overline{\boldsymbol x}}
\newcommand{\ob}[1]{\overline{\boldsymbol{#1}}}
\newcommand{\orr}[1]{\overline{#1}^r}
\newcommand{\obr}[1]{\overline{\boldsymbol{#1}}^r}
\newcommand{\oh}[1]{\overline{#1}^h}
\newcommand{\obh}[1]{\overline{\boldsymbol{#1}}^h}

\newcommand{\as}{a^{snap}}
\newcommand{\CinvNabla}{C_{inv}^{\nabla}(r)}
\newcommand{\CinvDelta}{C_{inv}^{\Delta}(r)}
\newcommand{\Deltar}{\Delta_r}
\newcommand{\Gs}{G^{snap}}
\newcommand{\ct}{\bu_h^{avg}}

\newcommand{\bus}{{\bf u}^*}
\newcommand{\By}{\mathcal B(\by)}
\newcommand{\eci}[1]{\mathcal E_{#1}}
\newcommand{\dpyi}[1]{\delta_{#1}^+(\by)}
\newcommand{\dmyi}[1]{\delta_{#1}^-(\by)}
\newcommand{\cA}{{\mathcal A(\by)}}
\newcommand{\dyi}[1]{\delta_{#1}(\by)}
\newcommand{\cG}{{\mathcal G(\bx,\by)}}
\newcommand{\cGi}[1]{{\mathcal G_{#1}(\bx,\by)}}
\newcommand{\pti}{\partial_i}
\newcommand{\ptii}[1]{\partial_{#1}}
\newcommand{\vertiii}[1]{{|\!|\!| #1 |\!|\!|}}

\newcommand{\half}{\frac{1}{2}}


\newcommand{\red}[1]{{\color{red}#1}}
\newcommand{\blue}[1]{{\color{blue}#1}}
\definecolor{vargreen}{rgb}{0.0, 0.5, 0.0}
\newcommand{\green}[1]{{\color{vargreen} #1}}

\newcommand{\todo}[1]{{\color{red}#1}}
\newcommand{\inserted}[1]{\blue{#1}}
\newcommand{\deleted}[1]{{}}

\newcommand{\TI}[1]{{\color{blue}TI: #1}}
\newcommand{\DRW}[1]{{\color{vargreen}DRW: #1}}
\newcommand{\ZW}[1]{{\color{red}ZW: #1}}
\newcommand{\XX}[1]{{\color{brown}XX: #1}}
\newcommand{\CM}[1]{{\color{purple}CM: #1}}